
\edef\inewcount{\noexpand\csname newcount\endcsname}
\edef\inewdimen{\noexpand\csname newdimen\endcsname}
\edef\inewskip{\noexpand\csname newskip\endcsname}
\edef\inewmuskip{\noexpand\csname newmuskip\endcsname}
\edef\inewbox{\noexpand\csname newbox\endcsname}
\edef\inewhelp{\noexpand\csname newhelp\endcsname}
\edef\inewtoks{\noexpand\csname newtoks\endcsname}
\edef\inewread{\noexpand\csname newread\endcsname}
\edef\inewwrite{\noexpand\csname newwrite\endcsname}
\edef\inewfam{\noexpand\csname newfam\endcsname}
\edef\inewlanguage{\noexpand\csname newlanguage\endcsname}
\edef\inewinsert{\noexpand\csname newinsert\endcsname}
\edef\inewif{\noexpand\csname newif\endcsname}


\countdef\ch=253
\ch="80 \loop\ifnum\ch<"100 \lccode\ch=0 \uccode\ch=0 \advance\ch1 \repeat 
\ch="80 \loop\ifnum\ch<"C0 \catcode\ch=11 \advance\ch1 \repeat 
\ch="C0 \loop\ifnum\ch<"100 \catcode\ch=15 \advance\ch1 \repeat 
\ch="C2 \loop\ifnum\ch<"F5 \catcode\ch=\active \advance\ch1 \repeat 
\ch="C2 \loop\ifnum\ch<"E2 \uccode\ch=1 \advance\ch1 \repeat 
\let\ch\underfined

\catcode0=12 \def\cdef#1#2{\begingroup\lccode0=`#2 \lowercase{\endgroup \def#1{^^@}}} \catcode0=9

\catcode`@=\active
\def@#1{\cdef\chr#1 \edef#1##1{\noexpand\csname\chr##1\noexpand\endcsname}} 
@^^c2@^^c3@^^c4@^^c5@^^c6@^^c7@^^c8@^^c9@^^ca@^^cb@^^cc@^^cd@^^ce@^^cf@^^d0@^^d1@^^d2@^^d3@^^d4@^^d5@^^d6@^^d7@^^d8@^^d9@^^da@^^db@^^dc@^^dd@^^de@^^df
\def@#1{\cdef\chr#1 \edef#1##1##2{\noexpand\csname\chr##1##2\noexpand\endcsname}} 
@^^e0@^^e1@^^e2@^^e3@^^e4@^^e5@^^e6@^^e7@^^e8@^^e9@^^ea@^^eb@^^ec@^^ed@^^ee@^^ef
\def@#1{\cdef\chr#1 \edef#1##1##2##3{\noexpand\csname\chr##1##2##3\noexpand\endcsname}} 
@^^f0@^^f1@^^f2@^^f3@^^f4
\let\chr\undefined
\let\cdef\undefined

\def\grabfuturelet{\futurelet\next\grabexamine}
\def\grabexamine{\ifx\next\csname\expandafter\grab\fi}
\obeylines \def\grab\csname#1\endcsname#2^^M{\expandafter\def\csname#1\endcsname{#2}\expandafter\grabfuturelet} \expandafter\grabfuturelet%
 ~
¢{\hbox{\rm\rlap/c}}
£{\it\$}
«\leftguillemet
­\-
»\rightguillemet
À{\`A}
Á{\'A}
Â{\^A}
Ã{\~A}
Ä{\"A}
Ç{\c C}
È{\`E}
É{\'E}
Ê{\^E}
Ë{\"E}
Ì{\`I}
Í{\'I}
Î{\^I}
Ï{\"I}
Ð\ETH
Ñ{\~N}
Ò{\`O}
Ó{\'O}
Ô{\^O}
Õ{\~O}
Ö{\"O}
Ù{\`U}
Ú{\'U}
Û{\^U}
Ü{\"U}
Ý{\'Y}
Þ\THORN
à{\`a}
á{\'a}
â{\^a}
ã{\~a}
ä{\"a}
ç{\c c}
è{\`e}
é{\'e}
ê{\^e}
ë{\"e}
ì{\`\i}
í{\'\i}
î{\^\i}
ï{\"\i}
ð\eth
ñ{\~n}
ò{\`o}
ó{\'o}
ô{\^o}
õ{\~o}
ö{\"o}
ù{\`u}
ú{\'u}
û{\^u}
ü{\"u}
ý{\'y}
þ\thorn
ÿ{\"y}
Ā{\=A}
ā{\=a}
Ă{\u A}
ă{\u a}
Ć{\'C}
ć{\'c}
Ĉ{\^C}
ĉ{\^c}
Ċ{\.C}
ċ{\.c}
Č{\v C}
č{\v c}
Ď{\v D}
ď{\v d}
Ē{\=E}
ē{\=e}
Ĕ{\u E}
ĕ{\u e}
Ė{\.E}
ė{\.e}
Ě{\v E}
ě{\v e}
Ĝ{\^G}
ĝ{\^g}
Ğ{\u G}
ğ{\u g}
Ġ{\.G}
ġ{\.g}
Ģ{\c G}
ģ{\c g}
Ĥ{\^H}
ĥ{\^h}
Ĩ{\~I}
ĩ{\~\i}
Ī{\=I}
ī{\=\i}
Ĭ{\u I}
ĭ{\u\i}
İ{\.I}
ĲIJ
ĳij
Ĵ{\^J}
ĵ{\^\j}
Ķ{\c K}
ķ{\c k}
Ĺ{\'L}
ĺ{\'l}
Ļ{\c L}
ļ{\c l}
Ľ{\v L}
ľ{\v l}
Ń{\'N}
ń{\'n}
Ņ{\c N}
ņ{\c n}
Ň{\v N}
ň{\v n}
Ō{\=O}
ō{\=o}
Ŏ{\u O}
ŏ{\u o}
Ő{\"O}
ő{\"o}
Ŕ{\'R}
ŕ{\'r}
Ŗ{\c R}
ŗ{\c r}
Ř{\v R}
ř{\v r}
Ś{\'S}
ś{\'s}
Ŝ{\^S}
ŝ{\^s}
Ş{\c S}
ş{\c s}
Š{\v S}
š{\v s}
Ţ{\c T}
ţ{\c t}
Ť{\v T}
ť{\v t}
Ũ{\~U}
ũ{\~u}
Ū{\=U}
ū{\=u}
Ŭ{\u U}
ŭ{\u u}
Ű{\H U}
ű{\H u}
Ŵ{\^W}
ŵ{\^w}
Ŷ{\^Y}
ŷ{\^y}
Ÿ{\"Y}
Ź{\'Z}
ź{\'z}
Ż{\.Z}
ż{\.z}
Ž{\v Z}
ž{\v z}
’'
‘`
”{''}
“{``}
‐-
–{--}
—{---}
¡{!`}
¿{?`}
−-
′'
ß\ss
æ\ae
Æ\AE
œ\oe
Œ\OE
ø\o
Ø\O
å\aa
Å\AA
ł\l
Ł\L
†\dag
‡\ddag
§\S
¶\P
©\copyright
…\dots
ı\relax\ifmmode\imath\else\i\fi
ȷ\relax\ifmmode\jmath\else\j\fi
¬\neg
△\bigtriangleup
÷\div
±\pm
◯\bigcirc
×\times
‖\|
←\leftarrow
→\rightarrow
∥\Vert
⟩\rangle
⟨\langle
∣|
∮\oint
≐\dot=
↑\uparrow
↓\downarrow
α\alpha
β\beta
γ\gamma
δ\delta
ϵ\epsilon
ζ\zeta
η\eta
θ\theta
ι\iota
κ\kappa
λ\lambda
μ\mu
ν\nu
ξ\xi
οo
π\pi
ρ\rho
σ\sigma
τ\tau
υ\upsilon
ϕ\phi
χ\chi
ψ\psi
ω\omega
ε\varepsilon
ϑ\vartheta
ϖ\varpi
ϱ\varrho
ς\varsigma
φ\varphi
Γ\Gamma
Δ\Delta
Θ\Theta
Λ\Lambda
Ξ\Xi
Π\Pi
Σ\Sigma
Υ\Upsilon
Φ\Phi
Ψ\Psi
Ω\Omega
ℵ\aleph
ℏ\hbar
ℓ\ell
℘\wp
ℜ\Re
ℑ\Im
∂\partial
∞\infty
∅\emptyset
∇\nabla
√\surd
⊤\top
⊥\bot
∠\angle
△\triangle
∀\forall
∃\exists
♭\flat
♮\natural
♯\sharp
♣\clubsuit
♢\diamondsuit
♡\heartsuit
♠\spadesuit
∐\coprod
⋁\bigvee
⋀\bigwedge
⨄\biguplus
⋂\bigcap
⋃\bigcup
∫\int
∏\prod
∑\sum
⨂\bigotimes
⨁\bigoplus
⨀\bigodot
⨆\bigsqcup
◁\triangleleft
▷\triangleright
▽\bigtriangledown
∧\wedge
∨\vee
∩\cap
∪\cup
⊓\sqcap
⊔\sqcup
⊎\uplus
⨿\amalg
⋄\diamond
∙\bullet
≀\wr
⊙\odot
⊘\oslash
⊗\otimes
⊖\ominus
⊕\oplus
∓\mp
∘\circ
○\Orb
∖\setminus
⋅\cdot
∗\ast
⨯\times
⋆\star
∝\propto
⊑\sqsubseteq
⊒\sqsupseteq
∥\parallel
∣\divides
⊣\dashv
⊢\vdash
↗\nearrow
↘\searrow
↖\nwarrow
↙\swarrow
⇔\Leftrightarrow
⇐\Leftarrow
⇒\Rightarrow
≠\neq
≤\leq
≥\geq
≻\succ
≺\prec
≈\approx
≽\succeq
≼\preceq
⊃\supset
⊂\subset
⊇\supseteq
⊆\subseteq
∈\in
∋\ni
≫\gg
≪\ll
↔\leftrightarrow
↦\mapsto
∼\sim
≃\simeq
⟂\perp
≡\equiv
≍\asymp
⌣\smile
⌢\frown
↼\leftharpoonup
↽\leftharpoondown
⇀\rightharpoonup
⇁\rightharpoondown
↪\hookrightarrow
↩\hookleftarrow
⋈\bowtie
⊨\models
⟹\Longrightarrow
⟶\longrightarrow
⟵\longleftarrow
⟸\Longleftarrow
⟼\longmapsto
⟷\longleftrightarrow
⟺\Longleftrightarrow
⋯\cdots
⋮\vdots
⋱\ddots
↕\updownarrow
⇑\Uparrow
⇓\Downarrow
⇕\Updownarrow
⌉\rceil
⌈\lceil
⌋\rfloor
⌊\lfloor
≅\cong
∉\notin
⇌\rightleftharpoons
≐\doteq
∄\not\exists
∌\not\ni
∔\dot+
∕/
∤\not|
∦\not\|
∬\int\!\!\!\int
∭\int\!\!\!\int\!\!\!\int
∮\oint
∸\dot-
≁\not\sim
≄\not\simeq
≆\not\cong
≇\not\cong
≉\not\approx
≔:=
≕=:
≢\not\equiv
≭\not\asump
≮\not<
≯\not<
≰\not\le
≱\not\ge
⊀\not\prec
⊁\not\succ
⊄\not\subset
⊅\not\supset
⊈\not\subseteq
⊉\not\supseteq
⊦\vdash
⊧\models
⊬\not\vdash
⊭\not\models
⊲\triangleleft
⊳\triangleright
⋠\not\preceq
⋡\not\succeq
⋤\not\sqsubseteq
⋥\not\sqsupseteq
⋪\not\triangleleft
⋫\not\triangleright
◻\square

\catcode`\^^M=5 %
\let\grabfuturelet\undefined \let\grabexamine\undefined \let\grab\undefined

\let\xcsname=\csname
\let\xendcsname=\endcsname
\def@#1{\def#1##1{\expandafter\ifx\csname\string#1##1\endcsname\relax\errmessage{Undefined UTF-8 sequence \string#1##1}\else\xcsname\string#1##1\xendcsname\fi}}
@^^c2@^^c3@^^c4@^^c5@^^c6@^^c7@^^c8@^^c9@^^ca@^^cb@^^cc@^^cd@^^ce@^^cf@^^d0@^^d1@^^d2@^^d3@^^d4@^^d5@^^d6@^^d7@^^d8@^^d9@^^da@^^db@^^dc@^^dd@^^de@^^df
\def@#1{\def#1##1##2{\expandafter\ifx\csname\string#1##1##2\endcsname\relax\errmessage{Undefined UTF-8 sequence \string#1##1##2}\else\xcsname\string#1##1##2\xendcsname\fi}}
@^^e0@^^e1@^^e2@^^e3@^^e4@^^e5@^^e6@^^e7@^^e8@^^e9@^^ea@^^eb@^^ec@^^ed@^^ee@^^ef
\def@#1{\def#1##1##2##3{\expandafter\ifx\csname\string#1##1##2##3\endcsname\relax\errmessage{Undefined UTF-8 sequence \string#1##1##2##3}\else\xcsname\string#1##1##2##3\xendcsname\fi}}
@^^f0@^^f1@^^f2@^^f3@^^f4
\let@\undefined
\catcode`@=12

\newif\ifscroll 
\newif\ifsuppressunusedbib 

\tracinglostchars=2

\def\printerr#1{\immediate\write17{#1}}
\def\warningline#1#2{\printerr{! #2}\printerr{l.#1}\printerr{}}
\def\ewarningline#1#2#3{\printerr{! #2}\printerr{l.#1 #3}\printerr{}}
\def\warning{\warningline{\the\inputlineno}}

\long\def\gobble#1{}
\ifx\gobbleinit\undefined{\long\gdef\gobbleinit#1\par{}}\fi
\def\expand#1{\edef\expandmacro{#1}\expandmacro\let\expandmacro\undefined}
\def\setetok#1#2{\expand{\noexpand#1{#2}}}
\def\expandtoks#1{\expandafter\edef\expandafter\expandmacro\expandafter{\the#1}#1\expandafter{\expandmacro}}
\def\appendexpand#1#2{\setetok#1{\the#1#2}}
\long\def\append#1#2{#1\expandafter{\the#1#2}}
\long\def\appendtoksexpand#1#2{#1\expandafter\expandafter\expandafter{\expandafter\the\expandafter#1\the#2}}
\long\def\appendonceexpand#1#2{#1\expandafter\expandafter\expandafter{\expandafter\the\expandafter#1#2}}


\def\link#1#2{\lhighlight{#2}}
\def\llink#1{\printlink{llink #1}\link{\ohash#1}}
\catcode`\#=11 \def\ohash{#}\catcode`\#=6
\catcode`\&=11 \def\ampersand{&}\catcode`\&=4
\def\anchor#1#2{\printlink{anchor #1 #2}#2}

\def\setpapersize#1#2{} 
\def\dumpbox#1#2#3{\shipout\vbox{\setpapersize{#1}{#2}\unvbox#3}}
\def\mps#1{\epsfbox{#1}}
\def\metadata#1#2{}
\def\src{} 

\newread\epsffilein
\newif\ifepsfbbfound\inewif\ifepsffilecont
\newdimen\epsfxsize\inewdimen\epsfysize
\newdimen\pspoints\pspoints1bp
\let\runmp\errmessage 
\def\epsfbox#1{\openin\epsffilein=#1 \ifeof\epsffilein\runmp{Could not open file #1}\else
	{\def\do##1{\catcode`##1=12}\dospecials\catcode`\ =10\epsffileconttrue
		\epsfbbfoundfalse
		\loop\read\epsffilein to\epsffileline \ifeof\epsffilein\epsffilecontfalse\else\expandafter\epsfaux\epsffileline :. \\\fi\ifepsffilecont\repeat
		\ifepsfbbfound\else\errmessage{No HiResBoundingBox comment found in file #1}\fi}%
	\closein\epsffilein
	\epsfysize\epsfury\pspoints \advance\epsfysize-\epsflly\pspoints
	\epsfxsize\epsfurx\pspoints \advance\epsfxsize-\epsfllx\pspoints
	\setbox0\hbox{\vbox to\epsfury\pspoints{\vfil\hbox to\epsfxsize{\dimen0=\epsfllx\pspoints \kern-\dimen0 \includegraphics{#1}\hfil}}}%
	\dp0=\epsflly\pspoints \dp0=-\dp0
	\box0 \fi}
\catcode`\%=12 \def\epsfbblit{
\let\dummybrace=} 
\def\epsfaux#1:#2\\{\def\testit{#1}\ifx\testit\epsfbblit \epsfgrab #2 . . . \\\epsffilecontfalse\epsfbbfoundtrue\fi}
\def\empty{}
\def\epsfgrab #1 #2 #3 #4 #5\\{\gdef\epsfllx{#1}\ifx\epsfllx\empty\epsfgrab #2 #3 #4 #5 .\\\else\gdef\epsflly{#2}\gdef\epsfurx{#3}\gdef\epsfury{#4}\fi} 

\newif\ifpdf \pdffalse \ifx\pdfoutput\undefined\else\ifx\pdfoutput\relax\else\ifnum\pdfoutput<1 \else\pdftrue\fi\fi\fi
\ifpdf
\pdfcompresslevel=0 
\pdfobjcompresslevel=0
\def\pdflink#1#2{\leavevmode \lhighlight{\pdfstartlink user { /Subtype /Link /Border [0 0 0] /A << /S #1 >> }#2\pdfendlink}}
\def\llink#1{\pdflink{/GoTo /D (#1)}}
\def\link#1{\pdflink{/URI /URI (#1)}}
\def\anchor#1#2{\pdfdest name {#1} xyz #2}

\def\setpapersize#1#2{\pdfpagewidth#1 \pdfpageheight#2 }
\def\dumpbox#1#2#3{\setpapersize{#1}{#2}\shipout\box#3}
\def\metadata#1#2{\pdfinfo{/Title (#1) /Author (#2)}}
\input supp-pdf.mkii 
\def\mps#1{\convertMPtoPDF{#1}{1}{1}}
\fi

\newtoks\buffertoks
\def\addcode{\immediate\write\mpout}
\def\addunexpandedcode#1{{\toks0={#1}\addcode{\the\toks0}}}
\def\addcodebuffer#1{\edef\tmp{#1}\buffertoks\expandafter\expandafter\expandafter{\expandafter\the\expandafter\buffertoks\tmp}}
\def\addunexpandedcodebuffer#1{\buffertoks\expandafter{\the\buffertoks#1}}

\def\grabcode{\catcode`\#=12 \endlinechar=10 
	\afterassignment\dumpcode\outtoks} 
\def\dumpcode{\addcode{\the\outtoks}\endinlinemp\gobble} 
\def\begininlinemp{\inimp\begingroup\catcode`\^=7 \iftypesetting\mps{\filestem.\the\figno}\let\addcode\gobble\fi \addcode{beginfig(\the\figno);}}
\def\endinlinemp{\addcode{endfig;}\addcode{}\endgroup\global\advance\figno1\relax}
\ifx\endprolog\undefined\let\endprolog\relax\fi

\def\plainfmtname{plain}\ifx\fmtname\plainfmtname\else
\edef\plainoutput{\the\output}
\global\chardef\itfam=4
\def\_{\leavevmode \kern.06em \vbox{\hrule width.3em}} 

\outputpenalty=0
\tracingstats=0
\newlinechar=-1
\maxdeadcycles=25
\showboxbreadth=5
\showboxdepth=3
\errorcontextlines=5
\overfullrule=5pt
\maxdepth=4pt
\parindent=20pt
\abovedisplayskip=12pt plus 3pt minus 9pt
\belowdisplayskip=12pt plus 3pt minus 9pt
\belowdisplayshortskip=7pt plus 3pt minus 4pt


\font\tensy=cmsy10

\catcode"18=12
\catcode`@=11
{
\global\let\end\@@end
\global\let\input\@@input
}
\def\eqalign#1{\null\,\vcenter{\openup\jot\m@th
  \ialign{\strut\hfil$\displaystyle{##}$&$\displaystyle{{}##}$\hfil
      \crcr#1\crcr}}\,}
\catcode`@=12
\fi

\font\tenmsa=msam10 \font\sevenmsa=msam7 \font\fivemsa=msam5 \newfam\msafam \textfont\msafam=\tenmsa \scriptfont\msafam=\sevenmsa \scriptscriptfont\msafam=\fivemsa 
\font\teneufm=eufm10 \font\seveneufm=eufm7 \font\fiveeufm=eufm5 \newfam\eufmfam \textfont\eufmfam=\teneufm \scriptfont\eufmfam=\seveneufm \scriptscriptfont\eufmfam=\fiveeufm 

\font\teneufb=eufb10 \font\seveneufb=eufb7 \font\fiveeufb=eufb5 \newfam\eufbfam \textfont\eufbfam=\teneufb \scriptfont\eufbfam=\seveneufb \scriptscriptfont\eufbfam=\fiveeufb 
\def\frakbf{\fam\eufbfam}
\font\teneurm=eurm10 \font\seveneurm=eurm7 \font\fiveeurm=eurm5 \newfam\eurmfam \textfont\eurmfam=\teneurm \scriptfont\eurmfam=\seveneurm \scriptscriptfont\eurmfam=\fiveeurm 
\def\eurm{\fam\eurmfam}
\font\teneurb=eurb10 \font\seveneurb=eurb7 \font\fiveeurb=eurb5 \newfam\eurbfam \textfont\eurbfam=\teneurb \scriptfont\eurbfam=\seveneurb \scriptscriptfont\eurbfam=\fiveeurb 

\font\teneusm=eusm10 \font\seveneusm=eusm7 \font\fiveeusm=eusm5 \newfam\eusmfam \textfont\eusmfam=\teneusm \scriptfont\eusmfam=\seveneusm \scriptscriptfont\eusmfam=\fiveeusm 

\font\teneusb=eusb10 \font\seveneusb=eusb7 \font\fiveeusb=eusb5 \newfam\eusbfam \textfont\eusbfam=\teneusb \scriptfont\eusbfam=\seveneusb \scriptscriptfont\eusbfam=\fiveeusb 
\def\eucalbf{\fam\eusbfam}

\font\tenss=cmss10 \font\sevenss=cmss7 \font\fivess=cmss5 \newfam\ssfam \textfont\ssfam\tenss \scriptfont\ssfam\sevenss \scriptscriptfont\ssfam\fivess 
\def\sf{\fam\ssfam}

\font\sevenit=cmti7 \scriptfont\itfam=\sevenit

\let\articletitle\seventeenss
\let\chaptertitle\twelvebf
\let\sectiontitle\tenbf
\let\subsectiontitle\tenbfit
\let\subsubsectiontitle\tenit
\let\contchaptertitle\tenbf 
\let\contsectiontitle\tenrm 
\let\contsubsectiontitle\sevenrm 
\let\contsubsubsectiontitle\fiverm 
\let\parnumfont\tenrm 
\let\parbackreffont\fiverm 
\let\proclaimfont\tenbf
\let\prooffont\tenit
\let\mainfont\tenrm


\def\lhighlight{} 
\def\suppresscomments#1#2#3{} 
\def\prohibitcomments#1#2#3{\errmessage{Draft comments are not allowed in the final version}} 

\ifx\format\undefined\else\format\fi 
\ifx\comment\undefined\def\comment{\prohibitcomments}\fi

\newdimen\hmargin
\hmargin=1in
\newdimen\vmargin
\vmargin=1in
\newdimen\plaintextwidth
\plaintextwidth=6.5in
\newdimen\plaintextheight
\plaintextheight=8.9in

\newdimen\totalht
\newdimen\totalwd
\def\shipbox#1{%
	\totalht\ht#1
	\advance\totalht2\vmargin
	\totalwd\wd#1
	\advance\totalwd2\hmargin
	\hoffset-1in
	\advance\hoffset\hmargin
	\voffset-1in
	\advance\voffset\vmargin
	\dumpbox\totalwd\totalht{#1}}

\newtoks\cont
\def\contents{\begingroup\suppressbackreftrue\let\anchor\gobble\the\cont\endgroup}
\let\printcont\gobble
\def\contlinechapt#1#2{\printcont{#2}\smallskip\everypar{}\noindent{\contchaptertitle\llink{chapter.#1}{#2}}\par} 
\def\contlinesect#1#2{\printcont{#1.  #2}\everypar{}\indent{\contsectiontitle\llap{#1.\enskip}\llink{section.#1}{#2}}\par} 
\def\contlineunsect#1#2{\printcont{#2}\everypar{}\noindent{\contsectiontitle\llink{section.#1}{#2}}\par} 
\def\contlinesubsect#1#2{\printcont{#1.  #2}\everypar{}\indent{\contsubsectiontitle{#1.\enskip}\llink{paragraph.#1}{#2}}\par} 
\def\contlinesubsubsect#1#2{\printcont{#1.  #2}\everypar{}\indent\indent\indent{\contsubsubsectiontitle\llap{#1.\enskip}\llink{paragraph.#1}{#2}}\par} 
\def\addcont#1#2#3{\append\cont{#1}\appendexpand\cont{{#2}}\append\cont{{#3}}}

\newif\ifpresec \presecfalse 
\def\chapter#1\par{%
	\def\chapname{#1}%
	\parn0
	\subparn0
	\presectrue
	\numbfalse
	\addcont\contlinechapt\chapname{#1}%
	\curverb{Chapter~}%
	\assignlabel{chapter}{\chapname}%
	\tchapter#1\par}
\def\tchapter#1\par{
	\everypar{\let\beforesect\beforesection}
	\chapbreak\bigbreak
	\centerline{\plabel\chaptertitle\anchor{chapter.#1}{#1}}
	\nobreak\medskip
	\let\beforesect\relax 
}

\let\chapbreak\relax 

\newcount\secn \secn0
\def\section#1\par{%
	\ifx\sectionid\undefined\advance\secn1 \edef\sectionid{\the\secn}\fi
	\everypar{\numpar}
	\parn0
	\subparn0
	\presecfalse
	\numbfalse
	\addcont\contlinesect\sectionid{#1}%
	\curverb{\S}%
	\assignlabel{section}{\sectionid}%
	\tsection#1\par
}
\def\unsection#1\par{%
	\everypar{\numpar}
	\parn0
	\subparn0
	\presecfalse
	\numbfalse
	\addcont\contlineunsect{#1}{#1}%
	\curverb{\S}%
	\assignlabel{section}{#1}%
	\tsection#1\par
}
\def\tsection#1\par{
	\beforesect\let\beforesect\beforesection
	\typesetsection{#1}%
	\aftersection	
	\let\sectionid\undefined
}
\def\beforesection{\vskip0pt plus.3\vsize \penalty-250 \vskip0pt plus-.3\vsize \bigskip \vskip\parskip}
\def\aftersection{\nobreak\smallskip}
\def\typesetsection#1{\leftline{\sectiontitle
	\ifx\sectionid\undefined
		\plabel\anchor{section.#1}{}%
	\else
		\hbox to \parindent{\hss\plabel\anchor{section.\sectionid}{\sectionid}\enspace\hfill}%
	\fi#1}} 
\let\beforesect\beforesection

\def\subsection#1\par{\bigbreak\numbtrue\curverb{\S}\subsectiontitle\noindent#1\/\mainfont\par\nobreak\medskip\addcont\contlinesubsect{\the\secn.\the\parn}{#1}}
\def\subsubsection#1\par{\bigbreak\numbtrue\curverb{\S}\subsubsectiontitle\noindent#1\/\mainfont\par\nobreak\medskip\addcont\contlinesubsubsect{\the\secn.\the\parn}{#1}}

\def\sskip#1{\ifdim\lastskip<\medskipamount \removelastskip\penalty#1\medskip\fi}
\def\slug{\hbox{\kern1.5pt\vrule width2.5pt height6pt depth1.5pt\kern1.5pt}}
\newif\ifqed \newif\ifneedqed
\def\qed{\unskip\nobreak\ \slug\ifhmode\spacefactor3000 \fi\global\qedtrue}
\def\proclaim{\medbreak\atendpar{\sskip{55}}\numbtrue\gproclaim\proclaimfont} 
\def\proof{\medbreak\atendpar{\sskip{-55}}\needqedtrue\gproclaim\prooffont}
\def\xproclaim#1.{\medbreak\atendpar{\sskip{55}}{\everypar{}\noindent}{\proclaimfont#1.\enspace}\ignorespaces} 
\let\abstract\xproclaim 

\def\ppar{\endgraf{\everypar{}\indent}} 
\newtoks\atendpar 
\newtoks\atendbr 
\newif\ifnumb \numbfalse 
\def\finishpar{\ifhmode\ifneedqed\ifqed\else\qed\fi\qedfalse\needqedfalse\fi\iflist\endlist\fi\the\atendbr\atendbr{}\endgraf\the\atendpar\atendpar{}\numbfalse\fi}
\def\endlist{\iflist\listfalse\endgraf{\parskip\smallskipamount\everypar{}\noindent}\fi} 
\newcount\parn
\newcount\subparn
\newif\ifparbref
\def\nextpar{\ifnumb\advance\parn1 \printlabel{advancing paragraph number to \the\parn}\def\brt{}%
	\ifparbref\edef\cseq{\csname backreference-list.paragraph.\the\secn.\the\parn\endcsname}%
	\expandafter\ifx\cseq\relax\else\edef\brt{\cseq}\printbackref{back references for paragraph.\the\secn.\the\parn: \cseq}\fi\fi
	\assignlabel{paragraph}{\the\secn.\the\parn}\fi}
\def\numpar{\ifnumb\nextpar\expand{\noexpand\typesetparnum{\the\secn.\the\parn}\noexpand\typesetpbr{\brt}}\fi}
\newdimen\pindent \pindent\parindent

\ifx\draftnum\undefined 
\def\gproclaim#1#2.{\curverb{#2~}
	\ifnumb\nextpar\fi
	{\everypar{}\noindent}%
	\plabel#1#2%
	\ifnumb\ \anchor{paragraph.\the\secn.\the\parn}{}\the\secn.\the\parn\fi.\enspace\mainfont 
	\ifnumb\expand{\noexpand\typesetpbr{\brt}}\fi
	\ignorespaces}
\def\typesetparnum#1{\ifnumb{\plabel\parnumfont\anchor{paragraph.#1}{}#1.\enspace}\fi} 
\def\typesetpbr#1{\ifnumb\def\brtext{#1}\ifx\brtext\empty\else\setetok\atendbr{{\parbackreffont Used in \noexpand\stripcomma\brtext.}}\fi\fi}
\else 
\let\numbfalse\relax \numbtrue
\def\gproclaim#1#2.{\curverb{#2~}\medbreak\noindent#1#2.\enspace\mainfont\ignorespaces}
\inewdimen\brwidth \brwidth.6in
\parindent0pt \parskip1ex plus 1ex minus 1ex
\def\typesetparnum#1{\ifnumb\llap{\plabel\anchor{paragraph.#1}{}\parnumfont#1\enspace}\fi} 
\def\typesetpbr#1{\ifnumb\def\brtext{#1}\ifx\brtext\empty\else
	\llap{\smash{\vtop{\everypar{}\raggedright\rightskip0pt plus 0pt \leftskip0pt plus 1fill \hsize\brwidth
	\parnumfont \strut \break 
	\parbackreffont\stripcomma#1}}\enspace}\fi\fi}
\fi

\def\hang{\hangindent\pindent}
\newif\iflist
\def\textindent#1{{\everypar{}\parindent\pindent\indent}\llap{#1\enspace}\listtrue\ignorespaces}

\def\li{\item{$\bullet$}}

\def\item{\endgraf\hang\textindent}
\def\filbreak{\endgraf\vfil\penalty-200\vfilneg}

\def\eject{\endgraf\break}
\def\supereject{\endgraf\penalty-20000}
\def\smallbreak{\endgraf\ifdim\lastskip<\smallskipamount
	\removelastskip\penalty-50\smallskip\fi}
\def\medbreak{\endgraf\ifdim\lastskip<\medskipamount
	\removelastskip\penalty-100\medskip\fi}
\def\bigbreak{\endgraf\ifdim\lastskip<\bigskipamount
	\removelastskip\penalty-200\bigskip\fi}


\let\printbackref\gobble
\def\predefbackref#1{%
	\printbackref{defining back reference list backreference-list.#1}%
	\expandafter\gdef\expandafter\cseq\expandafter{\csname backreference-list.#1\endcsname}
	\expandafter\ifx\cseq\relax\expandafter\gdef\cseq{}\else\printbackref{duplicate omitted}\fi
}
\newcount\backref \backref0
\newif\ifsuppressbackref \suppressbackreffalse
\def\firstletter#1#2\endletter{#1}
\newtoks\backreflist
\newif\ifaddbr \addbrfalse
\def\recordbackref#1{%
	\edef\params{{\ifpresec\expandafter\firstletter\chapname\endletter\else\the\secn\fi.\the\parn\ifnumb\else*\fi}{\the\backref}{\the\inputlineno}}%
	\printbackref{recording back reference \string#1 for future processing with params \params}%
	\edef\tmp{\the\backreflist\noexpand\processbackref\noexpand#1\params}%
	\global\backreflist\expandafter{\tmp}%
	\printbackref{new content of backreflist: \the\backreflist}} 
\def\processbackref#1#2#3#4{%
	\edef\key{\expandafter\gobble\string#1}%
	\edef\cseq{\csname id.\key\endcsname}%
	\expandafter\ifx\cseq\relax \ewarningline{#4}{Undefined reference \string#1.}{\string#1}\else
		\edef\lseq{\csname backreference-list.\cseq\endcsname}%
		\printbackref{processing back reference number #3 \string#1, originating from #2, at line #4; adding to \lseq}%
		\edef\lastnumber{\csname lastnumber.\cseq\endcsname}%
		\edef\newnumber{#2}%
		\addbrtrue
		\printbackref{lastnumber: \lastnumber; newnumber: \newnumber;}%
		\expandafter\ifx\csname lastnumber.\cseq\endcsname\relax 
		\else\ifx\lastnumber\newnumber 
			\printbackref{Suppressing duplicate back reference #2.}%
			\addbrfalse
		\fi\fi
		\expandafter\edef\csname lastnumber.\cseq\endcsname{#2}
		\ifaddbr
			\expandafter\expandafter\expandafter\gdef\expandafter\expandafter\csname backreference-list.\cseq\endcsname\expandafter{\lseq, \llink{backreference.#3}{#2}}
		\fi
	\fi
}
\newtoks\labelinitlist
\def\xxstripcomma, {}
\def\xstripcomma{\if\ntok,\let\xcont\xxstripcomma\else\let\xcont\relax\fi\xcont}
\def\stripcomma{\futurelet\ntok\xstripcomma}

\inewif\ifrecorddups \recorddupstrue
\def\checkduplicates#1#2{\edef\key{\expandafter\gobble\string#1}%
	\iftypesetting\else
	\expandafter\ifx\csname line:\key\endcsname\relax\printlabel{keydefline: relax}\ifrecorddups\expandafter\xdef\csname line:\key\endcsname{\the\inputlineno}\fi\else\edef\keydefline{\csname line:\key\endcsname}\errmessage{#2}\fi\fi}

\let\printlabel\gobble
\newtoks\curverb 
\ifx\draftlabel\undefined
\def\plabel{}
\else
\def\labeltext{} 
\def\plabel{\ifx\labeltext\empty\else\smash{\llap{\parbackreffont\labeltext\quad}}\gdef\labeltext{}\fi} 
\fi
\def\assignlabel#1#2{
	\ifx\lastlabel\undefined\else
	\edef\key{\expandafter\expandafter\expandafter\gobble\expandafter\string\lastlabel}
	\printlabel{label \key: id.\key\space = #1.#2, text.\key\space = #2}%
	\expandafter\xdef\csname id.\key\endcsname{#1.#2}
	\expandafter\xdef\csname text.\key\endcsname{#2}
	\edef\tmp{\the\curverb}%
	\ifx\tmp\empty\else 
	\printlabel{label v\key: id.v\key\space = #1.#2, text.v\key\space = \the\curverb#2}%
	\expandafter\xdef\csname id.v\key\endcsname{#1.#2}
	\expandafter\xdef\csname text.v\key\endcsname{\the\curverb#2}
	\fi
	\predefbackref{#1.#2}%
	\fi\let\lastlabel\undefined}
\newtoks\vlist 
\def\label#1{%
	\iftypesetting\else
	\checkduplicates#1{Label \string#1 was already defined at line \keydefline}%
	\addverunused#1\verifylabel
	\def\lastlabel{#1}%
	\fi
	\numbtrue 
}

\newif\ifyearkey 
\def\y{} 
\newdimen\bibindent 
\newtoks\bibt 
\def\tbib#1{
	\checkduplicates#1{Bibliographic reference \string#1 already defined at line \keydefline}%
	\addverunused#1\verifybib
	\edef\key{\expandafter\gobble\string#1}
	\printlabel{reference \key: id.\key\space = reference.\key, text.\key\space = \key}%
	\expandafter\edef\csname id.\key\endcsname{reference.\key}
	\expandafter\edef\csname text.\key\endcsname{\key}
	\predefbackref{reference.\key}%
	\ifyearkey\else\setbox0=\hbox{[\key]}\ifdim\bibindent<\wd0 \bibindent=\wd0 \fi \fi 
	\appendexpand\bibt{\noexpand\typesetbib\noexpand#1\src}%
	\ifyearkey\let\next\xxftbib\else\let\next\ftbibalpha\fi\next}
\def\ftbibalpha#1\par{\append\bibt{#1\par}}
\def\xxftbib{\futurelet\next\xftbib}
\def\xftbib{\if\next[\let\next\ftbibyear\else\let\next\ftbibnoyear\fi\next}
\def\ftbibyear[#1]{\edef\yearkey{#1}\expandafter\ftbibyearbis\ignorespaces}
\def\ftbibyearbis#1\par{\append\bibt{#1\par}\ftbibend}
\def\ftbibnoyear#1\par{\append\bibt{#1\par}\edef\yearkey{\extractyear#1\par}\ftbibend}
\def\ftbibend{\expandafter\ifx\csname year.\yearkey\endcsname\relax
		\edef\yearindex{0}%
	\else
		\edef\yearindex{\csname year.\yearkey\endcsname}%
	\fi
	{\count0=\yearindex
	\advance\count0 by 1 %
	\expandafter\xdef\csname year.\yearkey\endcsname{\the\count0 }%
	\xdef\alphakey{\ifcase \count0 ?\or a\or b\or c\or d\or e\or f\or g\or h\or i\or j\or k\or l\or m\or n\or o\or p\or q\or r\or s\or t\or u\or v\or w\or x\or y\or z\else .\the\count0 \fi}}%
	\printlabel{key \key, year key \yearkey, alpha key \alphakey}%
	\expandafter\edef\csname text.\key\endcsname{\noexpand\typesetyearalpha{\key}{\yearkey}{\alphakey}}%
	\printlabel{reference \key\space adjustment: id.\key\space = reference.\key, text.\key\space = \yearkey.\alphakey}%
	\setbox0=\hbox{[\yearkey.\alphakey]}\ifdim\bibindent<\wd0 \bibindent=\wd0 \fi
}
\def\typesetyearalpha#1#2#3{%
	\edef\yearindex{\csname year.#2\endcsname}%
	\ifnum\yearindex=1 #2\else#2.#3\fi}
\def\extractyear#1\y#2#3\par{#2}
\newif\iftype
\def\typesetbib#1#2\par{\edef\key{\expandafter\gobble\string#1}%
	\edef\bibbr{\csname backreference-list.reference.\key\endcsname}%
	\typetrue\ifsuppressunusedbib\ifx\bibbr\empty\typefalse
	\fi\fi
	\iftype
	\noindent\hbox to \bibindent{[\anchor{reference.\key}{\csname text.\key\endcsname}]\hfil}#2
	\ifx\bibbr\empty\else\expandafter\stripcomma\bibbr.\fi
	\hangindent\bibindent\filbreak
	\fi}

\let\printverify\gobble
\def\addverunused#1#2{\appendexpand\vlist{\noexpand#2\noexpand#1{\the\inputlineno}}}
\def\verifyref#1#2#3{\printverify{verifying for #3 \string#1 (line #2)}%
	\edef\key{\expandafter\gobble\string#1}%
	\edef\cseq{\csname id.\key\endcsname}%
	\edef\tmp{\csname backreference-list.\cseq\endcsname}%
        \ifx\tmp\empty\ewarningline{#2}{#3 \string#1.}{\string#1}\fi}
\def\verifylabel#1#2{\verifyref#1{#2}{Unused label}}
\def\verifybib#1#2{\verifyref#1{#2}{Unused reference}}

\let\printurl\gobble
\newtoks\urltext
\newtoks\urlt
\newif\ifpunct
\def\urldash{-}
\def\urltilde{{\tensy^^X}} 
\def\ndash{\def\urldash{--}}
\def\http://{\hfil\penalty900\hfilneg\urltext={http://}\urlt={http:/\negthinspace/}\punctfalse\urlgrab}
\def\https://{\hfil\penalty900\hfilneg\urltext={https://}\urlt={https:/\negthinspace/}\punctfalse\urlgrab}
\def\urlgrab{\catcode`\#=11 \catcode`\&=11 \futurelet\ntok\urldispatch}
\def\urldispatch{%
	\ifx\ntok~\let\proceed\urlcont\else
	\ifcat\noexpand\ntok\space\let\proceed\urlfinish\else
	\ifcat\noexpand\ntok\relax\let\proceed\urlfinish\else
	\let\proceed\urlcont
	\fi\fi\fi\proceed}
\def\urlcont#1{\ifpunct\appendexpand\urltext\punctc\appendexpand\urlt\punctc\punctfalse\fi
	\ifx\ntok~\appendexpand\urltext{\noexpand~}\appendexpand\urlt\urltilde
	\else\if\ntok\ampersand\appendexpand\urltext{&}\appendexpand\urlt{\&}%
	\else\if\ntok\ohash\appendexpand\urltext\ohash\appendexpand\urlt\#%
	\else\if\ntok_\appendexpand\urltext_\appendexpand\urlt\_%
	\else\if\ntok-\appendexpand\urltext-\appendexpand\urlt\urldash
	\else\if\ntok.\puncttrue\def\punctc{.}%
	\else\if\ntok,\puncttrue\def\punctc{,}%
	\else\if\ntok;\puncttrue\def\punctc{;}%
	\else\appendexpand\urltext{#1}\appendexpand\urlt{#1}%
	\fi\fi\fi\fi\fi\fi\fi\fi\urlgrab}
\def\urlfinish{\catcode`\#=6 \catcode`\&=4 \hbox{\printurl{\the\urltext}\link{\the\urltext}{\the\urlt}}\ifpunct\punctc\punctfalse\fi\def\urldash{-}}
\def\idgrab{\futurelet\ntok\iddispatch}
\def\iddispatch{\ifcat\noexpand\ntok\space\let\proceed\urlfinish
		\else\if\ntok,\let\proceed\urlfinish
		\else\let\proceed\idcont
		\fi\fi\proceed}
\def\idcont#1{\ifpunct\appendexpand\urltext.\appendexpand\urlt.\punctfalse\fi
	\if\ntok.\puncttrue\def\punctc{.}%
	\else\if\ntok_\appendexpand\urltext_\appendexpand\urlt\_%
	\else\appendexpand\urltext{#1}\appendexpand\urlt{#1}%
	\fi\fi\idgrab}

\let\printgrab\gobble
\newtoks\grabname
\newtoks\grabtoks 
\newtoks\grabcseq 
\newtoks\subsuptoks 
\newtoks\dtoks 
\newcount\grabsize
\newif\ifgrabsubscript 
\newif\ifgrabsupscript 
\def\grabsequence{\bgroup 
	\grabsubscripttrue\grabsupscripttrue 
	\grabstring}
\def\grabalpha{\bgroup 
	\grabsubscriptfalse\grabsupscriptfalse 
	\grabstring}
\def\grabingroup{\ifinfont\errmessage{Already inside a math token}\fi\append\grabtoks{\bgroup\grablink}\infonttrue}
\def\graboutgroup{\ifinfont\append\grabtoks{\endgrablink\egroup}\infontfalse\fi}
\def\grabstring#1#2#3{
	\let\specialhat^
	\catcode`\^=7
	\aftergroup#3 
	\ifx\specialaddon\undefined\else\expandafter\aftergroup\specialaddon\let\specialaddon\undefined\fi
	\inewif\ifdefine \inewif\ifinfont
	\printgrab{}\printgrab{grab a string of type #1, typeset using font \string#2, with postcommand \string#3}%
	\grabname{#1}\def\grabfont{#2}\grabsize0 \grabtoks={}\grabingroup \append\grabtoks{#2}\grabcseq={}%
	\futurelet\ntok\grabdeflookahead}
\def\grabdeflookahead{\if=\noexpand\ntok 
	\definetrue\printgrab{defining}\expandafter\grabgobblefuturelet
	\else\printgrab{referencing}\definefalse\expandafter\grablookahead\fi}
\def\grabgobblefuturelet#1{\futurelet\ntok\grabtestforsilent} 
\newif\ifsilentgrab
\def\grabtestforsilent{\if=\noexpand\ntok \silentgrabtrue \let\ncom\grabsilenteq \else \silentgrabfalse \let\ncom\grablookahead \fi \ncom}
\def\grabsilenteq={\grabfuturelet}
\def\grabfuturelet{\futurelet\ntok\grablookahead}
\def\grablookahead{\printgrab{futurelet token meaning: \meaning\ntok}%
	\let\ncom\grabfinish
	\if\bgroup\noexpand\ntok \printgrab{left brace, terminating}%
	\else \if\egroup\noexpand\ntok \printgrab{right brace, terminating}%
	\else \if\space\noexpand\ntok \printgrab{blank space, terminating}%
	\else \let\ncom\grabexamine \fi\fi\fi \ncom}
\def\grabexamine#1{\printgrab{grabexamine argument: \string#1, meaning \meaning#1}%
	\def\ncom{\grabfinish#1}%
	\ifcat$\ifcat*\string#1\fi$
		\ifcat _\noexpand#1 \ifgrabsubscript\printgrab{subscript, continuing}%
			\graboutgroup \append\grabtoks{#1}\subsuptoks{#1}\def\ncom{\grabsubsupfuturelet}%
						\else\printgrab{subscript, terminating}\fi
		\else \ifcat ^\noexpand#1 \ifgrabsupscript\printgrab{superscript, continuing}%
			\graboutgroup \append\grabtoks{#1}\subsuptoks{#1}\def\ncom{\grabsubsupfuturelet}%
						\else\printgrab{superscript, terminating}\fi
		\else \ifx \specialhat#1 \ifgrabsupscript\printgrab{specialhat superscript, continuing}%
			\graboutgroup \append\grabtoks{#1}\subsuptoks{#1}\def\ncom{\grabsubsupfuturelet}%
						\else\printgrab{superscript, terminating}\fi
		\else \ifcat\noexpand~\noexpand#1 \printgrab{active character \string#1, examining further}%
			\ifnum1=\uccode`#1 \printgrab{UTF-8 letter, continuing}%
				\advance\grabsize1 \append\grabtoks{#1}\appendexpand\grabcseq{\string#1}\def\ncom{\grabfuturelet}%
			\else
				\ifnum\the\grabsize=0 \printgrab{Nothing grabbed so far, continuing}%
					\advance\grabsize1 \append\grabtoks{#1}\appendexpand\grabcseq{\string#1}\def\ncom{\grabfuturelet}%
				\else\printgrab{Not a UTF-8 letter and not the first character in a string, terminating}%
				\fi
			\fi
		\else \ifcat a\noexpand#1 \printgrab{letter #1, continuing}%
			\advance\grabsize1 \append\grabtoks{#1}\append\grabcseq{#1}\def\ncom{\grabfuturelet}%
		\else\printgrab{nonactive character \string#1}%
			\ifnum\the\grabsize=0 \printgrab{sole argument, adding and terminating}%
				\advance\grabsize1 \append\grabtoks{#1}\append\grabcseq{#1}\def\ncom{\grabfinish}%
			\else\printgrab{terminating}\fi
		\fi\fi\fi\fi\fi
	\else \printgrab{command sequence \string#1, terminating}\fi
	\ncom}
\def\grabsubsupfuturelet{\futurelet\ntok\grabsubsuplookahead}
\newcount\dig
\newif\ifdigit
\def\grabsubsuplookahead{\printgrab{subsup futurelet token meaning: \meaning\ntok}%
	\if\bgroup\noexpand\ntok \printgrab{left brace, continuing}\let\ncom\grabentiresubsup%
	\else \if\egroup\noexpand\ntok \printgrab{right brace, continuing}\let\ncom\grabentiresubsup%
	\else \if\space\noexpand\ntok \errmessage{Blank space after \the\subsuptoks}%
	\else \let\ncom\grabsubsupexamine \fi\fi\fi \ncom}
\def\grabentiresubsup#1{\printgrab{subsup entire group added}\grabingroup\append\grabtoks{#1}\graboutgroup\grabfuturelet}
\def\grabsubsupexamine#1{\printgrab{examining subsup argument \string#1, meaning \meaning#1}%
	\ifcat$\ifcat*\string#1\fi$
		\ifcat\noexpand~\noexpand#1 \printgrab{active character \string#1, continuing}%
			\grabingroup\appendexpand\grabtoks{\grabfont\noexpand#1}\let\ncom\grabfuturelet
		\else\ifnum"8000=\the\mathcode`#1 \printgrab{math active character \string#1, continuing}%
			\let\specialaddon\egroup
			\def\ncom{#1}%
			\grabtypeset
		\else\ifcat a\noexpand#1 \printgrab{letter #1, checking whether single or not}%
			\dtoks{#1}\def\ncom{\futurelet\ntok\grabsubsupsecondletterlookahead}%
		\else\printgrab{something else, inserting a single-character sub/superscript, continuing}
			\appendexpand\grabtoks{\bgroup\grabfont\noexpand#1\egroup}%
			\advance\grabsize2 \appendexpand\grabcseq{\the\subsuptoks\string#1}
			\def\ncom{\grabfuturelet}\fi\fi\fi
	\else \printgrab{command sequence \string#1, continuing}%
		\append\grabtoks{#1}\let\ncom\grabfuturelet\fi
	\ncom}
\def\grabsubsupsecondletterlookahead{\def\ncom{\appendexpand\grabtoks{\the\dtoks}\grabfuturelet}%
	\ifcat a\noexpand\ntok \printgrab{not a single letter, grabbing the entire subsupscript}%
		\advance\grabsize1 \appendexpand\grabcseq{\expandafter\string\the\subsuptoks}
		\grabingroup
		\appendexpand\grabtoks{\grabfont\the\dtoks}%
		\appendexpand\grabcseq{\the\dtoks}%
		\def\ncom{\grabfuturelet}%
	\else \printgrab{single letter, continuing}\fi\ncom}
\def\grabfinish{\printgrab{grabfinish}\graboutgroup\grabtypeset\egroup}
\def\grabtypeset{\printgrab{grabtypeset grabsize=\the\grabsize, grabtoks=\the\grabtoks, grabcseq=\the\grabcseq}%
	\def\grablink##1\endgrablink{##1}%
	\ifnum\the\grabsize=0 \errmessage{No string to grab}\fi
	\ifnum\the\grabsize>1 
		\ifdefine 
			\expandafter\checkduplicates\csname\the\grabname.\the\grabcseq\endcsname{Mathematical identifier \key\space already defined at line \keydefline}%
			\iftypesetting 
				\ifsilentgrab
					\expandafter\gdef\csname silent:\the\grabname.\the\grabcseq\endcsname{}
				\else
					\anchor{\the\grabname.\the\grabcseq}{}
				\fi 
			\else
				\expandafter\xdef\csname id.\the\grabname.\the\grabcseq\endcsname{paragraph.\the\secn.\the\parn}
				\predefbackref{paragraph.\the\secn.\the\parn}%
  			\fi
		\else 
			\iftypesetting 
				\global\advance\backref1
				\expandafter\ifx\csname line:\the\grabname.\the\grabcseq\endcsname\relax 
					\warning{Undefined mathematical identifier \the\grabname.\the\grabcseq}%
					\expandafter\gdef\csname\the\grabname.\the\grabcseq\endcsname{\relax}
				\else 
					\expandafter\ifx\csname silent:\the\grabname.\the\grabcseq\endcsname\empty
						\edef\grablink##1\endgrablink{{##1}}%
					\else 
						\edef\grablink##1\endgrablink{\noexpand\llink{\the\grabname.\the\grabcseq}{##1}
							\noexpand\anchor{backreference.\the\backref}{}}
					\fi
				\fi
			\else 
				\ifsuppressbackref\else
					\global\advance\backref1
					\printbackref{math back reference to \the\grabname.\the\grabcseq: backref.\the\backref\space at line \the\inputlineno}%
					\expandafter\recordbackref\csname\the\grabname.\the\grabcseq\endcsname
				\fi
			\fi
		\fi
	\fi
	\ifsilentgrab\else\the\expandafter\grabtoks\fi}

\inewif\ifsilent \inewif\ifanchor
\newif\ifsuppresscs
\catcode`\^=7   \def\specialhat{\ifmmode\def\next{^}\else\let\next\beginxref\fi\next} \catcode`\^=\active \let^=\specialhat
\def\silentxref#1{\futurelet\next\silentxrefswitch}
\def\silentxrefswitch{\silenttrue\xref}
\def\beginxref{\futurelet\next\beginxrefswitch}
\def\beginxrefswitch{\ifx\next\specialhat\let\next\silentxref \else\silentfalse\let\next\xref\fi \next}
\def\xref{\leavevmode\futurelet\next\xrefswitch}
\def\xrefswitch{\ifx\next!\let\next\verbalxref \else \ifx\next=\let\next\anchorxref \else \anchorfalse \let\next\normalxref \fi \fi \next}
\newtoks\vtoksl
{\count0="C2 \loop\ifnum\count0<"F5 \catcode\count0=11 \advance\count0 by 1 \repeat 
\gdef\plainaccents{\suppresscstrue%
	\def\`##1{##1\empty ̀}%
	\def\'##1{##1\empty ́}%
	\def\^##1{##1\empty ̂}%
	\def\"##1{##1\empty ̈}%
	\def\~##1{##1\empty ̃}%
	\def\=##1{##1\empty ̄}%
	\def\.##1{##1\empty ̇}%
	\def\u##1{##1\empty ̆}%
	\def\v##1{##1\empty ̌}%
	\def\H##1{##1\empty ̋}%
	\def\t##1{##1\empty ͡}%
}}
\def\plainaccents{\let\xcsname=\empty \let\xendcsname=\empty}
\def\verbalxref!{\begingroup\plainaccents\verbalxrefaux}
\def\verbalxrefaux#1{%
	\lowercase{\vtoksl{#1}}%
	\expandtoks\vtoksl
	\iftypesetting
		\expandafter\gdef\expandafter\cseq\expandafter{\csname verbal.\the\vtoksl\endcsname}%
		\printlabel{verbal xref \the\vtoksl}%
		\expandafter\ifx\cseq\relax\errmessage{Undefined reference to \the\vtoksl}\else\cseq\fi
	\else
		\ifsuppressbackref\else
			\global\advance\backref1 %
			\printbackref{label verbal.#1: backref.\the\backref\space at line \the\inputlineno}%
			\blah
			\expandafter\recordbackref\csname verbal.\the\vtoksl\endcsname
		\fi
	\fi
	\endgroup}
\def\initverballabelcommand#1{%
	\printlabel{initializing verbal label #1 (\the\curverb\the\secn.\the\parn)}%
	\ifx\sectionid\undefined
		\expandafter\xdef\csname id.verbal.#1\endcsname{paragraph.\the\secn.\the\parn}%
	\else
		\expandafter\xdef\csname id.verbal.#1\endcsname{section.\the\secn}%
	\fi
	\ifnum\parn=0 %
	\expandafter\xdef\csname text.verbal.#1\endcsname{\the\curverb\the\secn}%
	\else
	\expandafter\xdef\csname text.verbal.#1\endcsname{\the\curverb\the\secn.\the\parn}%
	\fi
	\expandafter\initlabelcommand\csname verbal.#1\endcsname
}
\def\anchorxref={\anchortrue\futurelet\next\anchorxrefswitch}
\def\anchorxrefswitch{\ifx\next:\let\next\nonitalicanchor\else\italictrue\let\next\normalxref\fi \next}
\def\nonitalicanchor:{\italicfalse\normalxref}
\newtoks\firsttoks \newtoks\secondtoks \inewif\ifplural \inewif\ifitalic
\def\parseplural#1[#2|#3]{\let\next\parseplural\ifx\hfuzz#2\hfuzz\ifx\hfuzz#3\hfuzz\let\next\relax\else\pluraltrue\fi\else\pluraltrue\fi
	\append\firsttoks{#1#2}\append\secondtoks{#1#3}\next}
\newtoks\firsttoksl
\newtoks\secondtoksl
\newtoks\nexttoks
\def\normalxref{\begingroup\plainaccents\normalxrefaux}
\def\normalxrefaux#1{\firsttoks{}\secondtoks{}\pluralfalse\parseplural#1[|]%
	\lowercase\expandafter{\expandafter\firsttoksl\expandafter{\the\firsttoks}}%
	\expandtoks\firsttoksl
	\lowercase\expandafter{\expandafter\secondtoksl\expandafter{\the\secondtoks}}%
	\expandtoks\secondtoksl
	\ifanchor
		\iftypesetting\else
			\initverballabelcommand{\the\firsttoksl}%
			\ifplural\initverballabelcommand{\the\secondtoksl}\fi
			\predefbackref{paragraph.\the\secn.\the\parn}%
			\expandafter\checkduplicates\csname\the\firsttoksl\endcsname{Verbal label \the\firsttoksl\space was already defined at line \keydefline}
			\ifplural\expandafter\checkduplicates\csname\the\secondtoksl\endcsname{Verbal label \the\secondtoksl\space was already defined at line \keydefline}\fi
		\fi
		\anchor{verbal.\the\firsttoksl}{}
		\ifplural\anchor{verbal.\the\secondtoksl}{}\fi 
	\else
		\ifsuppressbackref\else
			\global\advance\backref1
			\printbackref{label normal.#1: backref.\the\backref\space at line \the\inputlineno}%
			\iftypesetting
			\else
				\expandafter\recordbackref\csname verbal.\the\secondtoksl\endcsname
			\fi
			\anchor{backreference.\the\backref}{}
		\fi
	\fi
	\nexttoks{}%
	\ifsilent
		\nexttoks{\ignorespaces}%
	\else
		\ifanchor
			\ifitalic\nexttoks\expandafter{\expandafter\bgroup\expandafter\it\the\firsttoks\italcorr}%
			\else\nexttoks\expandafter{\the\firsttoks}%
			\fi
		\else
			\edef\tmp{{verbal.\the\secondtoksl}}%
			\nexttoks\expandafter{\expandafter\llink\tmp}
			\nexttoks\expandafter\expandafter\expandafter{\expandafter\the\expandafter\nexttoks\expandafter{\the\firsttoks}}%
		\fi
	\fi
	\expandafter\endgroup\the\nexttoks}
\def\italcorr{\futurelet\next\italcorrtest}
\def\italcorrtest{\if,\noexpand\next\else\if.\noexpand\next\else\/\fi\fi\egroup}



\def\gen:{\http://libgen.rs/book/index.php?md5=}
\def\jstor:{\https://www.jstor.org/stable/}
\def\eudml:{\https://eudml.org/doc/}

\def\arXiv:{\urltext={https://arxiv.org/abs/}\urlt={arXiv:}\punctfalse\idgrab}

\def\Zbl:{\urltext={https://zbmath.org/?q=an:}\urlt={Zbl:}\punctfalse\idgrab}
\def\doi:{\ndash\urltext={https://doi.org/}\urlt={doi:}\punctfalse\urlgrab}

\def\sqr#1#2{{\thinspace\vbox{\hrule height.#2pt \hbox{\vrule width.#2pt height#1pt \kern#1pt \vrule width.#2pt} \hrule height0pt depth.#2pt}\thinspace}}
\def\square{\mathchoice\sqr64\sqr64\sqr{4.2}3\sqr33}

\def\ltoarr#1{\mathop{\count0=#1 \loop\ifnum\count0>0 \smash-\mkern-7mu \advance\count0 -1 \repeat \mathord\rightarrow}\limits} 
\def\lto#1#2{\mathrel{\ltoarr{#1}^{#2}}} 
\def\longto#1^#2_#3{\mathrel{\ltoarr{#1}^{#2}_{#3}}} 
\def\lgetsarr#1{\mathop{\mathord\leftarrow \count0=#1 \loop\ifnum\count0>0 \mkern-7mu\smash-\advance\count0 -1 \repeat}\limits} 
\def\longgets#1^#2_#3{\mathrel{\lgetsarr{#1}\limits^{#2}_{#3}}} 


\def\toto{\mathrel{\vcenter{\hbox{$\to$}\kern-1.5ex \hbox{$\to$}}}}
\def\prearrfill{\smash-\mkern-7mu}
\def\postarrfill{\mkern-7mu\smash-}
\def\midarrfill#1{\cleaders\hbox{$\mkern-2mu\smash-\mkern-2mu$}\hskip0pt plus #1fil}
\def\rightarrfill{\mkern-7mu\mathord\rightarrow}
\def\leftarrfill{\mathord\leftarrow\mkern-7mu\midarrfill1\postarrfill}
\def\ltoto#1#2#3{\ifinner
	\mathrel{\vcenter{\hbox to #1em{$\prearrfill\midarrfill1{\scriptstyle#2}\midarrfill3 \rightarrfill$}%
		\kern-1.5ex \hbox to #1em{$\prearrfill\midarrfill3{\scriptstyle#3}\midarrfill1 \rightarrfill$}}}%
	\else
	\mathrel{\mathop{\vcenter{\hbox to #1em{\rightarrowfill}%
		\kern-1.5ex \hbox to #1em{\rightarrowfill}}}\limits^{#2}_{#3}}%
	\fi}
\def\ltogets#1#2#3{\ifinner
	\mathrel{\vcenter{\hbox to #1em{$\prearrfill\midarrfill1{\scriptstyle#2}\midarrfill3 \rightarrfill$}%
		\kern-1.5ex \hbox to #1em{$\leftarrfill\midarrfill3{\scriptstyle#3}\midarrfill1 \postarrfill$}}}%
	\else
	\mathrel{\mathop{\vcenter{\hbox to #1em{\rightarrowfill}%
		\kern-1.5ex \hbox to #1em{\leftarrowfill}}}\limits^{#2}_{#3}}%
	\fi}

\def\ltogetscore#1#2{\dimen0=\fontdimen6 #1 2 \divide\dimen0 by 2 \multiply\dimen0 by #2 \vcenter{\hbox to \dimen0{\rightarrowfill}\kern-1.8ex \hbox to \dimen0{\leftarrowfill}}}
\def\ltogets#1#2#3{\mathrel{\mathop{\mathchoice{\ltogetscore\textfont{#1}}{\ltogetscore\textfont{#1}}{\ltogetscore\scriptfont{#1}}{\ltogetscore\scriptscriptfont{#1}}}^{#2}_{#3}}}

\def\rx#1#2{\rlap{\kern #1pt \raise#1pt \hbox{#2}}}
\def\dottednearrow{\rx{-8}. \rx{-6}. \rx{-4}. \rx{-2}. \rx0. \rx2. \rx4. \kern6pt \raise7.7pt \hbox{$\nearrow$}}

\def\gmatrix#1#2{\null\,\vcenter{\normalbaselines
	\ialign{#1\crcr
		\mathstrut\crcr\noalign{\kern-\baselineskip}
		#2\crcr\mathstrut\crcr\noalign{\kern-\baselineskip}}}\,}
\def\cdmatrix{\gmatrix{\hfil$##$\hfil&&\enspace\hfil$##$\hfil\enspace&\hfil$##$\hfil}}
\def\sqmatrix{\gmatrix{\hfil$##$&\enspace\hfil$##$\hfil\enspace&$##$\hfil}}
\def\cdbl{\def\normalbaselines{\baselineskip20pt \lineskip3pt \lineskiplimit3pt }}

\def\cd{\cdbl\cdmatrix}
\def\sqcd{\cdbl\let\vagap\;\sqmatrix}

\newcount\arrowsize \arrowsize3
\def\mapright#1{\smash{\lto\arrowsize{#1}}}

\def\rvagap{\vagap}  \def\rvaskip{\vaskip}  \def\vaskip{} \def\vagap{}
\def\mapdown#1{\rvagap\Big\downarrow\rlap{$\vcenter{\hbox{$\scriptstyle#1$}}$}\rvaskip}

\newcount\forno \forno0
\def\arrno#1#2{\global\advance\arr1 \edef\eeqnno{\the\arr}%
	\global\advance\forno1 \edef\eforno{\the\forno}%
	\xdef#2{\noexpand\llink{equation.\eforno}{\eeqnno}}
	#1{(\anchor{equation.\eforno}{\eeqnno})}} 
\newbox\mdiag
\def\wrapdiagram{%
	\setbox\mdiag\vtop\bgroup
	\null 
	\vskip\baselineskip
	\inewcount\arr \arr0
	\baselineskip0pt
	\lineskip4pt
	\lineskiplimit4pt
	\let\par\cr
	\obeylines
	\halign\bgroup\hfil$\displaystyle##$\hfil\cr
	\ewrapdiagram}

\def\ewrapdiagram#1{#1
	\egroup
	\egroup
	\vskip0pt plus \dp\mdiag \penalty-250 \vskip0pt plus-\dp\mdiag 
	\hangafter-\dp\mdiag
	\divide\hangafter\baselineskip
	\advance\hangafter-2
	\hangindent-\wd\mdiag
	\advance\hangindent-2em
	\hbox to\hsize{\hfil\dp\mdiag0pt \box\mdiag}%
	\ignorespaces}

\newlinechar=10
\inewtoks\preamble
{\endlinechar=10 \catcode`#=12 \global\preamble{
prologues := 3;

verbatimtex
\let\endprolog

\expandafter\gobbleinit\input }\global\appendonceexpand\preamble\jobname\global\append\preamble{
\catcode`\^=7
etex

input cmarrows
setup_cmarrows(arrow_name = "texarrow"; parameter_file = "cmr10.mf"; macro_name = "drawarrow");
setup_cmarrows(arrow_name = "doublearrow"; parameter_file = "cmr10.mf"; macro_name = "drawdarrow");
def drawmarrow expr p = _apth:=p; _finmarr enddef;
rule_thickness#:=.4pt#;    
def _finmarr text t_ =
  drawarrow subpath(0, 0.5 * length(_apth)) of _apth t_;
  draw subpath(0.5 * length(_apth), length(_apth)) of _apth withpen pencircle scaled rule_thickness# t_;
enddef;

def object(suffix O)(expr x,y)(expr l) =
  save O;       
  pair O;
  O := (x,y) * u;
  picture O.tx;    
  O.tx := thelabel(l,O);
  draw O.tx;                         
enddef;                                    

def smorphism(suffix A,B) =
  save ss, tt;
  ss := xpart ((A..B) intersectiontimes bbox A.tx);
  tt := xpart ((A..B) intersectiontimes bbox B.tx);
  drawarrow subpath(ss,tt) of (A..B);
enddef;

def morphism(suffix A,B)(expr l)(expr f)(suffix $) =
  smorphism(A, B);
  label.$(l, point f[ss, tt] of (A..B));
enddef;
}}%

\newcount\vertex                                      
\ifx\mathspecials\undefined\def\mathspecials{}\fi
\def\diagram{\begininlinemp\buffertoks={}\grabdiagram}
\def\grabdiagram{\mathspecials\grabdiagramaux}                                                 
\def\grabdiagramaux#1;{\vertex0 \dcom#1,*}                          
\def\dcom{\futurelet\next\dcomswitch}                                          
\def\dcomswitch{\ifcat a\noexpand\next \let\next\dcomalpha                                       
        \else\if*\noexpand\next \addcode{\the\buffertoks}\endinlinemp\let\next\gobble                                                 
        \else\if[\noexpand\next \let\next\grabscale                                                                
        \else\errmessage{Unrecognized diagram command \next}\let\next\relax\fi\fi\fi\next}
\def\grabscale[#1]{\dimen0=#1 \addcode{save u; u = \the\dimen0;}\dcom}
\def\dcomalpha#1{\def\objectname{#1}\futurelet\next\dcomalphaswitch}
\def\dcomalphaswitch{\ifcat a\noexpand\next \let\next\dcommorphism\else                                                          
        \expandafter\edef\csname vertex:\objectname\endcsname{\the\vertex}%
        \if:\noexpand\next\let\next\grabcoords                                              
        \else\if=\noexpand\next\let\next\grabobjectequ                                     
        \else\errmessage{Expected : or = while processing a diagram object, got \meaning\next}\let\next\relax\fi\fi\fi\next}
\def\grabcoords:#1,#2={\toks\vertex{#1,#2}\grabobjectlabel}
\def\grabobjectequ={\edef\tmp{\the\toks\vertex}\ifx\tmp\empty\errmessage{No coordinates specified for vertex \the\vertex: \objectname}\fi\grabobjectlabel}
\def\grabobjectlabel#1,{\addcode{object(\objectname, \the\toks\vertex, }\addunexpandedcode{btex #1 etex);}\advance\vertex1 \dcom}
\def\dcommorphism#1{\def\tobjectname{#1}%
        \def\labelpos{.5}%
        \futurelet\next\dcommorphismswitch}                     
\def\dcommorphismswitch{\if.\noexpand\next \expandafter\grabmorphismdir \else \setupmorphismlabel \expandafter\grabmorphismpos\fi}
\def\setupmorphismlabel{
        \edef\vlabeldira{direction.\objectname.\tobjectname}%
        \edef\vlabeldirb{direction.\tobjectname.\objectname}%
        \expandafter\ifx\csname\vlabeldira\endcsname\relax
        	\expandafter\ifx\csname\vlabeldirb\endcsname\relax
		        \edef\tmpa{\csname vertex:\objectname\endcsname}\expandafter\ifx\tmpa\relax\errmessage{No such vertex: \objectname}\fi
		        \edef\tmpb{\csname vertex:\tobjectname\endcsname}\expandafter\ifx\tmpb\relax\errmessage{No such vertex: \tobjectname}\fi
		        \ifnum \tmpa<\tmpb \edef\vlabeldir{direction.\tmpa.\tmpb}\else \edef\vlabeldir{direction.\tmpb.\tmpa}\fi
		        \expandafter\ifx\csname\vlabeldir\endcsname\relax\errmessage{No label direction specified for morphism \objectname->\tobjectname}\fi
		        \edef\labeldir{\csname\vlabeldir\endcsname}%
		\else
        		\edef\labeldir{\csname\vlabeldirb\endcsname}%
		\fi
	\else
        	\edef\labeldir{\csname\vlabeldira\endcsname}%
	\fi
}                                                                                        
\def\grabmorphismdir.{\def\labeldir{}\futurelet\next\grabmorphismdirec}
{\catcode`\@=13                                                               
\gdef\grabmorphismdirec{\ifx@\next \let\next\grabmorphismposition
        \else\if=\noexpand\next \let\next\grabmorphismequ               
        \else \let\next\grabmorphismdirect\fi\fi\next}                       
\gdef\grabmorphismdirect#1{\edef\labeldir{\labeldir#1}\futurelet\next\grabmorphismdirec}
\gdef\grabmorphismpos{\ifx@\next \expandafter\grabmorphismposition \else \expandafter\grabmorphismequ\fi} 
\gdef\grabmorphismposition@#1={\def\labelpos{#1}\grabmorphismlabel}
}
\def\grabmorphismequ={\grabmorphismlabel}
\def\grabmorphismlabel#1,{%
        \addcodebuffer{morphism(\objectname, \tobjectname, }%
        \addunexpandedcodebuffer{btex \everymath{\scriptstyle}#1 etex, }%
        \addcodebuffer{\labelpos, \labeldir);}%
        \dcom}

\let\printextract\gobble
\def\hinitlabelcommand#1{\printextract{initializing label command \string#1}%
	\gdef#1{\printlabel{invoked label \string#1}
		\ifsuppressbackref
			\edef\key{\expandafter\gobble\string#1}%
			\llink{\csname id.\key\endcsname}{\csname text.\key\endcsname}
		\else
			\global\advance\backref1 %
			\printbackref{label \expandafter\gobble\string#1: backref.\the\backref\space at line \the\inputlineno}%
			\iftypesetting
			\else
				\blah
				\recordbackref#1
			\fi
			\anchor{backreference.\the\backref}{}
			\edef\key{\expandafter\gobble\string#1}%
			\llink{\csname id.\key\endcsname}{\csname text.\key\endcsname}
		\fi
		}}
\let\initlabelcommand\hinitlabelcommand
\def\pinitlabelcommand#1{\printextract{initializing label command \string#1}%
	\gdef#1{\printlabel{invoked plain label \string#1}
		\iftypesetting
			\edef\key{\expandafter\gobble\string#1}%
			\csname text.\key\endcsname
		\else
			\blah
		\fi}}
\newread\labelin
\newif\iflabelcont
\let\terminate=\relax 
\long\def\labelauxaux#1\terminate
{}
\def\preprocesslabel#1{%
	\printextract{Processing label \string#1}%
	\edef\key{\expandafter\gobble\string#1}%
	\initlabelcommand#1
	\expandafter\initlabelcommand\csname v\key\endcsname
	\labelauxaux
}
\def\preprocessbib#1{%
	\printextract{Processing bib \string#1}%
	\initlabelcommand#1
	\labelauxaux
}
\long\def\labelaux#1{\ifx#1\label\let\next\preprocesslabel\else\ifx#1\bib\let\next\preprocessbib\else\let\next\labelauxaux\fi\fi\next}
\def\processoneline{\expandafter\labelaux\labelline\relax\relax\relax\terminate}
\def\preprocess#1{
	\openin\labelin=#1 \labelconttrue
	\loop
		\read\labelin to\labelline
		\ifeof\labelin\let\next\labelcontfalse\else\let\next\processoneline\fi
		\next
	\iflabelcont\repeat
	\expandafter\gobbleinit\input#1
	\ifhmode\par\fi\vfill\eject 
	\ifhmode\par\fi\vfill\supereject
}

\def\importlabels#1{ 
	\recorddupsfalse
	\let\initlabelcommand\pinitlabelcommand
	\preprocess{#1}%
	\let\initlabelcommand\hinitlabelcommand
	\recorddupstrue
	\cont={}
	\vlist={}
	\let\chapname\undefined \secn0 \parn0 \backref0 
}


\newif\iftypesetting 

\typesettingfalse 
\def\blah{blah} 
\output{\setbox0\box255 \setbox0\box\footins \deadcycles0 } 
\let\bib\tbib 
\let\refs\relax 
\everypar{\numpar}
\parn0
\def\par{\finishpar} 
\newif\ifmetapost \metapostfalse
\edef\filestem{\jobname.gen}
\newcount\figno
\newwrite\mpout
\newtoks\outtoks
\def\inimp{\ifmetapost\else\global\metaposttrue\immediate\openout\mpout=\filestem.mp \addcode{\the\preamble}\fi}
\figno0
\newdimen\oldhsize
\oldhsize=\hsize
\newdimen\oldvsize
\oldvsize=\vsize
\hsize\maxdimen
\vsize\maxdimen
\hbadness10000
\preprocess\jobname 
\hsize\oldhsize
\vsize\oldvsize
\ifx\fmtname\plainfmtname
\hsize\plaintextwidth
\vsize\plaintextheight
\fi
\hbadness1000
\printerr{(\jobname.tex}
\printbackref{list of all back references: \the\backreflist}%
\the\backreflist 
\ifsuppressunusedbib\def\verifybib#1#2{}\fi
\the\vlist 
\printerr{)}
\ifmetapost
\addcode{end}
\immediate\closeout\mpout
\immediate\write18{mpost -interaction nonstopmode \filestem.mp}
\ifeof18 \warning{Compile the METAPOST file \filestem.mp manually using mpost \filestem.mp}\let\runmp\warning\fi
\fi
\let\addcode\gobble 

\typesettingtrue
\def\importlabels#1{} 
\ifscroll \vsize\maxdimen\inewbox\abox\output{\setbox\abox\vbox{\unvbox255\unskip}\shipbox\abox}%
\else 
	\ifx\fmtname\plainfmtname
	\totalht\plaintextheight
	\advance\totalht.1in
	\advance\totalht2\vmargin
	\totalwd\plaintextwidth
	\advance\totalwd2\hmargin
	\setpapersize\totalwd\totalht
	\hoffset-1in
	\advance\hoffset\hmargin
	\voffset-1in
	\advance\voffset\vmargin
	\fi
	\output{\plainoutput}
\fi
\let\blah\undefined 
\def\bib#1\par{} 
\def\refs{\raggedright\rightskip0em plus \maxdimen \advance\bibindent1em \everypar{}\the\bibt\ignorespaces} 
\def\addverunused#1#2{} 
\def\addcont#1#2#3{} 
\def\prepass#1{} 
\parn0
\parbreftrue 
\everypar{\numpar}
\let\chapname\undefined \secn0 \backref0 
\figno0
\expandafter\gobbleinit\input\jobname\relax
\ifhmode\par\fi\vfill\eject 
\ifhmode\par\fi\vfill\supereject
\end


\catcode`^=7
\mathcode`@="8000 \catcode`@=13 \def@{\grabalpha{cat}\sf\empty} \catcode`@=12 
\mathcode`"="8000 \catcode`"=13 \def"{\grabalpha{fun}\eurm\empty} \catcode`"=12 
\def\mathspecials{\catcode`@=13 \catcode`"=13 \relax}
\catcode`"=12

\mathchardef\togets="181D
\def\newunicodechar#1#2{\expandafter\def\csname\string#1\endcsname{#2}}

\chardef\lasyfam=8 
\font\tenlasy=lasy10 \font\sevenlasy=lasy7 \font\fivelasy=lasy5 \textfont\lasyfam=\tenlasy \scriptfont\lasyfam=\sevenlasy \scriptscriptfont\lasyfam=\fivelasy

\mathchardef\unlhd="3802
\mathchardef\unrhd="3804

\catcode`^=\active

\mathcode`\:="603A 
\mathchardef\colon="303A 

{\catcode`\'=\active \gdef'{^\bgroup\primes}}
\def\primes{\prime\futurelet\next\primesbis}
\def\primesbis{\ifx'\next\let\nxt\primesbisbis
	\else\ifx^\next\let\nxt\primet
	\else\ifx\specialhat\next\let\nxt\primet
  	\else\let\nxt\egroup
	\fi\fi\fi \nxt}
\def\primesbisbis#1{\primes}
\def\primet#1#2{#2\egroup}


\inewtoks\titletoks
\inewtoks\authortoks
\titletoks{Combinatorial model categories are equivalent to presentable quasicategories}
\authortoks{Dmitri Pavlov}

\def\unicodechars{
\newunicodechar{⇄}{\togets}
\newunicodechar{⊴}{\unlhd}
\newunicodechar{⊵}{\unrhd}
}
\def\step#1.{\smallbreak\noindent{\bf#1.\enspace}\ignorespaces}

\unicodechars
\metadata{\the\titletoks}{\the\authortoks}
{\articletitle \leftskip0pt plus 20em \rightskip0pt plus 20em \parfillskip0pt \parindent0pt \baselineskip20pt \the\titletoks \bigskip}
{\tabskip0pt plus 1fil \let\par\cr \obeylines \halign to\hsize{\hfil#\hfil
\chaptertitle \the\authortoks\vadjust{\medskip}
Department of Mathematics and Statistics, Texas Tech University
\https://dmitripavlov.org/
}}

\begingroup\setbox0=\hbox{$@==op "==Sd "==Ex "==op @==Cat$}\endgroup


\def\define#1{\expandafter\def\csname#1\endcsname{{\sf#1}}}
\def\absconv{
	\define{sph}
	\define{odisk}
	\define{ndisk}
	\define{disks}
	\define{relh}
	\define{old}
	\define{new}
}
\absconv

\def\colim{\mathop{\rm colim}} 
\def\Clu{{\sf Clu}} 
\def\op{{\sf op}} 
\def\id{{\rm id}} 
\def\N{{\rm N}} 
\def\rN{{\cal N}}
\def\rK{{\cal K}}
\def\id{{\rm id}}
\def\RN{{\eucalbf N}}
\def\RK{{\eucalbf K}}
\def\frc{{\frakbf R}} 
\def\FRC{{\eucalbf R}} 
\def\ham{{\eucalbf H}} 
\def\K{\mathop{\sf K\null}\nolimits} 

\def\zz{{\eucalbf Z}}
\def\Ind{{\sf Ind}} 
\def\Exi{\mathop{\sf Ex^\infty}\nolimits} 

\endprolog

\hyphenation{
Grothen-dieck
co-do-main
com-po-nent-wise
lem-ma
pre-stack
pre-stacks
pre-sheaf
pre-sheaves
man-u-script
pseudo-functor
trac-ta-ble
mon-oidal
cat-e-gor-i-cal
pseu-do-func-tor
}


\abstract Abstract.
We establish a Dwyer–Kan equivalence of relative categories
of combinatorial model categories, presentable quasicategories, and other models for locally presentable $(∞,1)$-categories.
This implies that the underlying quasicategories of these relative categories are also equivalent.

\tsection Contents

\contents

\section Introduction

Combinatorial model categories and presentable quasicategories are the two most used formalisms for locally presentable $(∞,1)$-categories.
It has long been conjectured that these formalisms should be equivalent in a certain sense,
see, for example, Problems 8 and~11 on Hovey's algebraic topology problem list [\Hovey].
^!{main theorem} is a solution to (one precise formulation of) these problems.

Partial results in this direction existed for a long time, see, in particular, the work of Dugger [\Pres] and Lurie [\HTT].
For example, we know that the underlying quasicategory of a combinatorial model category is presentable,
and up to an equivalence of quasicategories, every presentable quasicategory arises in this manner.
Likewise, the underlying functor of quasicategories of a left Quillen functor
is a left adjoint functor between presentable quasicategories,
and up to an equivalence of functors, every left adjoint functor between presentable quasicategories arises in such a manner.

However, locally presentable $(∞,1)$-categories can themselves be organized into an $(∞,1)$-category,
so it is natural to inquire whether the resulting $(∞,1)$-categories of combinatorial model categories and presentable quasicategories are equivalent.
In this article, we formalize these $(∞,1)$-categories as relative categories and prove the following result.

\proclaim Theorem.
^^={main theorem}
The following relative categories are Dwyer–Kan equivalent.
In particular, their underlying quasicategories and homotopy (2,1)-categories are equivalent.
\li The relative category $@CMC$ of combinatorial model categories, left Quillen functors, and left Quillen equivalences.
\li The relative category $@CRC$ of combinatorial relative categories, homotopy cocontinuous relative functors, and ^{Dwyer–Kan equivalences}.
\li The relative category $@PrL$ of presentable quasicategories, left adjoint functors, and equivalences.
\endlist
These equivalences are implemented in two flavors:
\li Working in the Zermelo–Fraenkel set theory,
we have a ^{Dwyer–Kan equivalence} of relative categories $@CMC$ (^!{CMC}), $@CRC$ (^!{CRC}), $@PrL$ (^!{PrL}).
\li Assuming the existence of a strongly inaccessible cardinal~$U$,
we have a ^{Dwyer–Kan equivalence} of relative categories $@CMC'_U$ (^!{universe CMCU}), $@CRC'_U$ (^!{universe CRCU}), $@PrL'_U$ (^!{universe PrLU}).
\endlist
Furthermore, these equivalences are compatible with each other, as explained in the proof.

\proof Proof.
Combine ^!{Cof equivalence} and ^!{derived relative nerve equivalence} to establish the ^{Dwyer–Kan equivalences} $"Cof:@CMC→@CRC$ (^!{Cof}) and $\RN:@CRC→@PrL$ (^!{derived relative nerve}).
We also have ^{Dwyer–Kan equivalences} $@CMC_U→@CMC'_U$ (^!{smaller universe CMCU}), $@CRC_U→@CRC'_U$ (^!{smaller universe CRCU}), $@PrL_U→@PrL'_U$ (^!{smaller universe PrLU}),
where $@CMC_U$ (^!{smaller CMCU}), $@CRC_U$ (^!{smaller CRCU}), $@PrL_U$ (^!{smaller PrLU})
are certain full subcategories of $@CMC$, $@CRC$, and $@PrL$.
Restricting the ^{Dwyer–Kan equivalences} $"Cof$ and $\RN$
to the corresponding full subcategories
establishes ^{Dwyer–Kan equivalences} $"Cof_U:@CMC_U→@CRC_U$ and $\RN_U:@CRC_U→@PrL_U$.
As shown in ^!{universe Cof RInd} and ^!{smaller CRCU universe PrLU}, these equivalences are compatible
with certain naturally defined relative functors $"Cof_U:@CMC'_U→@CRC'_U$ and $\RN_U:@CRC'_U→@PrL'_U$.

The following theorem is a solution to (one precise formulation of) Problem~9 on Hovey's list [\Hovey].

\proclaim Theorem.
^^={left proper main theorem}
(See ^!{main theorem simplicial}.)
The following relative categories are Dwyer–Kan equivalent (and hence also equivalent to the categories in ^!{main theorem}).
\li The relative category $@CMC$ of combinatorial model categories and left Quillen functors.
\li Left proper combinatorial model categories and left Quillen functors.
\li Simplicial combinatorial model categories and simplicial left Quillen functors.
\li Simplicial left proper combinatorial model categories and simplicial left Quillen functors.
\endlist
In all four cases, weak equivalences are given by left Quillen equivalences.

\proclaim Theorem.
(See ^!{main theorem cartesian}.)
The following relative categories are Dwyer–Kan equivalent.
\li Cartesian combinatorial model categories, left Quillen functors, and left Quillen equivalences.
\li Same as the previous item, but additionally required to be simplicial (with simplicial left Quillen functors), or left proper, or both.
\li Cartesian closed presentable quasicategories.

We also compare the resulting constructions to derivators.
Derivators by their nature are not fully homotopy coherent, so some truncation must be performed.
Given a relative category~$C$ we can extract from it its homotopy (2,1)-category
by taking the hammock localization $\ham_C$ of~$C$
and replacing each simplicial hom-object $\ham_C(X,Y)$ with its fundamental groupoid.

\proclaim Theorem.
The following (2,1)-categories are equivalent.
\li The homotopy (2,1)-category of presentable quasicategories.
\li The homotopy (2,1)-category of combinatorial model categories.
\li The homotopy (2,1)-category of combinatorial left proper model categories.
\li The (2,1)-category of presentable derivators, left adjoints, and isomorphisms.

\proof Proof.
The equivalence of the first and second (2,1)-categories follows from ^!{main theorem}.
The equivalence of the second and third (2,1)-categories follows from ^!{left proper main theorem}.
For the equivalence of the third and fourth (2,1)-categories, see Renaudin [\Derivators, Theorem~3.4.4].

\proclaim Remark.
The ^{Dwyer–Kan equivalence} between $@CRC$ and $@PrL$ shown in ^!{derived relative nerve equivalence} works abstractly
with any pair of Quillen equivalent models for $(∞,1)$-categories,
since all what is used in ^!{derived relative nerve equivalence} is a Quillen equivalence $\rN⊣\rK$
together with a fibrant replacement functor $\frc$ and a cofibrant replacement functor
(i.e., the identity functor for the Joyal model structure, but a nontrivial functor for other models).
In particular, the same proof establishes ^{Dwyer–Kan equivalences}
between appropriate versions of relative categories
of complete Segal spaces, Segal categories, simplicial categories, etc.
We do not include proofs in this paper because
doing so would require us to develop notions of homotopy colimits, homotopy ind-completions, and homotopy local presentability
in each of these settings, and then show their compatibility with each other.
However, one can also transport these notions from a model where they are already
developed (such as quasicategories) along derived Quillen equivalences
connecting quasicategories to whatever models we are interested in.
Indeed, this is essentially how we defined objects, morphisms, and weak equivalences of $@CRC$.
With this convention, ^!{derived relative nerve equivalence} immediately yields Dwyer—Kan equivalences
of the relative subcategories of complete Segal spaces, Segal categories, simplicial categories, marked simplicial sets, quasicategories,
relative categories, and other models of $(∞,1)$-categories,
once we replace $@CRC$
with an analogously defined relative category where we take as objects the relevant model of a homotopy locally presentable category
and as morphisms the relevant model of a homotopy cocontinuous functor.

\subsection Previous work

Lurie [\HTT, Proposition~A.3.7.6] shows that any presentable quasicategory is equivalent to
the homotopy coherent nerve of the category of bifibrant objects of a combinatorial simplicial model category.
In the same proposition, he shows that the underlying quasicategory of a model category of simplicial presheaves
is a presentable quasicategory.
Combined with Dugger [\Pres, Propositions 3.2 and 3.3], this shows that the underlying quasicategory
of a combinatorial model category is a presentable quasicategory.
For another exposition, see Cisinski [\HCHA, Theorem~7.11.16, Remark~7.11.17].

The work of Quillen [\HoAlg],
Maltsiniotis [\ThQ],
Lurie [\HTT, Corollary~A.3.1.12],
Hinich [\DKLR, Proposition~1.5.1],
Mazel-Gee [\QAdj, Theorem~2.1]
shows that a Quillen adjunction between model categories (with finite limits and finite colimits)
induces an adjunction of quasicategories.
For another exposition, see Cisinski [\HCHA, Theorem~7.5.30].

Renaudin [\Derivators, Theorem~3.4.4]
proves that the functor from the localization of the 2-category
of combinatorial left proper model categories
at left Quillen equivalences
to the 2-category
of presentable derivators,
left adjoints,
and modifications
is an equivalence of 2-categories.
Arlin [\Preder, Theorems 4.1, 5.1,~6.4] establishes analogous results for quasicategories.
Low [\Hobi, Theorem~4.15] establishes an equivalence of bicategories of complete Segal spaces and quasicategories,
based on the work of Riehl–Verity [\TwoCat] on the 2-category of quasicategories.
Szumiło [\HTCQ] establishes an equivalence between the fibration categories
of cocomplete quasicategories and cofibration categories.

Rezk–Schwede–Shipley [\SimStruc, Theorem 1.1]
show that any left proper cofibrantly generated model category
that satisfies a certain realization axiom introduced there,
is Quillen equivalent to a simplicial model category.
In Theorem~1.2 they show that the existence of a Quillen equivalence between such model categories
implies the existence of a simplicial Quillen equivalences.
Dugger [\Replacing, Theorem~1.2] proves that a left proper combinatorial model category is Quillen equivalent to a simplicial left proper combinatorial model category.
Dugger [\Pres, Corollary~1.2] drops the left properness assumption from the previous theorem.



\subsection Prerequisites

We assume familiarity with basics of the following topics from homotopy theory.
Appropriate references will be given throughout the text.

\li Locally presentable and accessible categories,
including regular cardinals, $λ$-filtered colimits, $λ$-accessible categories, $λ$-accessible functors, locally $λ$-presentable categories,
$λ$-presentable objects (denoted by $\K_λ$), $λ$-ind-completions (denoted by $\Ind^λ$),
the sharp ordering of regular cardinals (denoted by $κ⊲λ$).
See Gabriel–Ulmer [\LPK], Makkai–Paré [\AccCat], and Adámek–Rosický [\LPAC].
\li Simplicial homotopy theory, including simplicial sets, simplicial maps, simplicial weak equivalences, and the simplicial Whitehead theorem (^!{simplicial Whitehead}).
See Goerss–Jardine [\SHT] and Dugger–Isaksen [\WESP].
\li Model categories, including model structures, left Quillen functors, projective model structures on presheaves, Reedy model structures, left Bousfield localizations.
See Hovey [\MC], Hirschhorn [\MCL], and Barwick [\LR].
\li Relative categories, including relative functors, simplicial categories, hammock localizations (denoted by $\ham$).
See Dwyer–Kan [\SimpLoc, \Calc, \Func] and Barwick–Kan [\RelCat].
\li Quasicategories, including the Joyal model structure, limits and colimits,
ind-completions, and presentable quasicategories.
See Joyal [\QCKC], Lurie [\HTT, \Kerodon], and Cisinski [\HCHA].

\subsection Further directions

We expect the methods developed in this paper to be applicable to other similar statements, some of which are indicated in ^!{conjecture symmetric monoidal} and ^!{conjecture enriched}.
Considerations of length prevent us from including proofs in this article.

\proclaim Conjecture.
^^={conjecture symmetric monoidal}
The relative functor
from the relative category of combinatorial symmetric monoidal model categories
to the relative category of closed symmetric monoidal presentable quasicategories
that sends a monoidal model category to its underlying symmetric monoidal quasicategory
is a ^{Dwyer–Kan equivalence} of relative categories.
In particular, the underlying quasicategories are also equivalent.

Nikolaus–Sagave [\NS, Theorem~1.1, Theorem~2.8] show that the underlying symmetric monoidal quasicategory functor
is homotopy essentially surjective and homotopy full on 1-morphisms.

\proclaim Conjecture.
^^={conjecture enriched}
Fix a combinatorial symmetric monoidal model category~$V$.
The relative functor
from the relative category of combinatorial $V$-enriched model categories
to the relative category of presentable $V$-enriched quasicategories
that sends an enriched model category to its underlying enriched quasicategory
is a ^{Dwyer–Kan equivalence} of relative categories.
In particular, the underlying quasicategories are also equivalent.

Haugseng [\RectEC, Theorem~5.8] shows that the underlying quasicategory of the relative category of $V$-enriched small categories, $V$-enriched functors,
and Dwyer–Kan equivalences
is equivalent to the quasicategory of $\cal V$-enriched small quasicategories,
where $\cal V$ is the underlying symmetric monoidal quasicategory of~$V$.

\subsection Acknowledgments

I thank Urs Schreiber for a discussion that led to this paper.

\section Preliminaries

In this paper, we adopt a convention that a category need not be small or locally small.

\proclaim Definition.
A ^={categor[y|ies]} is given by a class~$O$ of objects, a class~$M$ of morphisms, together with source, target, identity, and composition maps that satisfy the usual axioms.
In particular, morphisms $X→Y$ between objects $X$,~$Y$ in a category~$C$ can form a proper class $C(X,Y)$.
A ^={locally small categor[y|ies]} is a category~$C$ such that for any objects $X,Y∈C$, the class $C(X,Y)$ is a set.
A ^={small categor[y|ies]} is a category~$C$ such that the class of objects~$O$ and the class of morphisms~$M$ are both sets.
An ^={essentially small categor[y|ies]} is a category~$C$ that is equivalent to a small category.

\proclaim Definition.
Suppose $λ$ is a regular cardinal.
A $λ$-^={small set[|s]} is a set~$X$ of cardinality strictly less than~$λ$.
A {\it$λ$-small category\/} is a small category~$C$ whose set of morphisms is a $λ$-small set.

\subsection Accessible categories
^^={accessible categories}

We now review some of the more specialized definitions from Low [\Heart].
A good example of a category~$C$ to keep in mind is the category of $λ$-presentable objects in some locally presentable category, which is an essentially small category,
in fact a small category according to ^!{smaller locally presentable}.
Thus, ^!{moderate objects} can be seen as defining analogues of the usual notions (like that of a $κ$-presentable object, $κ$-accessible category, locally $κ$-presentable category)
in the setting of small categories whose objects are limited in size by some larger regular cardinal~$λ$.
This relationship is further clarified by ^!{accessibly generated characterization}.

\proclaim Definition.
^^={moderate objects}
(Low [\Heart, Definition~1.2].)
Given regular cardinals $κ≤λ$,
a $(κ,λ)$-^={presentable object[|s]} in a ^{locally small category}~$C$ is an object $A∈C$ such that the functor $C(A,-):C→@Set$ preserves $λ$-small $κ$-filtered colimits.
The full subcategory of $(κ,λ)$-presentable objects in~$C$ is denoted by $\K_κ^λ(C)$.
A {\it$κ$-presentable object\/} is an object that is $(κ,λ)$-presentable for all regular cardinals~$λ$.
The full subcategory of $κ$-presentable objects in~$C$ is denoted by $\K_κ(C)$.
A $(κ,λ)$-^={accessibly generated categor[y|ies]} (Low [\Heart, Definition~3.2]) is an essentially small category
that admits $λ$-small $κ$-filtered colimits
and every object is the colimit of some $λ$-small $κ$-filtered diagram of $(κ,λ)$-presentable objects.
A ^={locally $(κ,λ)$-presentable category} is a $(κ,λ)$-accessibly generated category that admits $λ$-small colimits.

\proclaim Proposition.
^^={accessibly generated characterization}
(Low [\Heart, Theorem~3.11].)
If $κ⊴λ$ are regular cardinals,
then for an idempotent-complete essentially small category~$C$ the following conditions are equivalent:
\li $C$ is a $(κ,λ)$-^{accessibly generated category};
\li $\Ind^λ(C)$ is a $κ$-accessible category;
\li $C$ is equivalent to $\K_λ(D)$ for some $κ$-accessible category~$D$.

\proclaim Proposition.
^^={locally presentable characterization}
If $κ≤λ$ are regular cardinals,
then for an idempotent-complete essentially small category~$C$ the following conditions are equivalent:
\li $C$ is a ^{locally $(κ,λ)$-presentable category};
\li $\Ind^λ(C)$ is a locally $κ$-presentable category;
\li $C$ is equivalent to $\K_λ(D)$ for some locally $κ$-presentable category~$D$.

\proclaim Definition.
(Low [\Heart, Definition~2.2].)
Given a regular cardinal~$κ$, a
^={strongly $κ$-accessible functor[|s]}
^^={strongly $λ$-accessible}
^^={strongly accessible functor[|s]}
^^={strongly accessible}
is a functor between $κ$-accessible categories
that preserves $κ$-filtered colimits and $κ$-^{presentable objects}.
Given regular cardinals $κ≤λ$,
a ^={strongly $(κ,λ)$-accessible functor[|s]}
is a functor between $(κ,λ)$-accessibly generated categories
that preserves $λ$-small $κ$-filtered colimits and $(κ,λ)$-presentable objects.

\proclaim Proposition.
^^={accessible functors are strong}
(Adámek–Rosický [\LPAC, Theorem~2.19].)
Every accessible functor~$"F$ is ^{strongly $λ$-accessible} for arbitrarily large regular cardinals~$λ$:
if $κ$ is a regular cardinal, there is a regular cardinal~$λ⊵κ$ such that $"F$ is ^{strongly $λ$-accessible}.

\proclaim Proposition.
^^={accessible functors small}
Given a regular cardinal~$κ$ and $κ$-accessible categories $C$ and $D$,
the functors $\Ind^κ$ and $\K_κ$ induce an equivalence of groupoids between the groupoid of functors $\K_κ(C)→\K_κ(D)$
and the groupoid of strongly $κ$-accessible functors $C→D$.
Given a regular cardinal~$κ$ and locally $κ$-presentable categories $C$ and $D$,
the functors $\Ind^κ$ and $\K_κ$ induce an equivalence of groupoids between the groupoid of functors $\K_κ(C)→\K_κ(D)$ that preserve $κ$-small colimits
and the groupoid of strongly $κ$-accessible left adjoint functors $C→D$.

\proclaim Proposition.
^^={bounded accessible functors small}
Given regular cardinals~$κ⊲λ$ and $κ$-accessible categories $C$ and $D$,
the functors $\Ind^κ_λ$ and $\K_κ^λ$ induce an equivalence of groupoids between the groupoid of functors $\K_κ^λ(C)→\K_κ^λ(D)$
and the groupoid of strongly $(κ,λ)$-accessible functors $C→D$.
Given regular cardinals~$κ≤λ$ and locally $(κ,λ)$-presentable categories $C$ and $D$,
the functors $\Ind^κ_λ$ and $\K_κ^λ$ induce an equivalence of groupoids between the groupoid of functors $\K_κ^λ(C)→\K_κ^λ(D)$ that preserve $κ$-small colimits
and the groupoid of strongly $(κ,λ)$-accessible left adjoint functors $C→D$.

\subsection Size aspects
^^={size aspects}

In this article we use the Zermelo–Fraenkel set theory with the axiom of choice.
In particular, we do not assume any large cardinal axioms,
since we intend the results of this paper
to be usable in papers that do not assume any additional axioms.

One subtlety that emerges from this decision is that
the three main relative categories of this paper
(combinatorial model categories, combinatorial relative categories, and presentable quasicategories)
must be defined with more care than usual,
since such categories typically have a proper class of objects,
so cannot themselves be elements of a class or objects in a category.

In what follows, we would like to use the functor $\K_λ$ to construct small categories.
A priori, if $C$ is a locally presentable category, then $\K_λ(C)$ is an essentially small category that is not necessarily small.
We circumvent this problem in ^!{small category of compact objects}.
To this end, we identify a full and essentially surjective subcategory $@Set$ of the category of sets
such that for any regular cardinal~$λ$, the full subcategory of $@Set$ consisting of $λ$-small sets in $@Set$ is a small category.

Recall that the {\it rank\/} of a set~$S$ is defined inductively on~$S$
as the smallest ordinal greater than the rank of all elements of~$S$.
The axiom of foundation guarantees that the induction makes sense.
Alternatively, the rank of~$S$ is the smallest ordinal~$α$ such that $S⊂V_α$,
where $V_α$ is von Neumann's cumulative hierarchy: $V_0=∅$, $V_{α+1}=2^{V_α}$,
$V_β=⋃_{α<β}V_α$, where $β$ is a limit ordinal.

\proclaim Definition.
^^={smaller sets}
Denote by $@=Set$ the full subcategory of the category of sets and maps of sets
on objects given by sets~$S$ whose rank does not exceed their cardinality.

\proclaim Remark.
The inclusion of the category $@Set$ (^!{smaller sets}) into the category of sets and maps of sets is an equivalence of categories,
since every set is bijective with a cardinal.

\proclaim Remark.
^^={smaller locally presentable}
Every locally presentable category is equivalent to a full subcategory~$C$ of a category of presheaves of sets on a small category.
If we interpret sets as objects of $@Set$ (^!{smaller sets}), then for every regular cardinal~$λ$, the category $\K_λ(C)$ is a small category.
Thus, from now on we require (without a loss of generality) that a locally presentable category~$C$
satisfies the following condition: for any regular cardinal~$λ$, the category $\K_λ(C)$ is a small category.
This convention will be used throughout this article.

Recall the following variant (due to Ehresmann, see Beurier–Pastor–Guitart [\Clusters, Definition~4.1])
of the free cocompletion construction of a locally small category~$C$, which defines a (strict) endofunctor
on locally small categories.

\proclaim Definition.
^^={strict cocompletion}
Given a locally small category~$C$, the {\it strict free cocompletion\/} of~$C$ is a category $\Clu(C)$ defined as follows (Beurier–Pastor–Guitart [\Clusters, Theorem~3.9]).
Objects are small diagrams $I→C$ (where a small category~$I$ is a category internal to $@Set$ as in ^!{smaller sets}).
Morphisms are ^={cluster[s|]}, as defined in Beurier–Pastor–Guitart [\Clusters, Definition~3.1].

\proclaim Remark.
The set of morphisms $P→Q$ of clusters is canonically isomorphic
to $$\lim_{p∈P}\colim_{q∈Q} C(P(p),Q(q))$$ (Beurier–Pastor–Guitart [\Clusters, Proposition~3.11]).
The resulting category $\Clu(C)$ is equivalent to the category of small presheaves on~$C$.
The construction presented here has an advantage that it is manifestly strictly functorial in~$C$.

\proclaim Definition.
^^={compact strict cocompletion}
If $μ$ is a regular cardinal and $C$ is a locally small category,
we denote by $\Clu_μ(C)$ the full subcategory of $\Clu(C)$
on diagrams whose indexing category is $μ$-small.

\proclaim Remark.
^^={small category of compact objects}
If $μ$ is a regular cardinal and $C$ is a locally small category,
the canonical inclusion of the category $\Clu_μ(C)$ into the full subcategory of $μ$-presentable objects in $\Clu(C)$ is an equivalence of categories.
If $C$ is a small category, then the domain $\Clu_μ(C)$ is a small category thanks to ^!{smaller sets},
whereas the codomain $\K_μ(\Clu(C))$ is not a small category, but merely an essentially small category,
since any diagram indexed by a category with a terminal object is a compact object,
so there is a proper class of compact objects.
In particular, if $C$ is a small category, the category $\Clu(C)$ is not locally presentable in the sense of ^!{smaller locally presentable}.

\proclaim Definition.
^^={smaller ind-completion}
Suppose $λ$ is a regular cardinal and $C$ is a category.
We define the locally $λ$-presentable (in the sense of ^!{smaller locally presentable}) category $\Ind^λ(C)$ as the full subcategory of the strict free cocompletion $\Clu(C)$ of~$C$ (^!{strict cocompletion})
on diagrams whose indexing category is a small $λ$-filtered category.

We need a size-restricted variant of ^!{smaller ind-completion}.

\proclaim Definition.
^^={small ind-completion}
Suppose $λ$ and $μ$ are regular cardinals, $λ≤μ$, and $C$ is a small category that admits $λ$-small colimits.
We define the category $\Ind^λ_μ(C)$ as the Gabriel–Zisman category of fractions
of the category $\Clu_μ(C)$ (^!{compact strict cocompletion}) with respect to morphisms inverted by the left adjoint functor
$$\Clu_μ(C)→\Ind^λ_μ(C)$$
induced by the universal property of $\Clu_μ(C)$ from the canonical inclusion $C→\Ind^λ_μ(C)$.
The functor $\smash{\Ind^λ_μ}$ is a (strict) functor from the category of small $λ$-cocomplete categories to the category of small $μ$-cocomplete categories.
We refer to $\Ind^λ_μ(C)$ as the {\it $(λ,μ)$-ind-completion of~$C$}.
The canonical inclusion $C→\Ind^λ_μ(C)$ is given by the constant diagram functor.

\proclaim Remark.
The canonical functor $\Ind^λ_μ(C)→\K_μ(\Ind^λ(C))$ is an equivalence of categories.
We could also define $\Ind^λ_μ(C)$ without appealing to categories of fractions
as the full subcategory of $\Clu(C)$ on $μ$-small $λ$-filtered diagrams,
but such a definition would make ^!{U-ind-completion} invalid, since $λ$-filtered colimits are typically not $μ$-filtered if $μ>λ$.

We conclude this section by examining the 1-categorical analogue of the three main relative categories of this article:
$@CMC$ (^!{CMC}), $@CRC$ (^!{CRC}), and $@PrL$ (^!{PrL}).
Informally, we want to define the category $@LPC$ of locally presentable categories and left adjoint functors.

Except for trivial cases, a left Quillen equivalence between combinatorial model categories never admits an inverse that is a left Quillen functor.
Thus, when we later define $@CMC$ in ^!{CMC} we are naturally forced to use the notion of a category equipped with a subcategory of weak equivalences,
i.e., a relative category (Barwick–Kan [\RelCat]).
Furthermore, the objects we are interested in have a 2-categorical nature
and we must take into account the notion of a equivalence between morphisms, e.g., left adjoint functors can be naturally isomorphic.
A common approach to this is to use 2-categories, defined by Bénabou in 1965,
which would lead us to develop a theory of relative 2-categories.
However, relative categories themselves can encode higher homotopy groups for hom-objects
by virtue of using appropriately chosen weak equivalences.
Thus, we stay in the realm of 1-categories and encode all structures as relative categories.

The relative category $@LPC$
can be informally described as the relative category
of locally presentable categories, left adjoint functors, and equivalences of categories.
This naive definition does not make sense in the usual ZFC set theory without large cardinal axioms
because proper classes (such as the class of objects of a locally presentable category that is not a poset)
cannot be elements of other classes.
We circumvent the problem by observing that a locally presentable category or a left adjoint functor between locally presentable categories
can be specified using sets only, without referring to classes.
The two fundamental facts that we need are as follows (see ^!{locally presentable characterization}, ^!{accessible functors are strong}, and ^!{accessible functors small} for a precise formulation):
\li Any locally presentable category is the $λ$-ind-completion of a small category~$C$ that admits $λ$-small colimits, for some regular cardinal~$λ$.
\li Any left adjoint functor between locally presentable categories is the $μ$-ind-completion of a functor between small categories, for some (possibly larger) regular cardinal~$μ$.

\proclaim Definition.
^^={smaller LPC}
The relative category $@=LPC$ of locally presentable categories is defined as follows.
Objects are pairs $(λ,C)$, where $λ$ is a regular cardinal and $C$ is a small category that admits $λ$-small colimits.
Morphisms $(λ,C)→(μ,D)$ exist if $λ≤μ$, in which case they are functors $C→D$ that preserve $λ$-small colimits.
Morphisms are composed by composing their underlying functors.
Weak equivalences are morphisms $(λ,C)→(μ,D)$ such that the functor $C→D$
is equivalent to the canonical inclusion $C→\Ind^λ_μ(C)$ (^!{small ind-completion}).

In particular, for $λ=μ$, a weak equivalence $(λ,C)→(μ,D)$ is simply an equivalence of categories,
since the $(λ,λ)$-ind-completion is equivalent to the idempotent completion and categories that admit $λ$-small colimits are automatically idempotent complete.

\proclaim Remark.
Given $(λ,C),(μ,D)∈@LPC$, the hom-object $\ham_@LPC((λ,C),(μ,D))$
is weakly equivalent to the nerve of the groupoid
whose objects are $λ$-cocontinuous functors $C→\Ind^μ(D)$
and morphisms are natural isomorphisms.
The groupoid of $λ$-cocontinuous functors $C→\Ind^μ(D)$ is equivalent to the groupoid of left adjoint functors $\Ind^λ(C)→\Ind^μ(D)$.
(This statement must be interpreted in terms of relevant constructions
implementing the functors in both directions and the unit and counit isomorphisms,
since left adjoint functors $\Ind^λ(C)→\Ind^μ(D)$ cannot be organized into a class.)
In particular, $@LPC$ indeed behaves like the purported $(2,1)$-category of locally presentable categories,
left adjoint functors, and natural isomorphisms.
We omit the proof of this claim since it is not used in the rest of the paper.

\subsection Universes

Although we do not assume any large cardinal axioms for our main results,
we find it useful to formulate explicit comparison results
for existing definitions of the categories $@CMC'_U$ (^!{universe CMCU}), $@CRC'_U$ (^!{universe CRCU}), $@PrL'_U$ (^!{universe PrLU}) that use large cardinals.
For the purposes of formulating these three comparison results, it suffices to assume the existence of a strongly inaccessible cardinal, i.e., a Grothendieck universe.

We start with the simpler definition, assuming the existence of a strongly inaccessible cardinal~$U$.
The following definition collects the pertinent adjustments to the notions of category theory that rely on the distinction between sets and classes.

\proclaim Definition.
^^={size notions}
Suppose $U$ is a strongly inaccessible cardinal.
\li A ^={$U$-small set[|s]} is a set of rank less than~$U$.
\li A ^={$U$-small class[|es]} is a set whose elements are $U$-small sets.
\li A ^={$U$-small categor[y|ies]} is a category whose classes of objects and morphisms are $U$-small sets.
\li A ^={locally $U$-small categor[y|ies]} is a category whose classes of objects and morphisms are $U$-small classes, and hom-classes between any pair of objects are $U$-small sets.
\li A ^={$U$-essentially $U$-small categor[y|ies]} is a locally $U$-small category that is equivalent to a $U$-small category.
\li A ^={$U$-locally $U$-presentable categor[y|ies]} is a locally $U$-small category~$C$ such that for some regular cardinal $κ<U$ the category $\K_κ^U(C)$ is a $U$-essentially $U$-small category that admits $κ$-small colimits and the inclusion $\K_κ^U(C)→C$ is a $(κ,U)$-ind-cocompletion functor.
\li Using $U$-small categories, we define
$U$-small limits and colimits (using $U$-small diagrams),
$U$-complete and $U$-cocomplete categories (using $U$-small limits and colimits),
free strict $U$-cocompletion (^!{compact strict cocompletion}).
\li The notions of ^!{accessible categories} specialized to the case $λ=U$
yield appropriate notions of 
$(κ,U)$-presentable objects,
$(κ,U)$-ind-completions,
$(κ,U)$-accessible categories,
strongly $(κ,U)$-accessible functors,
locally $(κ,U)$-presentable categories.
\li The above notions are extended to quasicategories
by defining a ^={$U$-small quasicategor[y|ies]} to be a quasicategory~$X$ such that for every $n≥0$ the set $X_n$ is $U$-small,
a ^={locally $U$-small quasicategor[y|ies]} to be a quasicategory~$X$ such that for every $n≥0$ the set $X_n$ is a $U$-small class
and the fibers of the vertex map $X_n→X_0^{n+1}$ are $U$-small sets,
and promoting the remaining definitions to the setting of quasicategories.

We now use these definitions to define a simpler version $@LPC'_U$ of the relative category $@LPC$ (^!{smaller LPC})
whose objects are actual categories, as opposed to pairs $(λ,C)$ that we used for $@LPC$.
The price to pay is that the relative category $@LPC'_U$ depends on the strongly inaccessible cardinal~$U$
in an essential way (^!{universe warning}).

\proclaim Definition.
^^={universe LPCU}
Assuming $U$ is a strongly inaccessible cardinal,
the relative category $@LPC'_U$ is defined as follows.
Objects are $U$-locally $U$-presentable categories (^!{size notions}).
Morphisms are left adjoint functors.
Weak equivalences are equivalences of categories.

\proclaim Remark.
^^={U-presentable}
We remark that for every strongly inaccessible cardinal~$U$,
every $C∈@LPC'_U$, and every regular cardinal $κ<U$, the category $\K_κ^U(C)$ is a $U$-essentially $U$-small category (^!{size notions}).
Conversely, if $D$ is a $U$-essentially $U$-small category,
then $\Ind^λ_U(D)$ is a $U$-locally $U$-presentable category
because $U$-small diagrams in a locally $U$-small category form a locally $U$-small category.

\proclaim Remark.
^^={universe adjustment}
Suppose $U<U'$ are strongly inaccessible cardinals.
The category $@LPC'_U$ is equivalent (as a relative category) to the subcategory $@LPC'^U_{U'}$ of the category $@LPC'_{U'}$
whose objects are categories $C∈@LPC'_{U'}$ such that $\K_U^{U'}(C)∈@LPC'_U$
and morphisms are strongly $(U,U')$-accessible left adjoint functors.
The functor $@LPC'^U_{U'}→@LPC'_U$
sends $C↦\K_U^{U'}(C)$ and a strongly $(U,U')$-accessible left adjoint functor to its restriction to $(U,U')$-presentable objects.
The functor $@LPC'_U→@LPC'^U_{U'}$
sends a category~$C$ to a variant of the $(U,U')$-ind-completion of~$C$ (^!{small ind-completion})
given by taking the full subcategory of $\Clu_{U'}(C)$ on $U$-small diagrams
as well as $U'$-small $U$-filtered diagrams whose indexing category does not contain a final subcategory of smaller cardinality.
(The latter condition ensures that $(U,U')$-presentable objects in the resulting category are given by $U$-small diagrams, hence form a locally $U$-small category.)

\proclaim Warning.
^^={universe warning}
Suppose $U<U'$ are strongly inaccessible cardinals.
If a category~$C$ belongs to both $@LPC'_U$ and $@LPC'_{U'}$, then it is a preorder.
Also, as is clear from ^!{universe adjustment}, the categories $@LPC'_U$ and $@LPC'_{U'}$ are not equivalent.
This explains why in ^!{smaller universe LPCU}, both relative categories must depend on~$U$.

In order to compare the relative categories $@LPC$ (^!{smaller LPC}) and $@LPC'_U$ (^!{universe LPCU}), we must take ^!{universe warning} into account
and modify the relative category $@LPC$ to ensure that the resulting relative category $@LPC_U$ can be weakly equivalent to the category $@LPC'_U$.

\proclaim Definition.
^^={smaller LPCU}
Given a strongly inaccessible cardinal~$U$,
the relative category $@LPC_U$ is defined as
the full subcategory of $@LPC$ (^!{smaller LPC})
on objects $(λ,C)$, where
$λ<U$ and $C$ is a $U$-essentially $U$-small category (^!{size notions}).

In order to compare the relative categories $@LPC_U$ (^!{smaller LPCU}) and $@LPC'_U$ (^!{universe LPCU}), we define a comparison functor between them.
The need to define $\Ind_U$ in ^!{U-ind-completion} as a strict functor justifies the somewhat convoluted construction in ^!{small ind-completion} of the $(λ,μ)$-ind-completion $\Ind^λ_μ(C)$
as a category of fractions of $\Clu_μ(C)$, instead of constructing it directly as the full subcategory of $\Clu(C)$ on $μ$-small $λ$-filtered colimits.

\proclaim Definition.
^^={U-ind-completion}
Assuming $U$ is a strongly inaccessible cardinal,
the relative functor
$$\Ind_U:@LPC_U→@LPC'_U$$
between the relative categories $@LPC_U$ (^!{smaller LPCU}) and $@LPC'_U$ (^!{universe LPCU}) is defined
as follows.
The functor $\Ind_U$ sends an object $(λ,C)∈@LPC_U$ to the category $\Ind^λ_U(C)$ (^!{small ind-completion}),
which is $U$-locally $U$-presentable by ^!{U-presentable}.
The functor $\Ind_U$ sends a morphism $"G:(λ,C)→(μ,D)$ to the functor $\Ind^λ_U(C)→\Ind^μ_U(D)$
that sends a $U$-small diagram $d:I→C$ to the $U$-small diagram $"G∘d:I→D$.
The functor $\Ind_U$ is a relative functor because it sends a weak equivalence $(λ,C)→(μ,\Ind^λ_μ(C))$ to the equivalence $\Ind^λ_U(C)→\Ind^μ_U(\Ind^λ_μ(C))$.

The following proposition and its proof serve as a base for the three comparison results ^!{smaller universe CMCU}, ^!{smaller universe CRCU}, ^!{smaller universe PrLU}.

\proclaim Proposition.
^^={smaller universe LPCU}
Assuming $U$ is a strongly inaccessible cardinal,
the relative functor
$$\Ind_U:@LPC_U→@LPC'_U$$
of ^!{U-ind-completion} is a ^{Dwyer–Kan equivalence} of relative categories.

\proof Proof.
^{Dwyer–Kan equivalences} of relative categories are stable under filtered colimits.
We introduce filtrations on $@LPC_U$ and $@LPC'_U$
that are respected by the functor $\Ind_U$.
We then show that $\Ind_U$ induces a ^{Dwyer–Kan equivalence} on every step of the filtration.
\ppar
Fix a regular cardinal~$ν$.
Define $@LPC_{U,ν}$ as the full subcategory of $@LPC_U$ consisting of objects $(λ,C)$ for which $λ≤ν$.
Define $@LPC'_{U,ν}$ as the full subcategory of $@LPC'_U$ consisting of objects given by locally $(ν,U)$-presentable categories~$C$
and morphisms given by strongly $(ν,U)$-accessible functors (^!{strongly accessible functor}).
By ^!{U-presentable} and ^!{bounded accessible functors small}, the functor $\Ind_U$ restricts to a functor $$\Ind_{U,ν}:@LPC_{U,ν}→@LPC'_{U,ν}.$$
\ppar
To show that $\Ind_U$ is a ^{Dwyer–Kan equivalence}, it suffices to construct a functor
$$\K_{U,ν}:@LPC'_{U,ν}→@LPC_{U,ν}$$
together with natural weak equivalences
$$η:\id_{@LPC_{U,ν}}→\K_{U,ν}∘\Ind_{U,ν}$$
and
$$ε:\id_{@LPC'_{U,ν}}∘ι→\Ind_{U,ν}∘\K_{U,ν}∘ι,$$
where $ι$ is a Dwyer–Kan equivalence of relative categories constructed below.
\ppar
The functor $\K_{U,ν}$ sends $C∈@LPC'_{U,ν}$ to the object $(ν,\K_ν^U(C))∈@LPC_{U,ν}$, where $\K_ν^U(C)$ is $U$-essentially $U$-small by definition of $@LPC'_{U,ν}$.
(At this point, the presence of a filtration is crucial: without having $ν$ at our disposal, we would not be able to define the first component of an object in $@LPC_U$ in a functorial way.)
The functor $\K_{U,ν}$ sends a functor $"F:C→D$ in $@LPC'_{U,ν}$ to the restriction $$(ν,\K_ν^U(C))→(ν,\K_ν^U(D)),$$
which is well-defined because $"F$ is a strongly $(ν,U)$-accessible functor.
\ppar
The natural weak equivalence $$η:\id_{@LPC_{U,ν}}→\K_{U,ν}∘\Ind_{U,ν}$$
is given on an object $(λ,C)∈@LPC_{U,ν}$
by the embedding $$(λ,C)→(ν,\K_{U,ν}(\Ind_{U,ν}(λ,C)))=(ν,\K_ν^U(\Ind^λ_U(C)))$$
that sends an object $X∈C$ to the singleton diagram $X:1→C$.
The inclusion functor
$$C→\K_ν^U(\Ind^λ_U(C))≃\Ind^λ_ν(C)$$
is the $(λ,ν)$-ind-completion functor,
hence the constructed morphism is indeed a weak equivalence.
\ppar
We would like to construct a natural weak equivalence $$ε:\id_{@LPC'_{U,ν}}→\Ind_{U,ν}∘\K_{U,ν}$$
by sending an object $D∈@LPC'_{U,ν}$
to the embedding $$D→\Ind_{U,ν}(\K_{U,ν}(D))=\Ind_U(ν,\K_ν^U(D))=\Ind^ν_U(\K_ν^U(D))$$
that sends an object $d∈D$ to its canonical diagram indexed by the comma category $\K_ν^U(D)/d$.
Since $D$ is locally $(ν,U)$-presentable, this morphism is indeed an equivalence of categories.
Unfortunately, $\K_ν^U(D)$ is not $U$-small in general, only $U$-essentially $U$-small,
which means that the comma category $\K_ν^U(D)/d$ does not produce an object in $\Ind^ν_U(\K_ν^U(D))$.
\ppar
Consider the full relative subcategory $ι:@LPC''_{U,ν}→@LPC'_{U,ν}$
on objects $D∈@LPC'_{U,ν}$ such that $D$ is a skeletal category (in particular, $\K_ν^U(D)$ is a $U$-small category).
The above construction produces a natural weak equivalence
$$ε':ι→\Ind_{U,ν}∘\K_{U,ν}∘ι.$$
\ppar
It remains to show that $ι$ is a Dwyer–Kan equivalence of relative categories.
By construction, $ι$ is homotopically essentially surjective.
By ^!{hammock construction hocolim}, it suffices to show that for every zigzag type~$Z$ and objects $X,Y∈@LPC''_{U,ν}$,
the induced map
$$ι^Z_{X,Y}:(@LPC''_{U,ν})^Z_{X,Y}→(@LPC'_{U,ν})^Z_{X,Y}$$
induces a weak equivalence on nerves.
Pick an arbitrary object $A:Z→@LPC'_{U,ν}$ in the codomain.
\ppar
We claim that the comma category $B=A/(@LPC''_{U,ν})^Z_{X,Y}$ is filtered,
therefore by Quillen's Theorem~A, the nerve of $ι^Z_{X,Y}$ is a weak equivalence.
Indeed, $B$ is nonempty: an object $A→A'$ in~$B$ can be constructed as follows.
For every object~$X$ in the zigzag~$A$,
construct its skeleton~$X'$ together with some inverse equivalences $X'→X$ and $X→X'$.
For a morphism $X→Y$ in the zigzag~$A$,
construct a morphism $X'→Y'$ as the composition $X'→X→Y→Y'$.
The resulting morphisms form a zigzag~$A'$ in $(@LPC''_{U,ν})^Z_{X,Y}$,
and the maps $X→X'$ provide a natural transformation $A→A'$ of zigzags.
\ppar
Next, if $A→A_1$ and $A→A_2$ are objects in~$B$, i.e., natural weak equivalences of $Z$-indexed zigzags,
then $A_1$ is weakly equivalent to $A_2$, and since both are skeletal by definition of $@LPC''_{U,ν}$,
$A_1$ is isomorphic to~$A_2$.
In particular, we have morphisms $A_1→A_2$ and $A_2→A_2$, showing that any two objects admit a pair of arrows to a third object.
\ppar
Finally, if $a_1:A→A_1$ and $a_2:A→A_2$ are objects in~$B$ and $f,g:A_1→A_2$ are a pair of parallel arrows in~$B$,
then $fa_1=ga_1=a_2$.
Since $a_1$ and $a_2$ are equivalences, we deduce that $f$ is naturally isomorphic to~$g$.
Since $A_1$ and $A_2$ are skeletal, we infer that $f=g$.

\subsection Relative categories

Consistent with our convention for categories,
we do not require relative categories or simplicial categories to be locally small.
In particular, in a simplicial category~$C$ the hom-object $C(X,X')$ for any objects $X,X'∈C$
can have a proper class $C(X,X')_n$ of $n$-simplices for any $n≥0$.
Thus, a simplicial category is a category enriched in simplicial classes.

The notion of a Dwyer–Kan equivalence of simplicial categories (Dwyer–Kan [\Func, §2.4])
continues to make sense for simplicial categories that are not locally small:
a simplicial functor $"F:C→D$ is a Dwyer–Kan equivalence if any object in~$D$ is homotopy equivalent to $"F(X)$ for some object~$X∈C$
and for any object $X,X'$, the induced map $C(X,X')→D("F(X),"F(X'))$ is a simplicial weak equivalence of simplicial classes.
The latter can be defined, for example, by adopting the statement of the simplicial Whitehead theorem for nonfibrant simplicial sets (^!{nonfibrant simplicial Whitehead}) as a definition.

The hammock localization construction of Dwyer–Kan [\Calc, §2.1] continues to make sense for relative categories that are not small or locally small.

\proclaim Remark.
^^={hammock construction}
The ^={hammock localization[|s]} of a relative category~$@C$
is a simplicial category~$\ham_@C$ with the same objects as~$@C$.
Given objects $X,Y∈@C$, the simplicial class $\ham_@C(X,Y)$
is constructed as the colimit
$$\colim_{Z∈\zz^\op} \N(@C^Z_{X,Y}).$$
Here $\zz^\op$ is Dwyer–Kan's indexing category~$\bf II$ [\Calc, §4.1].
(We prefer to work with the opposite category~$\zz$ since it is more directly related to categorical constructions;
the difference is analogous to how Segal's category~$Γ$ can be described by an ad hoc construction,
or as the opposite category of the category of finite pointed sets.)
We refer to the objects of~$\zz$ as {\it zigzag types}.
Objects of~$\zz$ are relative categories freely generated by finite sequences of morphisms like $←←→←→→←$,
where all left-pointing arrows $←$ are weak equivalences.
Morphisms of~$\zz$ are relative functors that preserve the natural ordering on objects and preserve the leftmost and rightmost objects.
Furthermore, $\N$ denotes the nerve functor and $@C^Z_{X,Y}$ is the category of relative functors $Z→@C$
that map the leftmost and rightmost objects of~$Z$ to $X$ and~$Y$ respectively.

\proclaim Remark.
^^={hammock construction hocolim}
By Dwyer–Kan [\Calc, Propositions 4.5, 5.4, 5.5], the colimit over~$Z$ in ^!{hammock construction} computes the homotopy colimit.

\proclaim Definition.
A ^={Dwyer–Kan equivalence[|s]} of relative categories is a relative functor whose hammock localization is a Dwyer–Kan equivalence of simplicial categories.

\proclaim Definition.
A relative functor $"F:C→D$ is ^={homotopically essentially surjective}
^^={homotopically essentially surjective functor[|s]}
if any object in~$D$ is weakly equivalent to an object in the image of~$"F$.
A relative functor $"F:C→D$ is ^={homotopically fully faithful}
^^={homotopically fully faithful functor[|s]}
if for any objects $X,X'∈C$ the induced simplicial map
$$\ham_C(X,X')→\ham_D("F(X),"F(X'))$$
is a simplicial weak equivalence.

\proclaim Proposition.
^^={criterion for relative equivalences}
A relative functor that is ^{homotopically essentially surjective} and ^{homotopically fully faithful}
is a ^{Dwyer–Kan equivalence}.

\proclaim Remark.
In the context of Dugger [\Pres, Definition~3.1],
homotopically surjective functors used in that definition
coincide with ^{homotopically essentially surjective functors}.

\proclaim Definition.
(Barwick–Kan [\RelCat, §3.3].)
A ^={homotopy equivalence[|s] of relative categories} is a relative functor $"F:C→D$
such that there is a relative functor $"G:D→C$
together with zigzags of natural weak equivalences $η:\id_C→"G∘"F$ and $ε:"F∘"G→\id_D$.

Every ^{homotopy equivalence of relative categories} is a ^{Dwyer–Kan equivalence}
because its hammock localization is a Dwyer–Kan equivalence of simplicial categories,
but the converse need not hold.

The literature on homotopy limits and colimits in relative categories is sparse.
While Dwyer–Hirschhorn–Kan–Smith [\DHKS, Chapter VIII] do provide an account
of homotopy limits and colimits in relative categories whose class of weak equivalences
satisfies the 2-out-of-6 property, it does not readily extend to all relative categories.

Instead, we apply the right Quillen equivalence from small relative categories to simplicial sets equipped with the Joyal model structure,
and then use the theory of limits and colimits in quasicategories.

\proclaim Notation.
^^={relative nerve notation}
Recall (Barwick–Kan [\RelCat, Corollary~6.11(i)]) that the model category of small relative categories
is Quillen equivalent to the Joyal model structure on simplicial sets
via a right Quillen equivalence, which we denote by $$\rN:@=RelCat→@sSet_@=Joyal.$$
The corresponding left Quillen equivalence will be denoted by $$\rK:@sSet_@Joyal→@RelCat.$$
The fibrant replacement functor on $@RelCat$ will be denoted by $$\frc: @RelCat → @RelCat.$$

\proclaim Definition.
^^={relative homotopy limit}
Suppose $D$ is a small category.
A small relative category~$C$ {\it admits $D$-indexed homotopy colimits\/}
if the small quasicategory $\rN\frc C$ admits $D$-indexed quasicategorical colimits.
A small relative category~$C$ {\it admits $D$-indexed homotopy limits\/}
if the small quasicategory $\rN\frc C$ admits $D$-indexed quasicategorical limits.

\proclaim Remark.
A ^{Dwyer–Kan equivalence} $"F:C→C'$ of relative categories induces
an equivalence $$\rN\frc("F):\rN\frc C→\rN\frc C'$$ of quasicategories.
Thus, $C$ admits $D$-indexed homotopy (co)limits if and only if $C'$ does.

\proclaim Definition.
^^={preservation of relative homotopy limits}
Suppose $D$ is a small category
and $C$ and $C'$ are small relative categories that admit $D$-indexed homotopy colimits (respectively limits).
A relative functor $"F:C→C'$
{\it preserves $D$-indexed homotopy colimits\/}
if the functor $$\rN\frc("F):\rN\frc(C)→\rN\frc(C')$$
preserves $D$-indexed quasicategorical colimits.
A relative functor $"F:C→C'$
{\it preserves $D$-indexed homotopy limits\/}
if the functor $\rN\frc("F)$
preserves $D$-indexed quasicategorical limits.

\subsection Simplicial Whitehead theorem

\let\ssim\Delta
\def\sbdr{\partial\ssim}

\def\dom{\mathop{\rm dom}}
\def\codom{\mathop{\rm codom}}

\def\diconv{
	\def\sph{\sbdr^n}
	\def\odisk{\ssim^n}
	\def\ndisk{\odisk}
	\def\disks{\odisk\sqcup_{\sph}\ndisk}
	\def\relh{\odisk\times\ssim^1\sqcup_{\sph\times\ssim^1}\sph}
}


\def\whitehead#1{\sqcd{
	\sph&\mapright{}&\ndisk\cr
	\mapdown{}&&\mapdown{}\cr
	\odisk&\mapright{}&\relh#1\cr
}}

\def\smap{\sSet^\to}
\def\sSet{{\sf sSet}}
\def\Map{\mathop{\rm Map}}
\def\RMap{{\bf R}\Map}

\proclaim Definition.
^^={simplicial map category}
Denote by $\smap$ the category of functors $\{0\to1\}\to\sSet$.
Objects are simplicial maps (depicted vertically)
and morphisms are commutative squares,
where the two vertical maps are the source and target.
Equip $\smap$ with the projective model structure.

\proclaim Remark.
In the model category $\smap$ (^!{simplicial map category}),
projectively cofibrant objects are simplicial maps that are cofibrations.
Projective cofibrations are commutative squares where the top map
and pushout product of left and top maps is a cofibration of simplicial sets.
Fibrant objects are simplicial maps whose domain and codomain are Kan complexes.

\proclaim Proposition.
^^={surjections detected in homotopy categories}
Fix a simplicial model category~$M$, such as $\smap$ (^!{simplicial map category}).
Suppose $\alpha$ is a cofibration between cofibrant objects in~$M$
and $\Omega$ is a fibrant object in~$M$.
The map of sets $\hom(\alpha,\Omega)$ is surjective if and only if
the map of sets $\pi_0\RMap(\alpha,\Omega)$ is surjective.
Here $\hom$ denotes mapping sets in~$M$
and $\RMap$ denotes derived mapping simplicial sets in~$M$.

\proof Proof.
Since $\alpha$ is a cofibration between cofibrant objects and the object $\Omega$ is fibrant,
the simplicial map $\Map(\alpha,\Omega)$ is a fibration between fibrant objects in simplicial sets.
A fibration of simplicial sets is surjective on 0-simplices if and only if the induced map on $\pi_0$ is a surjection.
Thus, the map of sets $\hom(\alpha,\Omega)$ is surjective if and only if
the map of sets $\pi_0\Map(\alpha,\Omega)$ is surjective.
The latter map is isomorphic to $\pi_0\RMap(\alpha,\Omega)$ because $\alpha$ is a cofibration between cofibrant objects and $Ω$ is fibrant.

\proclaim Proposition.
^^={Whitehead is model-independent}
Fix a simplicial model category~$M$, such as $\smap$ (^!{simplicial map category}).
Suppose $\alpha$ and $\beta$ are weakly equivalent cofibrations between cofibrant objects in~$M$
and $\Omega$ is a fibrant object in~$M$.
Then the map of sets
$\hom(\alpha,\Omega)$
is a surjection of sets
if and only if
$\hom(\beta,\Omega)$
is a surjection of sets.
Here $\hom$ denotes mapping sets in~$M$.

\proof Proof.
By ^!{surjections detected in homotopy categories}, the map $\hom(α,Ω)$ is surjective if and only if $\pi_0\RMap(α,Ω)$ is surjective.
Likewise, the map $\hom(β,Ω)$ is surjective if and only if $\pi_0\RMap(β,Ω)$ is surjective.
Since $\alpha$ is weakly equivalent to $\beta$, the map of sets $\pi_0\RMap(\alpha,\Omega)$ is isomorphic to the map of sets $\pi_0\RMap(\beta,\Omega)$,
which proves the lemma.

\proclaim Definition.
^^={Whitehead objects}
Denote by $\lambda$
the projective cofibration between projectively cofibrant objects in $\smap$
given by the commutative square on the right of the following diagram:
$$\whitehead{}\simeq\diconv\whitehead.$$
In what follows, $\iota$ refers to any projective cofibration with projectively cofibrant source in $\smap$ that is weakly equivalent to~$\lambda$,
as depicted on the left.
Here $\sph$ means ``sphere'', $\odisk$ means ``old disk'', $\ndisk$ means ``new disk'', $\relh$ means ``relative homotopy''.
The idea is that $\relh$ expresses a relative homotopy from the old disk $\odisk$ to the new disk $\ndisk$ relative boundary $\sph$, the sphere.
We also set $\disks=\odisk\sqcup_{\sph}\ndisk$, both disks combined, which is the boundary of~$\relh$.

The following lemma reformulates a criterion
due to Kan [\Kan, Theorem~7.2], originally due to Whitehead [\CHi, Theorem~1] in the case of topological spaces.

\proclaim Proposition.
^^={simplicial Whitehead}
(Dugger–Isaksen [\WESP, Proposition~4.1].)
A simplicial map~$p$ between Kan complexes is a weak equivalence
if and only if the map of sets
$$\hom(\lambda,p): \hom(\codom\lambda,p)\to\hom(\dom\lambda,p)$$
is a surjection.
Here $\hom$ denotes mapping sets in the category $\smap$.

The following corollary combines ^!{simplicial Whitehead} and ^!{Whitehead is model-independent}.

\proclaim Corollary.
^^={model-independent simplicial Whitehead}
Suppose $\iota$ is a projective cofibration between projectively cofibrant objects in $\smap$
that is weakly equivalent to the map~$\lambda$ in ^!{Whitehead objects}.
Then a simplicial map~$p$ between Kan complexes is a weak equivalence
if and only if the map of sets
$$\hom(\iota,p): \hom(\codom\iota,p)\to\hom(\dom\iota,p)$$
is a surjection of sets.
Here $\hom$ denotes mapping sets in the category $\smap$.

\proclaim Remark.
^^={unfolded model-independent simplicial Whitehead}
Expanding the statement of ^!{model-independent simplicial Whitehead},
a map~$p: A\to B$ of Kan complexes is a weak equivalence if and only if for any commutative square
$$\sqcd{
\sph&\mapright a&A\cr
\mapdown{}&&\mapdown p\cr
\odisk&\mapright b&B\cr
}$$
we can find maps $d:\ndisk\to A$
and $e:\relh\to B$
that make the following diagram commute:
$$\diagram[2em]A:4,4=$A$,B:4,0=$B$\rlap.,S:0,4=$\sph$,O:0,0=$\odisk$,N:2,3=$\ndisk$,R:2,1=$\relh$,AB.rt=$p$,SA.top=$a$,OB.bot=$b$,SO.lft=,SN.lft=,NA.lrt=$d$,OR.lft=,RB.urt=$e$,NR.lft=;$$

\proclaim Remark.
^^={nonfibrant simplicial Whitehead}
One way to expand ^!{model-independent simplicial Whitehead} (and ^!{unfolded model-independent simplicial Whitehead}) to the case when $A$ or $B$ is not a Kan complex
is to replace $f$ with the map $\Exi f$, where $\Exi$ denotes Kan's fibrant replacement functor for simplicial sets.
If the simplicial sets in the commutative square $ι$ of ^!{Whitehead objects} are compact (i.e., have finitely many nondegenerate simplices),
then the maps to $\Exi f$ will factor through $"Ex^k f$ for some $k≥0$,
which allows us to use the adjunction $"Sd^k⊣"Ex^k$ to keep $f$ intact.
See ^!{specialized simplicial Whitehead} for an example.
Somewhat more generally, we can formulate the following criterion that is independent of specific choices of models for spheres and disks.
A simplicial map~$p$ between simplicial sets is a weak equivalence
if and only if
for any simplicial map $σ:\sph→\odisk$ weakly equivalent to the inclusion $κ:∂Δ^n→Δ^n$
and a morphism $ψ:σ→p$ in $\smap$,
we can factor $ψ$ as the composition $ψ=χι$,
where $ι:σ→τ$ is some morphism in $\smap$ satisfying the conditions of ^!{Whitehead objects} and $χ:τ→p$ is some other morphism in $\smap$.
In fact, it suffices to examine a single representative $(σ,ψ)$ for every element in the set of morphisms $κ→p$ in the homotopy category of $\smap$.

\proclaim Corollary.
^^={specialized simplicial Whitehead}
Denote by~$Λ$ the simplicial subset of $Δ^2$ generated by the 1-simplices $0→2$ and $1→2$.
A simplicial map $p:A→B$ is a simplicial weak equivalence whenever for any $k≥0$, $n≥0$,
and a commutative square
$$\sqcd{
"Sd^k ∂Δ^n&\mapright{α}&A\cr
\mapdown{}&&\mapdown{p}\cr
"Sd^k Δ^n&\mapright{β}&B,\cr
}$$
we can construct maps
$$γ:"Sd^k Δ^n→A,\qquad Γ:Δ^1⨯"Sd^k Δ^n→B,$$
$$π:Λ⨯"Sd^k ∂Δ^n→A,\qquad Π:Δ^2⨯"Sd^k ∂Δ^n→B$$
such that
the map $Γ$ is a simplicial homotopy from $β$ to $p∘γ$,
the map $Π$ restricts to $p∘π$ on $Λ⨯"Sd^k ∂Δ^n$,
the map $π$ restricts to $α$ on $0⨯"Sd^k ∂Δ^n$,
the restrictions of $π$ to $1⨯"Sd^k ∂Δ^n$ and $γ$ to $"Sd^k ∂Δ^n$ coincide,
and the restrictions of $Π$ to $(0→1)⨯"Sd^k ∂Δ^n$ and $Γ$ to $Δ^1⨯"Sd^k ∂Δ^n$ coincide.

\proof Proof.
Specialize ^!{nonfibrant simplicial Whitehead} to the case when $$σ:"Sd^k ∂Δ^n→"Sd^k Δ^n$$
and $$τ:Λ⨯"Sd^k ∂Δ^n ⊔_{"Sd^k ∂Δ^n} "Sd^k Δ^n → Δ^2⨯"Sd^k ∂Δ^n ⊔_{Δ^1⨯"Sd^k ∂Δ^n} Δ^1⨯"Sd^k Δ^n$$
are the canonical inclusions.
(See the figure in ^!{specialized simplicial Whitehead figure} to see how different pieces stick together.)

\proclaim Remark.
^^={specialized simplicial Whitehead figure}
We illustrate ^!{specialized simplicial Whitehead} with the following diagrams, where the left diagram depicts $A$ and the right diagram depicts~$B$.
We depict only a single radius connecting a point $α$ (respectively $p∘α$) on the sphere $"Sd^k ∂Δ^n$ to the center of $"Sd^k Δ^n$, represented by $⋯$:
$$\diagram[3em]A:0,0=$α$,C:-2,-1=$π_2$,D:0,-2=$π_1$,E:2,-2=$⋯$,DE.bot=$γ$,AC.ulft@0.4=$π_{02}$,DC.llft@0.4=$π_{12}$,G:-0.7,-1=$π$;
\qquad\diagram[3em]A:0,0=$p∘α$,B:2,0=$⋯$,C:-2,-1=$p∘π_2$,D:0,-2=$p∘π_1$,E:2,-2=$⋯$,AB.top=$β$,DE.bot=$p∘γ$,AC.ulft@0.4=$p∘π_{02}$,DC.llft@0.4=$p∘π_{12}$,AD.lft=,F:1,-1=$Γ$,G:-0.7,-1=$Π$;$$
Thus, the left diagram depicts a sphere (represented by the single vertex~$α$)
being filled by a disk (represented by the bottom chain of morphisms going to $⋯$),
whereas the right diagram takes the image of the left diagram under~$p$,
and then homotopes it relative boundary to the map~$β$,
using the indicated triangle~$Π$ together with a finite collection of squares that look like~$Γ$.

\section Combinatorial model categories

See Beke [\SHM], Barwick [\LR], Lurie [\HTT, Appendix~A], Low [\Heart] for the background on combinatorial model categories.
We recall some basic definitions to fix terminology.

\proclaim Definition.
A ^={model structure[|s]} on a category~$C$ is a pair of weak factorization systems $(@=C,@=AF)$, $(@=AC,@=F)$
such that the class $@=W=@AF∘@AC$ satisfies the 2-out-of-3 property.
A ^={model categor[y|ies]} is a model structure on a category that admits finite limits and finite colimits.
A weak factorization system $(@A,@B)$ is ^{cofibrantly generated} if there is a set of morphisms~$I$ in~$@A$ such that the class of morphisms in~$C$ with a right lifting property with respect to every element of~$I$ coincides with the class~$@B$.
A model structure is ^={cofibrantly generated} if both weak factorization systems $(@C,@AF)$, $(@AC,@F)$ are cofibrantly generated.
A ^={combinatorial model categor[y|ies]}
is a cofibrantly generated model structure on a locally presentable category.
A ^={left Quillen functor} between model categories is a
left adjoint functor that preserves
elements of $@C$ (cofibrations) and $@AC$ (acyclic cofibrations).

We now review some of the more specialized definitions from Low [\Heart].
Similar notions can be found in Chorny–Rosický [\ClassComb].
The appearance of the sharp ordering $κ⊲λ$ in this section is dictated by ^!{accessibly generated characterization}.

\proclaim Definition.
(Low [\Heart, Definition~5.11].)
Given regular cardinals $κ⊲λ$, a $(κ,λ)$-^={miniature model categor[y|ies]} is a model category~$M$
such that there exist $λ$-small sets of morphisms in $\K_κ^λ(M)$ (^!{presentable object}) that cofibrantly generate the model structure of~$M$
and the underlying category of~$M$ satisfies the following properties:
\li $M$ is a $(κ,λ)$-^{accessibly generated category} (^!{accessibly generated category});
\li $M$ has finite limits and $λ$-small colimits;
\li Hom-sets in $\K_κ^λ(M)$ are $λ$-small.

\proclaim Definition.
(Low [\Heart, Definition~5.1].)
Given regular cardinals $κ⊲λ$, a
^={strongly $(κ,λ)$-combinatorial model categor[y|ies]}
^^={strongly combinatorial model categor[y|ies]}
is a combinatorial model category~$M$
such that there exist $λ$-small sets of morphisms in $\K_κ(M)$ that cofibrantly generate the model structure of~$M$
and the underlying category of~$M$ satisfies the following properties:
\li $M$ is a locally $κ$-presentable category;
\li $\K_λ(M)$ is closed under finite limits in~$M$;
\li Hom-sets in $\K_κ(M)$ are $λ$-small.

\proclaim Proposition.
^^={miniature strongly CMC}
(Low [\Heart, Proposition~5.12, Theorem~5.14(i)].)
The functor $\K_λ$ and functor $\Ind^λ$ establish a correspondence between
$(κ,λ)$-^{miniature model categories} and ^{strongly $(κ,λ)$-combinatorial model categories}.
In particular, we have equivalences of model categories $C→\Ind^λ(\K_λ(C))$ and $D→\K_λ(\Ind^λ(D))$ that preserve and reflect weak equivalences, fibrations, and cofibrations.
This correspondence preserves left Quillen equivalences.
Furthermore, any combinatorial model category is a strongly $(κ,λ)$-combinatorial model category for some regular cardinals $κ⊲λ$ (Low [\Heart, Proposition~5.6]).
Any strongly $(κ,λ)$-combinatorial model category is a strongly $(κ,μ)$-combinatorial model category for any $μ⊳λ$ (Low [\Heart, Remark~5.2]).

\proclaim Definition.
Suppose $"F:C→D$ is a left Quillen functor between combinatorial model categories.
If $λ$ is a regular cardinal,
we say that $"F$ is a ^={left $λ$-Quillen functor[|s]} if $"F$ is strongly $λ$-accessible (^!{strongly accessible})
and $C$ and $D$ are ^{strongly $(κ,λ)$-combinatorial model categories} (^!{strongly combinatorial model category}) for some regular cardinal~$κ$.

The following proposition follows from ^!{accessible functors are strong}.

\proclaim Proposition.
^^={left Quillen functors are $λ$-Quillen}
(Low [\Heart, Proposition~5.6, Lemma~2.5].)
For any left Quillen functor $"F:C→D$ between combinatorial model categories
and any regular cardinal~$κ$,
there is a regular cardinal $λ⊵κ$ such that
$"F$ is a left $λ$-Quillen functor (^!{left $λ$-Quillen functor}).

The relative category $@CMC$ can be informally described as follows.
Objects are combinatorial model categories.
Morphisms are left Quillen functors.
Weak equivalences are left Quillen equivalences.
To avoid size issues, we follow ^!{size aspects}.

\proclaim Definition.
^^={CMC}
The relative category $@=CMC$ is defined as follows.
Objects are pairs $(λ,C)$, where $λ$ is a regular cardinal and $C$ is a small $(κ,λ)$-^{miniature model category}
(Low [\Heart, Definition~5.11]), where $κ$ is some regular cardinal such that $κ⊲λ$.
Morphisms $(λ,C)→(μ,D)$ exist if $λ⊴μ$, in which case they are left Quillen functors $C→D$.
Weak equivalences are generated as a subcategory
by morphisms $(λ,C)→(μ,D)$ for which $λ=μ$ and $C→D$ is a left Quillen equivalence,
together with morphisms $(λ,C)→(μ,D)$
for which the left Quillen functor $C→D$ exhibits $D$ as the $(λ,μ)$-ind-completion of~$C$,
i.e., the small model category of $μ$-presentable objects (Low [\Heart, Proposition~5.12]) in the $λ$-ind-completion of~$C$ (Low [\Heart, Theorem~5.14]).

While we do not assume any large cardinal axioms for the main results of this paper,
we can ask whether in presence of a strongly inaccessible cardinal
our definition of~$@CMC$ is equivalent to the more obvious definition of~$@CMC$ that uses universes.
This is answered in the affirmative by the following definitions and proposition.

\proclaim Definition.
^^={smaller CMCU}
Given a strongly inaccessible cardinal~$U$,
the relative category $@CMC_U$ is defined as the full subcategory of~$@CMC$ (^!{CMC})
on objects $(λ,M)$ such that after discarding the model structure we have $(λ,M)∈@LPC_U$ (^!{smaller LPCU}).

\proclaim Definition.
^^={universe CMCU}
Suppose $U$ is a strongly inaccessible cardinal.
A {\it $U$-combinatorial model category\/} is a category in $@LPC'_U$ (^!{universe LPCU}) equipped with a model structure
that is cofibrantly generated by a $U$-small set of morphisms.
We define a relative category $@CMC'_U$
as the relative category of $U$-combinatorial model categories, left Quillen functors, and left Quillen equivalences.

\proclaim Proposition.
^^={smaller universe CMCU}
There is a ^{Dwyer–Kan equivalence} of relative categories
$$@CMC_U→@CMC'_U$$
of ^!{smaller CMCU} and ^!{universe CMCU}.

\proof Proof.
The functor $$@CMC_U→@CMC'_U$$
is constructed by promoting the functor $$\Ind_U:@LPC_U→@LPC'_U$$ (^!{smaller universe LPCU})
to a functor $$\Ind_U:@CMC_U→@CMC'_U,$$
as described in Low [\Heart, Theorem~5.14].
(Low's construction works with large $U$-ind-completions; we pass to the full subcategory of $U$-presentable objects to obtain the version stated above.)
By definition of~$@CMC$ (^!{CMC}), this functor preserves weak equivalences,
so it is a relative functor.
\ppar
In complete analogy to ^!{smaller universe LPCU},
we introduce filtrations on $@CMC_U$ and $@CMC'_U$ (indexed by a regular cardinal~$ν$)
that are respected by the functor $\Ind_U$
and show that $\Ind_U$ induces a ^{homotopy equivalence of relative categories} for every step of the filtration.
\ppar
Fix a regular cardinal~$ν$.
Define $@CMC_{U,ν}$ as the full subcategory of $@CMC_U$ consisting of objects $(λ,C)$ for which $λ≤ν$.
Define $@CMC'_{U,ν}$ as the full subcategory of $@CMC'_U$ on objects~$M$
such that $\Ind^U(M)$ (^!{smaller ind-completion}) is a strongly $(κ,ν)$-combinatorial model category for some regular cardinal $κ⊲ν$
(^!{strongly combinatorial model category})
and morphisms given by left Quillen functors~$"F$ such that $\Ind^U("F)$ is a $ν$-Quillen functor (^!{left $λ$-Quillen functor}).
By construction and Low [\Heart, Theorem~5.14], the functor $\Ind_U$ restricts to a functor $$\Ind_{U,ν}:@CMC_{U,ν}→@CMC'_{U,ν}.$$
We now show that $\Ind_{U,ν}$ is a ^{homotopy equivalence of relative categories}.
\ppar
The inverse functor is $$\K_{U,ν}:@CMC'_{U,ν}→@CMC_{U,ν},$$
which is well-defined by Low [\Heart, Proposition~5.12].
The functor $\K_{U,ν}$ sends
an object $C∈@CMC'_{U,ν}$ to the object $(ν,\K_ν^U(C))∈@CMC_{U,ν}$
and a functor $"F:C→D$ in $@CMC'_{U,ν}$ to the restriction $$(ν,\K_ν^U(C))→(ν,\K_ν^U(D)),$$
which is well-defined because $\Ind^U("F)$ is a left $ν$-Quillen functor.
\ppar
The natural weak equivalences $$\id_{@CMC_{U,ν}}→\K_{U,ν}∘\Ind_{U,ν}, \qquad \Ind_{U,ν}∘\K_{U,ν}→\id_{@CMC'_{U,ν}}$$
are inherited from the proof of ^!{smaller universe LPCU} and are weak equivalences
because their underlying functors are equivalences of categories
and the model structures coincide by Low [\Heart, Proposition~5.12 and Theorem~5.14].

\section Combinatorial relative categories

\proclaim Definition.
^^={CRC}
The relative category $@=CRC$ is defined as follows.
Objects are pairs $(λ,C)$, where $λ$ is a regular cardinal and $C$ is a small relative category that admits $λ$-small homotopy colimits (^!{relative homotopy limit}).
Morphisms $(λ,C)→(μ,D)$ exist if $λ≤μ$, in which case they are relative functors $C→D$ that preserve $λ$-small homotopy colimits (^!{preservation of relative homotopy limits}).
Weak equivalences $(λ,C)→(μ,D)$ are generated as a subcategory
by morphisms $(λ,C)→(μ,D)$ for which $λ=μ$ and $C→D$ is a ^{Dwyer–Kan equivalence},
together with morphisms $(λ,C)→(μ,D)$ for which the functor $"F:C→D$ exhibits $D$ as the $(λ,μ)$-ind-completion of~$C$,
namely,
the functor $\rN\frc("F)$ exhibits $\rN\frc D$ as the quasicategory of $μ$-presentable objects in the $λ$-ind-completion of the quasicategory~$\rN\frc C$
(see ^!{PrL} and ^!{relative nerve notation}).

\proclaim Remark.
In ^!{CRC}, we created the class of weak equivalences in $@CRC$
using the functor
$$\RN:@CRC→@PrL$$
of ^!{derived relative nerve}, where the weak equivalences of $@PrL$ are given in ^!{PrL}.
A more natural way to introduce weak equivalences in $@CRC$
is to define homotopy $μ$-presentable objects and homotopy $λ$-ind-completions of relative categories directly,
without referring to quasicategories.
Such an approach would produce exactly the same class of weak equivalences.
However, it would require us to introduce all the relevant definitions
and show their compatibility with analogous quasicategorical definitions,
further adding to the length of this article,
whereas in the quasicategorical context the necessary results are already available in Cisinski [\HCHA, Chapter~7].
Thus, we bypass the issue by transferring weak equivalences from $@PrL$.

\proclaim Definition.
^^={smaller CRCU}
Given a strongly inaccessible cardinal~$U$,
the relative category $@CRC_U$ is defined as the full subcategory of~$@CRC$ (^!{CRC})
on objects $(λ,C)$, where $λ<U$ and $C$ is $U$-essentially $U$-small (^!{size notions}).

\proclaim Definition.
^^={universe CRCU}
Given a strongly inaccessible cardinal~$U$,
the relative category $@CRC'_U$
is the relative category of
locally $U$-small relative categories~$C$ (^!{locally $U$-small category}) such that $\rN\frc C∈@PrL'_U$ (^!{universe PrLU}),
with morphisms given by relative functors~$"F$ such that $\rN\frc("F)$ is a morphism in $@PrL'_U$
and weak equivalences being ^{Dwyer–Kan equivalences} of relative categories.

\proclaim Proposition.
^^={smaller universe CRCU}
^^={smaller RInd}
The functor $$"=RInd_U="Reedy∘"MInd_U:@CRC_U→@CRC'_U,$$
where $"Reedy$ is as in ^!{Reedy} and $"MInd_U$ is as in ^!{smaller MInd}
is a ^{Dwyer–Kan equivalence} of relative categories.

\proof Proof.
The proof is very similar to the proofs of ^!{smaller universe LPCU} and ^!{smaller universe CMCU}.
We only briefly indicate the necessary modifications.
The functor $"RInd_U$ indeed lands in $@CRC'_U$ by ^!{Reedy MInd equivalence} (taking $ν=U$ and $C=@CRC_U$ there).
\ppar
Given a regular cardinal~$ν$, we take $@CRC_{U,ν}$ to be the full subcategory of $@CRC_U$ consisting of objects $(λ,C)$ with $λ≤ν$
and $@CRC'_{U,ν}$ to be the subcategory of $@CRC'_U$ consisting of objects $C$ such that $\Ind^U(\rN\frc C)$ is a $ν$-presentable quasicategory
and morphisms $"F$ such that $\Ind^U(\rN\frc "F)$ is a strongly $ν$-accessible left adjoint functor of quasicategories.
\ppar
The homotopy inverse to $"RInd_U$
is given by the functor $"=RK_ν^U$
that sends $C∈@CRC'_U$ to $(ν,C_ν)$, where $C_ν$ is the full subcategory of~$C$
on objects whose images in $\rN\frc C$
are $(ν,U)$-presentable objects (in the quasicategorical sense).
By ^!{Reedy MInd equivalence}, we have a natural weak equivalence $(κ,E)→(ν,"RK_ν^U("RInd_U(κ,E)))$ (where $(κ,E)∈@CRC_U$) that sends $e∈E$ to $n↦(Δ^n⊗"Y(e))$, where $"Y$ is the Yoneda embedding.
We also have a natural weak equivalence $C→"RInd_U(ν,"RK_ν^U(C))$ for $C∈@CRC'_U$,
which sends $X∈C$ to $n↦Δ^n⊗(A↦\ham_C(A,X))$.

\section Presentable quasicategories

Recall that a quasicategory is {\it presentable\/} if it is accessible (Lurie [\HTT, Definition~5.4.2.1]) and admits small colimits (Joyal [\QCKC, Definition~4.5]).
The relative category $@PrL$ can be informally described as the relative category of presentable quasicategories, left adjoint functors, and equivalences.
To avoid size issues, we follow ^!{size aspects}.

\proclaim Definition.
^^={PrL}
The relative category $@=PrL$ is defined as follows.
Objects are pairs $(λ,C)$, where $λ$ is a regular cardinal and $C$ is a small quasicategory that admits $λ$-small colimits.
Morphisms $(λ,C)→(μ,D)$ exist if $λ≤μ$, in which case they are functors $C→D$ that preserve $λ$-small colimits.
Weak equivalences are morphisms $(λ,C)→(μ,D)$ such that $C→D$ exhibits $D$ as the $(λ,μ)$-ind-completion of~$C$,
i.e., the quasicategory of $μ$-presentable objects (Lurie [\HTT, Definition~5.3.4.5]) in the $λ$-ind-completion of~$C$ (Lurie [\HTT, Definition~5.3.5.1]).

\proclaim Definition.
^^={smaller PrLU}
Given a strongly inaccessible cardinal~$U$,
the relative category $@PrL_U$ is defined as the full subcategory of~$@PrL$ (^!{PrL})
on objects $(λ,C)$, where $λ<U$ and $C$ is $U$-essentially $U$-small (^!{size notions}).

\proclaim Definition.
^^={universe PrLU}
Given a strongly inaccessible cardinal~$U$,
the relative category $@PrL'_U$
is the relative category of $U$-locally $U$-presentable quasicategories (^!{size notions}),
left adjoint functors, and equivalences of quasicategories, i.e., weak equivalences in the Joyal model structure.

\proclaim Proposition.
^^={smaller universe PrLU}
^^={smaller QInd}
The functor $$"=QInd_U=\rN∘\frc∘"Reedy∘"MInd_U∘\rK:@PrL_U→@PrL'_U,$$
where $\rK$, $\rN$, and $\frc$ are as in ^!{relative nerve notation}, $"Reedy$ is as in ^!{Reedy}, and $"MInd_U$ is as in ^!{smaller MInd},
is a ^{Dwyer–Kan equivalence} of relative categories $@PrL_U$ (^!{smaller PrLU}) and $@PrL'_U$ (^!{universe PrLU}).

\proof Proof.
The proof is very similar to the proofs of ^!{smaller universe LPCU}, ^!{smaller universe CMCU}, and ^!{smaller universe CRCU}.
We only briefly indicate the necessary modifications.
The somewhat cumbersome and roundabout definition of $"QInd_U$
is explained by the fact that we need a (strict) relative functor,
whereas the familiar quasicategorical constructions of ind-completions only provide homotopy coherent functors.
\ppar
Given a regular cardinal~$ν$, we take $@PrL_{U,ν}$ to be the full subcategory of $@PrL_U$ consisting of objects $(λ,C)$ with $λ≤ν$
and $@PrL'_{U,ν}$ to be the subcategory of $@PrL'_U$ on objects~$C$ such that $\Ind^U(C)$ is a $ν$-presentable quasicategory
and morphisms~$"F$ such that $\Ind^U("F)$ is a strongly $ν$-accessible left adjoint functor of quasicategories.
\ppar
The homotopy inverse to $"QInd_U$
is given by the functor $"=QK_ν^U$
that sends $C∈@PrL'_U$ to $(ν,C_ν)$, where $C_ν$ is the full subcategory of~$C$ on $(ν,U)$-presentable objects (in the quasicategorical sense).
We have a natural weak equivalence $(κ,E)→(ν,"QK_ν^U("QInd_U(κ,E)))$ (where $(κ,E)∈@PrL_U$) given by the quasicategorical variant of the Yoneda embedding.
We also have a natural weak equivalence $C→"QInd_U(ν,"QK_ν^U(C))$ for $C∈@PrL'_U$,
given by the quasicategorical variant of the restricted Yoneda embedding.

\proclaim Remark.
We can turn $@PrL'_U$ into a simplicial category $@PrL'_{U,Δ}$ by declaring the hom-object $@PrL'_U(A,B)$
to be the simplicial subset of the maximal Kan subcomplex in the mapping simplicial set~$B^A$,
comprising connected components of left adjoint functors.
The homotopy coherent nerve $\N@PrL'_{U,Δ}$ of this simplicial category is precisely the quasicategory $@PrL$ constructed by Lurie [\HTT, Definition~5.5.3.1].
The canonical functor $\N@PrL'_U→\N@PrL'_{U,Δ}$
descends to a functor $\N@PrL'_U[W^{-1}]→\N@PrL'_{U,Δ}$,
which is an equivalence of quasicategories.


\section From combinatorial model categories to combinatorial relative categories

In this section we define two weakly equivalent ^{Dwyer–Kan equivalences} $@CMC→@CRC$.
The first equivalence, $"Cof$, is defined in a straightforward way by restricting to full subcategories of cofibrant objects.
The second equivalence, $"Reedy$, is defined by taking the relative category of cosimplicial resolutions,
i.e., Reedy cofibrant cosimplicial objects whose cosimplicial structure maps are weak equivalences.
This enables us to construct left Quillen equivalences from simplicial presheaves on such categories of diagrams
to the original model category, replicating a construction of Dugger [\Pres].

\proclaim Definition.
^^={Cof}
The relative functor $$"=Cof:@CMC→@CRC$$ is defined as follows.
An object $(λ,C)∈@CMC$ is sent to $(λ,"cof(C))$, where $"=cof(C)$ is the relative category of cofibrant objects in~$C$ with induced weak equivalences.
A morphism $(λ,C)→(μ,D)$ given by a left Quillen functor $"F:C→D$
is sent to the morphism $(λ,"cof(C))→(μ,"cof(D))$ given by the restriction and corestriction of~$"F$.

\proclaim Proposition.
^^={Cof is correct}
^!{Cof} is correct.

\proof Proof.
Given an object $(λ,M)∈@CMC$, we have to show that $(λ,"cof(M))∈@CRC$.
That is, if $M$ is a small $(κ,λ)$-^{miniature model category}, we have to show that the small relative category $"cof(M)$ admits $λ$-small homotopy colimits.
By ^!{relative homotopy limit}, this means that
the small quasicategory $\rN\frc "cof(M)$ admits $λ$-small colimits.
Since $"cof(M)→M$ is a ^{Dwyer–Kan equivalence},
it suffices to show that the small quasicategory $\rN \frc M$
admits $λ$-small colimits.
The small quasicategory $\rN\frc M$ is a localization of~$M$ with respect to its class of weak equivalences
in the sense of Cisinski [\HCHA, Definition~7.1.2],
denoted by $L(M)$ there.
By Cisinski [\HCHA, Remark~7.9.10],
the small quasicategory $L(M)≃\rN \frc M$ admits $λ$-small colimits.
\ppar
Given a morphism $(λ,M)→(μ,N)$ in~$@CMC$, we have to show that the functor $\rN\frc M→\rN\frc N$ preserves $λ$-small colimits.
The latter functor is an induced functor between localizations of $M$ and~$N$
in the sense of Cisinski [\HCHA, Definition~7.1.2], denoted by $L(M)→L(N)$ there.
By Cisinski [\HCHA, Remark~7.9.10],
the map $L(M)→L(N)$ preserves $λ$-small colimits.
\ppar
The functor $"Cof$ preserves weak equivalences in $@CMC$.
Indeed, the latter are generated by left Quillen equivalences and ind-completions.
$"Cof$ maps left Quillen equivalences $(λ,M)→(λ,N)$ to ^{homotopy equivalences of relative categories} (^!{homotopy equivalence of relative categories});
the homotopy inverse is given by the right derived functor of the right adjoint, composed with a cofibrant replacement functor.
$"Cof$ maps a morphism $(λ,M)→(μ,\Ind^λ_μ(M))$
to the morphism
$$(λ,"cof(M))→(μ,"cof(\Ind^λ_μ(M))),$$
which is weakly equivalent to the morphism
$$(λ,M)→(μ,\Ind^λ_μ(M)).$$
Taking $N=\Ind^λ(M)$, we can identify the latter morphism with
$$(λ,\K_λ(N))→(μ,\K_μ(N)).$$
The left Quillen functor $\K_λ(N)→\K_μ(N)$ preserves weak equivalences.
Furthermore, in the model category~$N$ the $λ$-filtered or $μ$-filtered colimits are also homotopy colimits.
Thus, $\K_λ(N)$ respectively $\K_μ(N)$ comprise the homotopy $λ$-presentable respectively homotopy $μ$-presentable objects in~$N$.
Hence, applying the functor $\rN\frc "cof$ (equivalently, $\rN \frc$ or simply $\rN$) yields a functor of quasicategories
that is equivalent to the inclusion $\K_λ(\rN N)→\K_μ(\rN N)$.

Although the functor $"Cof$ is a Dwyer–Kan equivalence,
it is not quite sufficient for our purposes,
and we have to introduce a weakly equivalent functor $"Reedy$.
In particular, the proof of crucial ^!{MInd Reedy equivalence} does not work with $"Cof$ instead of $"Reedy$,
as explained in ^!{Reedy Cof difference}.

\proclaim Definition.
^^={Reedy}
The relative functor $$"=Reedy:@CMC→@CRC$$
is defined as follows.
An object $(λ,C)$ is sent to
the pair $(λ,"Reedy(C))$,
where $"Reedy(C)$ is the small relative category of cosimplicial resolutions in~$C$, i.e.,
Reedy cofibrant cosimplicial objects in~$C$ whose cosimplicial structure maps are weak equivalences.
We equip $"Reedy(C)$ with degreewise weak equivalences.
A morphism $(λ,C)→(μ,D)$ is sent to the morphism $$(λ,"Reedy(C))→(μ,"Reedy(D))$$
given by the relative functor $"Reedy(C)→"Reedy(D)$
itself induced by the left Quillen functor $C→D$.
\ppar
The natural weak equivalence $$"=ev_0:"Reedy→"Cof$$ ($"Cof$ was introduced in ^!{Cof}) sends an object $(λ,C)∈@CMC$
to the morphism $$"ev_0(λ,C):(λ,"Reedy(C))→(λ,"Cof(C))$$ induced by the relative functor
$$"Reedy(C)→"Cof(C)$$ that evaluates a Reedy cofibrant cosimplicial diagram at the simplex $[0]∈Δ$.

\proclaim Proposition.
^!{Reedy} is correct.

\proof Proof.
Given an object $(λ,M)∈@CMC$, we have to show that $(λ,"Reedy(M))∈@CRC$.
That is, if $M$ is a small $(κ,λ)$-^{miniature model category}, we have to show that the small relative category $"Reedy(M)$ admits $λ$-small homotopy colimits.
Since we have a ^{Dwyer–Kan equivalence} $"Reedy(M)→"cof(M)$, it suffices to recall (^!{Cof is correct}) that $"cof(M)$ admits $λ$-small homotopy colimits.
\ppar
A left Quillen functor $M→N$ induces a left Quillen functor $M^Δ→N^Δ$ between the corresponding Reedy model categories of cosimplicial objects.
Therefore, it induces a relative functor $"Reedy(M)→"Reedy(N)$ between the corresponding relative categories of cofibrant objects.
If $(λ,M)→(μ,N)$ is a morphism, then the relative functor $"Reedy(M)→"Reedy(N)$
is weakly equivalent to the relative functor $"cof(M)→"cof(N)$, which preserves $λ$-small homotopy colimits (^!{preservation of relative homotopy limits}).
Thus, a morphism $g:(λ,M)→(μ,N)$ in $@CMC$ is sent to a morphism~$h$ in $@CRC$.
Furthermore, if $g$ is a weak equivalence, then so is~$h$ by the 2-out-of-3 property.
\ppar
Finally, the natural transformation $"ev_0(λ,C):(λ,"Reedy(C))→(λ,"Cof(C))$ is a weak equivalence:
its weak inverse is a natural transformation that sends $X∈"Cof(C)$ to the Reedy cofibrant resolution of the constant cosimplicial object on~$X$.

\section From combinatorial relative categories to combinatorial model categories

In this section we introduce and study constructions that allow us to pass from the relative category $@CRC$ to the relative category $@CMC$.
The primary source of difficulty is the fact that the regular cardinal $λ$ may increase in an uncontrolled fashion.
This prevents us from defining a relative functor $@CRC→@CMC$.
Instead, we provide an ad hoc construction for every small subcategory of $@CRC$.

\proclaim Definition.
A ^={simplicial set[|s]} is a simplicial object in the category $@Set$ of ^!{smaller sets}.
The category of simplicial sets is denoted by $@=sSet$.

\proclaim Definition.
Given a small relative category~$C$,
the model category $"=sPSh(C)$ of simplicial presheaves on~$C$ is defined as follows.
Its underlying category is the category of simplicial
objects in the strict free cocompletion of~$C$ (^!{strict cocompletion}).
By abuse of language, we refer to objects of $"sPSh(C)$ as {\it simplicial presheaves on~$C$}.
The universal property of strict free cocompletions
constructs an equivalence of categories $"sPSh(C)→@Cat(C^@op,@sSet)$,
which allows us to define a projective model structure on $"sPSh(C)$.
The model structure on $"sPSh(C)$ is defined as the left Bousfield localization of the projective model structure
at morphisms of simplicial presheaves that are representable by a weak equivalence in~$C$.

Under the equivalence of $"sPSh(C)$ with functors $C^@op→@sSet$,
the fibrant objects in $"sPSh(C)$ are precisely the relative functors $C^@op→@sSet_@=Kan$.

\proclaim Definition.
^^={sPSh on functors}
Given a relative functor $"F:C→D$ between small relative categories,
the left Quillen functor $$"sPSh("F):"sPSh(C)→"sPSh(D)$$ is induced by the construction of ^!{strict cocompletion}.
It is a simplicial left Quillen functor that restricts to $"F$ on representable presheaves.
This construction yields a (strict) relative functor $$"sPSh:@RelCat→@=CombModCat.$$

\proclaim Remark.
^^={abuse}
^!{sPSh on functors} contains a considerable abuse of notation:
the category $@CombModCat$ is supposed to have combinatorial model categories as objects,
which is not possible since combinatorial model categories of simplicial presheaves have a proper class of objects.
However, we only need the functor $"sPSh$ to construct the functor $"MInd$ (^!{MInd}),
itself used to construct the functor $"MInd_ν$ (^!{smaller MInd})
landing in ^{miniature model categories} (^!{miniature model category}), which do form a relative category.
Thus, the functor $"MInd_ν$ is well-defined and the abuse of notation is harmless.

We now introduce the small model category $"MInd(λ,C)$, which models the homotopy $λ$-ind-completion $\Ind^λ C$ of a small relative category~$C$.
The 1-categorical construction that we imitate here presents $\Ind^λ C$
by the category of functors $C^@op→@Set$ that preserve $λ$-small limits, provided that $C$ admits $λ$-small colimits.
The latter category can be encoded in turn as the reflective localization of the category of functors $C^@op→@Set$ at morphisms of the form $\colim_I "Y∘D → "Y(\colim_I D)$
for small diagrams $D:I→C$.
In the model-categorical setting, reflective localizations become left Bousfield localizations
and we use quasicategories to define the class of localizing morphisms to avoid developing the relevant machinery of homotopy colimits directly for relative categories.

\proclaim Definition.
^^={MInd}
Given an object $(λ,C)∈@CRC$,
the model category $"=MInd(λ,C)$
is defined as the left Bousfield localization of $"sPSh(C)$
at the set of maps of the form $η_D$ (constructed in the next paragraph)
for a set of representatives~$D$ of weak equivalence classes of
diagrams $D:I→"sPSh(C)$ of weakly representable presheaves,
where $I$ is a $λ$-small relative category.
Since $C$ is a small relative category, such representatives form a set.
The resulting left Bousfield localization is independent of the choices of $D$ and~$η_D$.
\ppar
The morphism $η_D$ is constructed as follows.
Consider the adjunction of quasicategories
$$\rN\frc "sPSh(C) \ltogets{11}{"L}{\rN\frc "Y} \rN\frc C,$$
where $"Y:C→"sPSh(C)$ is the Yoneda embedding functor and $"L$ is the left adjoint of $\rN\frc "Y$.
Suppose $I$ is a $λ$-small relative category and $D:I→"sPSh(C)$ is a relative functor.
\ppar
Consider the induced diagram of quasicategories $$\rN\frc(D):\rN\frc I→\rN \frc "sPSh(C).$$
The unit map~$α_D$ of the object $$\colim (\rN\frc(D))∈\rN \frc "sPSh(C)$$ has the form
$$α_D:\colim (\rN\frc(D)) → \rN\frc "Y("L(\colim \rN\frc(D))).$$
Denote by $η_D$ some morphism in $"sPSh(C)$ whose image in $\rN\frc "sPSh(C)$ is equivalent to~$α_D$.
This completes the construction of~$η_D$ and the definition of $"MInd(λ,C)$.
\ppar
Given a morphism $"F:(λ,C)→(μ,D)$ in $@CRC$,
the left Quillen functor $$"MInd("F):"MInd(λ,C)→"MInd(μ,D)$$
coincides with $"sPSh("F)$ as a functor.
In particular, $"MInd$ is itself a functor, keeping in mind ^!{abuse}.

\proclaim Definition.
^^={smaller MInd}
Given an object $(λ,C)∈@CRC$ and a regular cardinal $μ⊵λ$ such that $"MInd(λ,C)$ is a strongly $(κ,μ)$-combinatorial model category (^!{strongly combinatorial model category})
for some regular cardinal~$κ$,
the small model category $"MInd_μ(λ,C)$
is defined as the model category $\K_μ("MInd(λ,C))$ (Low [\Heart, Proposition~5.12]),
which is guaranteed to be small by ^!{smaller locally presentable}.
\ppar
If $"MInd_ν(λ,C)$ and $"MInd_ν(μ,D)$ are defined and $"MInd("F)$ is a left $ν$-Quillen functor (^!{left $λ$-Quillen functor}),
then we denote by
$$"MInd_ν("F):"MInd_ν(λ,C)→"MInd_ν(μ,D)$$
the functor $\K_ν("MInd("F))$.
\ppar
By ^!{miniature strongly CMC}, the functor $"MInd_ν$ is a relative functor
from a (nonfull) subcategory of the relative category $@CRC$
to the relative category $@CMC$
if we decorate the resulting objects and morphisms with $ν$ as the first component.

\proclaim Proposition.
^^={smaller MInd is miniature}
Given an object $(λ,C)∈@CRC$,
there are arbitrarily large regular cardinals~$μ⊵λ$ such that
the small model category $"MInd_μ(λ,C)$ is defined and is a $(κ,μ)$-^{miniature model category} for some regular cardinal~$κ$.
\ppar
Given a morphism $(λ,C)→(μ,D)$ in~$@CRC$,
there are arbitrarily large regular cardinals~$ν$ such that
the left Quillen functor $"MInd_ν(λ,C)→"MInd_ν(μ,D)$ is defined.

\proof Proof.
Apply ^!{miniature strongly CMC} and ^!{left Quillen functors are $λ$-Quillen}.

\proclaim Definition.
^^={Reedy MInd}
Suppose $ι:C→@CRC$ is an inclusion of a small full subcategory~$C$
and $ν$ is a regular cardinal such that $"MInd_ν$ is defined for all objects and morphisms of~$C$.
(Such a regular cardinal always exists by ^!{smaller MInd is miniature}.)
The natural transformation $$η:ι → "Reedy∘"MInd_ν∘ι$$ sends
an object $(κ,E)∈C$ to the morphism
$$(κ,E)→(ν,"Reedy("MInd_ν(κ,E)))$$
induced by the canonical functor $$E→"Reedy("MInd_ν(κ,E)),\qquad X↦(n↦Δ^n⊗"Y(X)).$$

\proclaim Proposition.
^^={Reedy MInd equivalence}
The natural transformation~$η$ of ^!{Reedy MInd} is a natural weak equivalence.

\proof Proof.
Compose the morphism
$$(κ,E)→(ν,"Reedy("MInd_ν(κ,E)))$$
with the weak equivalence
$$"ev_0:(ν,"Reedy("MInd_ν(κ,E)))→(ν,"MInd_ν(κ,E)).$$
It remains to show that the composition
$$(κ,E)→(ν,"MInd_ν(κ,E))$$
is a weak equivalence.
\ppar
By Cisinski [\HCHA, Remark~7.9.10],
the functor $\rN\frc$ applied to the projective model structure on simplicial presheaves on~$E$
yields a quasicategory equivalent to the quasicategory of presheaves on the nerve of~$E$.
By Cisinski [\HCHA, Proposition~7.11.4],
the quasicategory $\rN\frc "sPSh(E)$
is equivalent to the reflective localization of the quasicategory of presheaves on the nerve of~$E$
with respect to weak equivalences of~$E$.
The latter localization is itself equivalent to the quasicategory of presheaves on $\rN\frc E$.
Furthermore, by the same proposition,
the left Bousfield localization $"MInd(κ,E)$ (^!{MInd}) of $"sPSh(E)$
is equivalent to the reflective localization of presheaves on $\rN\frc E$
at morphisms constructed in ^!{MInd}.
The latter localization is itself equivalent to the category of presheaves on $\rN\frc E$
that (as functors from $(\rN \frc E)^@op$ to spaces) preserve $κ$-small limits.
This is precisely the $κ$-ind-completion of the quasicategory $\rN\frc E$,
which shows that $(κ,E)→(ν,"MInd_ν(κ,E))$ is a weak equivalence by definition of~$@CRC$.

\proclaim Definition.
^^={MInd Reedy}
Suppose $ι:C→@CMC$ is an inclusion of a small full subcategory~$C$
into the relative category $@CMC$ (^!{CMC})
and $ν$ is a regular cardinal such that $"MInd_ν$ (^!{smaller MInd}) is defined for all objects and morphisms of the diagram $"Reedy∘ι$ (^!{Reedy}).
The natural transformation
$$"Re:"MInd_ν∘"Reedy∘ι → ι$$
sends an object $(λ,M)∈@CMC$
to the morphism
$$(ν,"MInd_ν(λ,"Reedy(M)))→(ν,\Ind^λ_ν(M))$$
given by the left Quillen functor
$$"=Re:"MInd_ν(λ,"Reedy(M))→\Ind^λ_ν(M)$$
induced by the functor $$Δ^@op⨯"Reedy(M)→M$$ that sends $([n],X)↦X_n$.

\proclaim Remark.
^^={Reedy Cof difference}
The proof of ^!{MInd Reedy equivalence} explains why $"Re$ is a left Quillen functor.
The proof uses the specific properties of the functor $"Reedy$ (^!{Reedy})
and does not work for the weakly equivalent functor $"Cof$ (^!{Cof}).
Indeed, the functor $"Re$ was defined by means of the functor $([n],X)↦X_n$,
and the only obvious analogue of this construction for $"Cof$ sends $([n],X)↦X$.
The latter formula, however, prevents $"Re$ from being a left Quillen functor
because the image of a generating cofibration $(∂Δ^1→Δ^1)⊗X$ is the map $X⊔X→X$,
which is rarely a cofibration.

\proclaim Proposition.
^^={MInd Reedy equivalence}
The natural transformation
$$"Re:"MInd_ν∘"Reedy∘ι → ι$$
of ^!{MInd Reedy} is a natural weak equivalence.

\proof Proof.
To show that for any $(λ,M)∈@CMC$ the left adjoint functor
$$"Re:"MInd_ν(λ,"Reedy(M))→\Ind^λ_ν(M)$$
is a left Quillen equivalence, it suffices to show that
the left adjoint functor
$$"RE:"MInd(λ,"Reedy(M))→\Ind^λ(M)$$
(defined in the same way as~$"Re$) is a left $ν$-Quillen equivalence (^!{left $λ$-Quillen functor}),
after which we can pass to the subcategories of $ν$-presentable objects to recover~$"Re$.
\ppar
Given a model category~$M$,
consider the left adjoint functor $$"=RE:"sPSh("Reedy(M))→\Ind^λ M$$
that sends $([n],X)↦X_n$.
This functor is a left Quillen functor because the image of some generating projective cofibration $(∂Δ^n→Δ^n)⊗X$ is precisely the $n$th latching map of~$X$,
which is a cofibration by definition of a Reedy cofibrant cosimplicial object.
Likewise, the image of some generating projective acyclic cofibration $(Λ^n_k→Δ^n)⊗X$ is a weak equivalence.
Finally, a weak equivalence $X→X'$ of representable presheaves is sent to the morphism $X_0→X'_0$ in $\Ind^λ M$, which is a weak equivalence by definition of $"Reedy(M)$.
\ppar
Next, observe that the left Quillen functor $"RE$
factors through the localization $$"sPSh("Reedy(M))→"MInd(λ,"Reedy(M)).$$
Indeed, suppose $D:I→"sPSh("Reedy(M))$ is a $λ$-small diagram of weakly representable simplicial presheaves
and consider the morphism $η_D$ constructed in ^!{MInd}.
To show that the left derived functor of~$"RE$ sends $η_D$ to a weak equivalence in $\Ind^λ M$,
pass to the setting of quasicategories by restricting to cofibrant objects and applying the functor $\rN\frc$,
which yields the functor of quasicategories
$$\rN\frc("Cof("sPSh("Reedy(M))))→\rN\frc("Cof("MInd(λ,"Reedy(M)))).$$
By ^!{MInd}, the image of~$η_D$ in the quasicategory $\rN\frc "sPSh("Reedy(M))$
is equivalent to the unit map
$$α_D:\colim (\rN\frc(D)) → \rN\frc "Y("L(\colim \rN\frc(D))),$$
and the functor $\rN\frc("Cof∘"RE)$ is equivalent to~$"L$.
By the triangle identity for quasicategorical adjunctions,
the map $"L(α_D)$ is equivalent to the identity map on the object $"L(\colim \rN\frc(D))$ in the quasicategory $\rN\frc "Reedy(M)$,
which shows that the left derived functor of~$"RE$ sends the map $η_D$ to a weak equivalence in $\Ind^λ M$.
\ppar
The functor~$"RE$ is ^{homotopically essentially surjective} (^!{homotopically essentially surjective}).
Indeed, given any object $X∈M$, take the Reedy cofibrant resolution~$R$ of the constant cosimplicial object on~$X$.
Then $"RE("Y(R))=R_0∈M⊂\Ind^λ(M)$, so every object in $M⊂\Ind^λ(M)$ is weakly equivalent to an object in the image of the left derived functor of~$"RE$.
Since the latter image is closed under small $λ$-filtered homotopy colimits in $\Ind^λ(M)$, its closure under weak equivalences must coincide with $\Ind^λ(M)$.
\ppar
The right adjoint of $"RE$ is the functor
$$R:\Ind^λ(M)→"MInd(λ,"Reedy(M)),\qquad X↦(([n],R)↦M(R_n,X)).$$
The functor~$R$ preserves $λ$-filtered colimits, hence its right derived functor preserves $λ$-filtered homotopy colimits.
\ppar
The regular cardinal~$ν$ satisfies the conditions of Dugger [\Pres, Proposition~3.2],
so the functor~$"RE$ is a left Quillen equivalence
once we show that
the derived unit map of any object $P∈"MInd(λ,"Reedy(M))$ is a weak equivalence.
Since the left derived functor of $"RE$ and the right derived functor of $R$ preserve $λ$-filtered homotopy colimits,
it suffices to
establish the case when $P$ is a $λ$-small homotopy colimit of representable presheaves in $"MInd(λ,"Reedy(M))$.
By construction of $"MInd(λ,"Reedy(M))$, any such homotopy colimit is weakly equivalent to the representable presheaf of some $Q∈"Reedy(M)$.
Without loss of generality we can assume $Q$ to be (the representable presheaf of) a Reedy bifibrant cosimplicial object in~$M$.
Now $"RE("Y(Q))=Q_0$ is bifibrant in $\Ind^λ(M)$,
so the derived unit map of~$Q$ is simply the ordinary unit map of~$Q$.
Its codomain is $$R("RE("Y(Q)))=R(Q_0)=(([n],R)↦M(R_n,Q_0)).$$
Observe that the simplicial set $([n],R)↦M(R_n,Q_0)$ is weakly equivalent to the derived mapping simplicial set from~$R_0$ to~$Q_0$,
since $R$ is a cosimplicial resolution of~$R$.
Thus, the simplicial presheaf $R("RE("Y(Q)))$ is weakly equivalent to the representable presheaf of~$Q_0$, hence also to the representable presheaf of~$Q$.

\section Equivalence of combinatorial model categories and combinatorial relative categories

\proclaim Theorem.
^^={Reedy equivalence}
The relative functor $$"Reedy:@CMC→@CRC$$ (^!{Reedy}) is a ^{Dwyer–Kan equivalence} of relative categories.

\proof Proof.
The functor $"Reedy$ is ^{homotopically essentially surjective} by ^!{Reedy homotopically essentially surjective} and ^{homotopically fully faithful} by ^!{Reedy homotopically fully faithful},
so by ^!{criterion for relative equivalences} it is a ^{Dwyer–Kan equivalence} of relative categories.

Somewhat more generally, we have the following result.

\proclaim Theorem.
^^={Reedy generalized equivalence}
Suppose $Λ:C→@CMC$ is a relative functor
such that the construction of $"MInd_ν$ (^!{smaller MInd})
as well as the constructions of ^!{Reedy MInd} and ^!{MInd Reedy} lift through~$Λ$,
and ^!{Reedy MInd equivalence} and ^!{MInd Reedy equivalence} continue to hold for these lifts.
Then the relative functor $"Reedy∘Λ$ (and hence also $"Cof∘Λ$) is a ^{Dwyer–Kan equivalence} of relative categories.
In particular, the relative functor~$Λ$ itself is a ^{Dwyer–Kan equivalence} of relative categories.
\ppar
More generally, suppose $Λ:C→@CMC$ is a relative functor and $Σ:D→@CRC$ is a relative inclusion
such that the functor $"Reedy∘Λ$ factors through the image of~$Σ$,
and the construction of $"MInd_ν$ (^!{smaller MInd})
as well as the constructions of ^!{Reedy MInd} and ^!{MInd Reedy} lift through~$Λ$
once we restrict them to the image of~$Σ$,
and ^!{Reedy MInd equivalence} and ^!{MInd Reedy equivalence} continue to hold for these lifts.
Then the relative functor $"Reedy∘Λ:C→D$ is a ^{Dwyer–Kan equivalence} of relative categories.

\proof Proof.
^!{Reedy homotopically essentially surjective} and ^!{Reedy homotopically fully faithful} continue to hold in this generality,
since their proofs use precisely the indicated properties of $"MInd$ and the natural transformations of ^!{Reedy MInd} and ^!{MInd Reedy}.

\proclaim Proposition.
^^={main theorem simplicial}
^!{Reedy generalized equivalence} is applicable to the following choices of~$C$, constructed exactly like $@CMC$ (^!{CMC}), but with the indicated changes to objects and morphisms:
\li left proper miniature model categories and left Quillen functors;
\li miniature simplicial model categories and simplicial left Quillen functors;
\li left proper miniature simplicial model categories and simplicial left Quillen functors;
\endlist
Here a $(κ,λ)$-miniature simplicial model category is a $(κ,λ)$-^{miniature model category} (^!{miniature model category})
enriched over the cartesian model category of $λ$-small simplicial sets.

\proof Proof.
This is an immediate consequence of the construction
of $"MInd_ν$ as a left Bousfield localization of the category of simplicial presheaves on a small category.
We remark that the notions of left properness and simpliciality for ^{miniature model categories} match
the same notions for combinatorial model categories:
see Low [\Heart, Remark~5.17] for left proper model categories
and Low [\Heart, Remark~5.19] for simplicial model categories.

\proclaim Proposition.
^^={main theorem cartesian}
^!{Reedy generalized equivalence} is applicable to the following choices of~$C$ and~$D$, constructed exactly like $@CMC$ (^!{CMC}) and $@CRC$ (^!{CRC}),
but with the indicated changes to objects and morphisms:
\li For $C$, we take cartesian combinatorial model categories, which we can require to be left proper, or cartesian, or both.
\li For $D$, we take relative categories~$(λ,C)$ such that the category~$C$ admits finite products and the quasicategory $\rN\frc C$ is cartesian closed.
\endlist
Furthermore, the relative categories $C$ and $D$ are Dwyer–Kan equivalent to the full subcategory of $@PrL$ (^!{PrL}) on cartesian closed presentable quasicategories.

\proof Proof.
Given $(λ,C)∈@CRC$, the model category $"MInd(λ,C)$ is cartesian whenever $C$ has finite products (which ensures the pushout product axiom for cofibrations in $"sPSh(C)$)
and the morphisms used for the left Bousfield localization of $"sPSh(C)$ are closed under derived pushout products.
By Cisinski [\HCHA, Proposition~7.11.4] this is true whenever the quasicategory $\rN\frc C$
is a reflective localization of the quasicategory of presheaves on a small quasicategory
with respect to a set of morphisms that are closed under pushout products.
This is true for any cartesian closed quasicategory in $@PrL$.

\proclaim Proposition.
^^={Reedy homotopically essentially surjective}
The relative functor $$"Reedy:@CMC→@CRC$$ (^!{Reedy}) is a ^{homotopically essentially surjective} relative functor of relative categories.

\proof Proof.
Given an object $(λ,C)∈@CRC$,
^!{Reedy MInd equivalence} supplies (for a sufficiently large regular cardinal $μ⊵λ$)
a weak equivalence $$(λ,C)→(μ,"Reedy("MInd_μ(λ,C))),$$
which establishes the homotopy essential surjectivity of the relative functor $"Reedy$.

We are now ready to prove the main technical result of the whole article: ^!{Reedy homotopically fully faithful},
which shows that the relative functor $"Reedy:@CMC→@CRC$
is ^{homotopically fully faithful}.
A common way to establish such statements is to construct a weak inverse,
i.e., a relative functor of the form $@CRC→@CMC$.
As explained in ^!{Reedy MInd equivalence} and ^!{MInd Reedy equivalence}, we can construct such an inverse (namely, $"MInd_ν$ for some regular cardinal~$ν$)
for any small subcategory $C$ of~$@CRC$.
However, since we have no control over~$ν$ (i.e., the required choice of $ν$ does not seem to depend functorially on $(λ,C)∈@CRC$),
we cannot promote these choices to a single functor $@CRC→@CMC$.
This necessitates the more complicated proof of ^!{Reedy homotopically fully faithful}.

\proclaim Proposition.
^^={Reedy homotopically fully faithful}
The relative functor $$"Reedy:@CMC→@CRC$$ (^!{Reedy}) is a ^{homotopically fully faithful} relative functor of relative categories:
for any objects $(λ,C),(μ,D)∈@CMC$, the induced map $$\ham_{@CMC}((λ,C),(μ,D))→\ham_{@CRC}("Reedy(λ,C),"Reedy(μ,D))$$ is a simplicial weak equivalence.

\proof Proof.
We invoke a variant of the simplicial Whitehead theorem (^!{specialized simplicial Whitehead}).
Suppose we are given a commutative square
$$\sqcd{
"Sd^k ∂Δ^n&\mapright{α}&\ham_{@CMC}((λ,C),(μ,D))\cr
\mapdown{}&&\mapdown{\ham_{"Reedy}}\cr
"Sd^k Δ^n&\mapright{β}&\ham_{@CRC}("Reedy(λ,C),"Reedy(μ,D)),\cr
}$$
where $"Sd$ denotes the
barycentric subdivision functor.
Denote by~$Λ$ the simplicial subset of $Δ^2$ generated by the 1-simplices $0→2$ and $1→2$.
We construct maps
$$γ:"Sd^k Δ^n→\ham_{@CMC}((λ,C),(μ,D)),\qquad Γ:Δ^1⨯"Sd^k Δ^n→\ham_{@CRC}("Reedy(λ,C),"Reedy(μ,D)),$$
$$π:Λ⨯"Sd^k ∂Δ^n→\ham_{@CMC}((λ,C),(μ,D)),\qquad Π:Δ^2⨯"Sd^k ∂Δ^n→\ham_{@CRC}("Reedy(λ,C),"Reedy(μ,D))$$
such that
the map $Γ$ is a simplicial homotopy from $β$ to $\ham_{"Reedy}∘γ$,
the map $Π$ restricts to $\ham_{"Reedy}∘π$ on $Λ⨯"Sd^k ∂Δ^n$,
the map $π$ restricts to $α$ on $0⨯"Sd^k ∂Δ^n$,
the restrictions of $π$ to $1⨯"Sd^k ∂Δ^n$ and $γ$ to $"Sd^k ∂Δ^n$ coincide,
and the restrictions of $Π$ to $(0→1)⨯"Sd^k ∂Δ^n$ and $Γ$ to $Δ^1⨯"Sd^k ∂Δ^n$ coincide.
The maps $Γ$, $γ$, $Π$, and~$π$ are constructed in the remainder of the proof.
All conditions required for $Γ$, $γ$, $Π$, and~$π$ will be satisfied automatically by construction.
We refer the reader to ^!{specialized simplicial Whitehead figure} for a pictorial representation of the maps $Γ$, $γ$, $Π$, and~$π$.
\ppar
\step Reduction to a fixed zigzag type.
Recall (^!{hammock construction}) that for a relative category~$@C$ with objects $X,Y∈@C$, the simplicial set $\ham_@C(X,Y)$
is constructed as the colimit
$$\colim_{Z∈\zz} \N(@C^Z_{X,Y}),$$
where $Z$ runs over the category of zigzag types (Dwyer–Kan [\Calc, §4.1]), $\N$ denotes the nerve functor,
and $@C^Z_{X,Y}$ is the category of relative functors $Z→@C$
that map the leftmost and rightmost objects of~$Z$ to $X$ and~$Y$ respectively.
\ppar
By ^!{hammock construction hocolim}, the colimit over~$Z$ computes the homotopy colimit.
Thus, it suffices to show that for every zigzag type~$Z$,
the above square with the right map replaced by $$@CMC_{(λ,C),(μ,D)}^Z→@CRC_{"Reedy(λ,C),"Reedy(μ,D)}^Z$$
is a simplicial weak equivalence.
From now on, we work with a fixed zigzag type~$Z$.
\ppar
Now, maps of simplicial sets
$S→\N(@C^Z_{X,Y})$
can be identified with functors $π_{≤1}S→@C^Z_{X,Y}$,
where $π_{≤1}$ denotes the fundamental category functor.
The latter functors can themselves be identified with relative functors $Z⨯π_{≤1}S→@C$
that are constant functors valued in $X$ respectively $Y$ when restricted to the leftmost respectively rightmost object of~$Z$.
From now on, we interpret existing simplicial maps and construct new simplicial maps to $\ham$ in this form, as diagrams given by relative functors $Z⨯π_{≤1}S→@C$.
Since the value of such a diagram on the leftmost and rightmost vertex of~$Z$ is prescribed,
in the remainder of the proof we construct relative functors $Z⨯π_{≤1}S→@C$ as follows:
we pick some interior vertex of~$Z$,
construct a functor $π_{≤1}S→@C$,
establish naturality with respect to morphisms in~$Z$,
and verify the fact that left-pointing maps are sent to weak equivalences.
\ppar
\step Selection of the regular cardinal~$ν$.
We now define the regular cardinal~$ν$ that will be used in constructions of the maps $Γ$, $γ$, $Π$, and~$π$.
Apply the functor $"MInd$ (^!{MInd}) to all vertices and edges of the diagram~$β$.
This produces a commutative diagram of combinatorial model categories.
Choose a regular cardinal~$ν$ such that all vertices in this diagram are strongly $(κ,ν)$-combinatorial model categories for some $κ⊲ν$ (^!{strongly combinatorial model category})
and all edges in this diagram are left $ν$-Quillen functors (^!{left $λ$-Quillen functor}).
Since $"Sd^k Δ^n$ has only finitely many nondegenerate vertices and edges,
such a regular cardinal~$ν$ exists by ^!{miniature strongly CMC} and ^!{left Quillen functors are $λ$-Quillen}.
\ppar
\step Construction of the map~$γ$.
Apply the functor $"MInd_ν$ (^!{smaller MInd}) to the diagram~$β$.
The choice of~$ν$ guarantees that $"MInd_ν$ is defined for all objects and morphisms of~$β$.
The resulting model categories are $(κ,ν)$-^{miniature model categories} by ^!{smaller MInd is miniature},
so we can interpret the result as a map
$$δ:"Sd^k Δ^n→\ham_{@CMC}((ν,"MInd_ν(λ,"Reedy(C))),(ν,"MInd_ν(μ,"Reedy(D)))).$$
\ppar
Define $Z'=→←Z→←$, i.e., the zigzag type~$Z'$ is obtained from~$Z$ by attaching 4 additional morphisms as indicated.
From now on, we will be constructing simplicial maps of zigzag type~$Z'$.
Where necessary, existing maps of zigzag type~$Z$ are silently promoted to the zigzag type~$Z'$ by adding identity morphisms.
Now produce a map
$$γ:"Sd^k Δ^n→\ham_{@CMC}((λ,C),(μ,D))$$
by attaching to every zigzag in~$δ$ the weak equivalences
$$(λ,C)→(ν,\Ind^λ_ν(C))←(ν,"MInd_ν(λ,"Reedy(C))),\qquad
(ν,"MInd_ν(μ,"Reedy(D)))→(ν,\Ind^μ_ν(D))←(μ,D).$$
Here the left Quillen functor $"MInd_ν(λ,"Reedy(C))→\Ind^λ_ν(C)$ is a left Quillen equivalence by ^!{MInd Reedy equivalence}.
\ppar
\step Construction of the map~$Γ$.
The map
$$Γ:Δ^1⨯"Sd^k Δ^n→\ham_{@CRC}("Reedy(λ,C),"Reedy(μ,D))$$
is a simplicial homotopy from $β$ to $\ham_{"Reedy}∘γ$
constructed as a natural transformation of diagrams of zigzag type~$Z'$, i.e., a functor
$$Z⨯π_{≤1}("Sd^k Δ^n)→@CRC^{π_{≤1}(Δ^1)}.$$
\ppar
First, promote $β$ to the zigzag type~$Z'$ by
precomposing with the relative functor $Z'→Z$ that collapses the outer two vertices on each side.
This amounts to attaching to every zigzag in~$β$ the identity maps
$$(λ,"Reedy(C))→(λ,"Reedy(C))←(λ,"Reedy(C)),\qquad (μ,"Reedy(D))→(μ,"Reedy(D))←(μ,"Reedy(D)),$$
ensuring that both $β$ and $\ham_{"Reedy}∘γ$ have the same zigzag type~$Z'$.
\ppar
Now we construct $Γ$ as a natural weak equivalence
from the diagram of~$β$ to the diagram of~$\ham_{"Reedy}∘γ$.
Following the tactic outlined in the paragraph on reduction to a fixed zigzag type,
we work with a fixed interior vertex $z∈Z'$ and construct a natural transformation of functors $π_{≤1}("Sd^k Δ^n)→@CRC$.
\ppar
If the vertex~$z$ belongs to~$Z⊂Z'$,
the value of~$Γ$ on some object $W∈π_{≤1}("Sd^k Δ^n)$
with $β(W)=(κ,E)∈@CRC$ is given by
the weak equivalence (^!{Reedy MInd equivalence}) in $@CRC$
$$Γ_{κ,E}:(κ,E)→(ν,"Reedy("MInd_ν(κ,E)))$$
whose underlying relative functor $$E→"Reedy("MInd_ν(κ,E))$$
sends an object $X∈E$ to the Reedy cofibrant cosimplicial diagram $n↦Δ^n⊗"Y(X)$.
\ppar
If the vertex~$z$ does not belong to~$Z⊂Z'$, then it is one of the two interior vertices added to the zigzag~$Z$.
Suppose $z$ is adjacent to the leftmost vertex of~$Z'$ (corresponding to $(λ,C)$); the other case (corresponding to $(μ,D)$) is treated symmetrically.
The resulting morphism does not depend on the choice of $W∈π_{≤1}("Sd^k Δ^n)$
and is given by the weak equivalence $$(λ,"Reedy(C))→(ν,"Reedy(\Ind^λ_ν(C)))$$
induced by the relative functor $$"Reedy(C)→"Reedy(\Ind^λ_ν(C))$$
obtained by applying $"Reedy$ to the canonical inclusion $$C→\Ind^λ_ν(C).$$
This completes the construction of~$Γ$.
\ppar
\step Construction of the maps $π$ and $Π$.
Next, we construct the maps
$$π:Λ⨯"Sd^k ∂Δ^n→\ham_{@CMC}((λ,C),(μ,D)),\qquad Π:Δ^2⨯"Sd^k ∂Δ^n→\ham_{@CRC}("Reedy(λ,C),"Reedy(μ,D))$$
using similar techniques.
As before, fix some interior vertex $z∈Z'$ and construct functors
$$π_{≤1}("Sd^k ∂Δ^n)→@CMC^{π_{≤1}Λ},\qquad π_{≤1}("Sd^k ∂Δ^n)→@CRC^{π_{≤1}Δ^2}.$$
If the vertex~$z$ belongs to~$Z⊂Z'$,
the value of~$π$ on some object $W∈π_{≤1}("Sd^k ∂Δ^n)$
with $α(W)=(κ,M)$ is given by the following object in $@CMC^{π_{≤1}Λ}$:
$$\diagram[3em]A:0,0=${(κ,M)}$,B:0,-2=${(ν,"MInd_ν(κ,"Reedy(M)))}${,},C:-3,-1=${(ν,\Ind^ κ_ν(M))\strut}$,BC.urt@0.45=$"ev$,AC.ulft@0.4=$ι$;$$
where the map~$ι$ is the canonical inclusion and the map $"ev$ is defined on representables via the formula $"ev(Δ^n⊗R)=R_n$, where $R∈"Reedy(M)$.
\ppar
Likewise, the map $Π$ is
given by the following object in $@CRC^{π_{≤1}Δ^2}$:
$$\diagram[4em]A:0,0=${(κ,"Reedy(M))}$,B:0,-2={${(ν,"Reedy("MInd_ν(κ,"Reedy(M))))}$,},C:-3,-1=${(ν,"Reedy(\Ind^ κ_ν(M)))\strut}$,AB.rt=$Γ_{κ,"Reedy(M)}$,BC.llft@0.5=$"Reedy("ev)$,AC.ulft@0.3=$"Reedy(ι)$;$$
where the map $Γ_{κ,"Reedy(M)}$ was defined in the previous part of the proof: it sends $R∈"Reedy(M)$ to the Reedy cofibrant object $n↦Δ^n⊗"Y(R)$.
\ppar
If the vertex~$z$ does not belong to~$Z⊂Z'$, then it is one of the two interior vertices added to the zigzag~$Z$.
Suppose $z$ is adjacent to the leftmost vertex of~$Z'$ (corresponding to $(λ,C)$); the other case (corresponding to $(μ,D)$) is treated symmetrically.
The resulting object in $@CMC^{π_{≤1}Λ}$ does not depend on the choice of $W∈π_{≤1}("Sd^k ∂Δ^n)$
and is given by the following diagram:
$$\diagram[2.5em]A:0,0=${(λ,C)}$,B:0,-2=${(ν,\Ind^ λ_ν(C))}${,},C:-4,-1=${(ν,\Ind^ λ_ν(C))}$,BC.llft@0.45=$\id$,AC.ulft@0.4=$ι$;$$
where $ι$ denotes the canonical inclusion.
Likewise, the map $Π$ is
given by the following object in $@CRC^{π_{≤1}Δ^2}$:
$$\diagram[3.5em]A:0,0=${(λ,"Reedy(C))}$,B:0,-2={${(ν,"Reedy(\Ind^ λ_ν(C)))}$,},C:-4,-1=${(ν,"Reedy(\Ind^ λ_ν(C)))\strut}$,AB.rt=$Γ_{λ,"Reedy(C)}$,BC.llft@0.4=$"Reedy(\id)$,AC.ulft@0.3=$"Reedy(ι)$;$$
where the map $Γ_{λ,"Reedy(C)}$ was defined in the previous part of the proof: evaluate $"Reedy$ on the canonical inclusion $C→\Ind^λ_ν(C)$.

\proclaim Theorem.
^^={Cof equivalence}
The relative functor $$"Cof:@CMC→@CRC$$ (^!{Cof}) is a ^{Dwyer–Kan equivalence} of relative categories.

\proof Proof.
The relative functor $"Cof$ is weakly equivalent to the relative functor $"Reedy$ (^!{Reedy}) via the natural weak equivalence $"ev_0:"Reedy→"Cof$ of ^!{Reedy}.
By ^!{Reedy equivalence}, $"Reedy$ is a ^{Dwyer–Kan equivalence}, hence so is $"Cof$.

The following proposition is not used anywhere else in the article.
It shows that the more straightforward way to define a ^{Dwyer–Kan equivalence} $@CMC'_U→@CRC'_U$
is weakly equivalent to the functor $"Cof$
under the ^{Dwyer–Kan equivalences} $@CMC_U→@CMC'_U$ (^!{smaller universe CMCU}) and $@CRC_U→@CRC'_U$ (^!{smaller universe CRCU}).

\proclaim Proposition.
^^={universe Cof RInd}
Suppose $U$ is a strongly inaccessible cardinal.
Consider the functor $"Cof_U:@CMC'_U→@CRC'_U$
that sends an object $M∈@CMC'_U$ to the relative category of cofibrant objects in~$M$
and a left Quillen functor $M→N$ in $@CMC'_U$ to the induced functor between the categories of cofibrant objects.
The functors $$"Cof_U ∘ \Ind_U: @CMC_U→@CRC'_U$$
and $$"RInd_U ∘ "Cof: @CMC_U → @CRC'_U$$
are naturally weakly equivalent.

\proof Proof.
The natural weak equivalence is given by the morphism
$$"Cof_U(\Ind_U(λ,M)) → "RInd_U("Cof(λ,M))$$
that sends a cofibrant object $A∈\Ind_U(λ,M)$ to $n↦Δ^n⊗(B↦\ham_M(B,A))$.

\section Equivalence of combinatorial relative categories and presentable quasicategories

\proclaim Definition.
^^={derived relative nerve}
The relative functor
$$\RN:@CRC→@PrL$$
between the relative categories $@CRC$ (^!{CRC}) and $@PrL$ (^!{PrL})
is defined as follows.
An object $(λ,C)$ is sent to $(λ,\rN \frc C)$.
A morphism $(λ,C)→(μ,D)$ given by a relative functor $"F:C→D$ is sent to the morphism $$(λ,\rN \frc C)→(μ,\rN \frc D)$$ given by the functor $$\rN \frc "F:\rN \frc C→\rN \frc D.$$
The relative functors $\rN$ and $\frc$ are introduced in ^!{relative nerve notation}.

\proclaim Proposition.
^!{derived relative nerve} is correct.

\proof Proof.
If $(λ,C)∈@CRC$, then the small quasicategory $\rN\frc C$ admits $λ$-small colimits by ^!{relative homotopy limit}.
Likewise, if $(λ,C)→(μ,D)$ is a morphism in~$@CRC$, the functor $\rN\frc C→\rN\frc D$ preserves $λ$-small colimits by ^!{preservation of relative homotopy limits}.
\ppar
We now show that $\RN$ preserves weak equivalences
by establishing this claim separately for each generating class.
If $(λ,C)→(μ,D)$ is a weak equivalence such that $λ=μ$ and $C→D$ is a ^{Dwyer–Kan equivalence},
then $\frc C → \frc D$ is a ^{Dwyer–Kan equivalence} between fibrant objects,
and $\rN \frc C → \rN \frc D$ is an equivalence of quasicategories because $\rN$ is a right Quillen functor.
For weak equivalences $(λ,C)→(μ,D)$ of the second generating class (i.e., involving ind-completions),
applying $\RN$ produces a weak equivalence in $@PrL$ by definition of~$@CRC$ (^!{CRC}),
since we defined the second generating class there as the preimage of corresponding weak equivalences in $@PrL$.

\proclaim Definition.
^^={relative realization}
The relative functor
$$\RK:@PrL→@CRC$$
between the relative categories $@PrL$ (^!{PrL}) and $@CRC$ (^!{CRC}) is defined as follows.
An object $(λ,C)$ is sent to $(λ,\rK C)$, where $\rK$ is the functor from ^!{relative nerve notation}.
A morphism $(λ,C)→(μ,D)$ given by a map of simplicial sets $"F:C→D$ is sent to the morphism $(λ,\rK C)→(μ,\rK D)$ given by the functor $\rK C→\rK D$.

\proclaim Proposition.
^!{relative realization} is correct.

\proof Proof.
Suppose $(λ,C)∈@PrL$.
Since the functor $C→\rN\frc \rK C$ is an equivalence,
by ^!{relative homotopy limit} the small relative category $\rK C$ admits $λ$-small homotopy colimits.
\ppar
Suppose $(λ,C)→(μ,D)$ is a morphism in~$@PrL$.
Since the morphism $\rN\frc\rK C → \rN\frc\rK D$ is weakly equivalent to $C→D$,
by ^!{preservation of relative homotopy limits} the relative functor $\rK C→\rK D$ preserves $λ$-small homotopy colimits.
\ppar
We show that $\RK$ preserves weak equivalences by establishing this claim separately for each generating class.
If $(λ,C)→(μ,D)$ is a weak equivalence such that $λ=μ$ and $C→D$ is an equivalence of quasicategories,
then the relative functor $\rK C → \rK D$ is a ^{Dwyer–Kan equivalence} because $\rK$ is a left Quillen functor and all simplicial sets are cofibrant in the Joyal model structure.
If $(λ,C)→(μ,D)$ is a weak equivalence such that $C→D$ exhibits $D$ as the quasicategory of $μ$-presentable objects in the $λ$-ind-completion of~$C$,
the morphism $(λ,\rK C) → (μ,\rK D)$ is a weak equivalence in $@CRC$
if its image under $\rN\frc$ is a weak equivalence in $@PrL$.
The resulting morphism $(λ,\rN\frc \rK C)→(μ,\rN\frc \rK D)$
is weakly equivalent to the original morphism $(λ,C)→(μ,D)$
via the derived unit map, which completes the proof.

\proclaim Definition.
^^={relative nerve relative realization}
The natural transformation
$$η: \id_{@PrL} → \RN ∘ \RK$$
from the identity functor on the relative category $@PrL$ (^!{PrL})
to the composition of relative functors $\RN$ (^!{derived relative nerve}) and $\RK$ (^!{relative realization})
is constructed as follows.
Given $(λ,C)∈@PrL$, we send it to the map
$$(λ,C)→(λ,\rN \frc \rK C)$$
given by composing the unit map $C → \rN \rK C$ with the map $\rN \rK C → \rN \frc \rK C$.
The relative functors $\rN$, $\rK$, and $\frc$ are introduced in ^!{relative nerve notation}.

\proclaim Proposition.
^^={relative nerve relative realization equivalence}
^!{relative nerve relative realization} is correct and the natural transformation~$η$ is a natural weak equivalence.

\proof Proof.
Suppose $(λ,C)→(μ,D)$ is a morphism in $@PrL$.
We must show that the square
$$\sqcd{
(λ,C)&\mapright{}&(μ,D)\cr
\mapdown{}&&\mapdown{}\cr
(λ,\rN \frc \rK C)&\mapright{}&(μ,\rN \frc \rK D)\cr
}$$
commutes, which follows from the commutativity of the following diagram:
$$\cd{
C&\mapright{}&D\cr
\mapdown{}&&\mapdown{}\cr
\rN \rK C&\mapright{}&\rN \rK D\cr
\mapdown{}&&\mapdown{}\cr
\rN \frc \rK C&\mapright{}&\rN \frc \rK D.\cr
}$$
The top square commutes because the unit is a natural transformation.
The bottom square commutes because $\frc$ is a functor and the fibrant replacement map $\id→\frc$ is a natural transformation.
\ppar
Finally, $η$ is a weak equivalence because $\rK$ and $\rN$ form a Quillen equivalence,
so the derived unit map of $\rK⊣\rN$ is a weak equivalence.

\proclaim Definition.
^^={fibrant CRC}
The relative endofunctor
$$\FRC: @CRC → @CRC$$
on the relative category $@CRC$ (^!{CRC})
is constructed as follows.
An object $(λ,C)∈@CRC$ is sent to $(λ,\frc C)$, where $\frc$ is the functor from ^!{relative nerve notation}.
A morphism $(λ,C)→(μ,D)$ given by a relative functor $"F:C→D$
is sent to the morphism $(λ,\frc C)→(μ,\frc D)$ given by the relative functor $\frc "F: \frc C→\frc D$.

\proclaim Definition.
^^={relative realization relative nerve}
The zigzag~$ε$ of natural transformations
$$\RK ∘ \RN → \FRC ← \id_{@CRC}$$
between functors $\RK∘\RN$ (^!{relative realization}, ^!{derived relative nerve}), $\FRC$ (^!{fibrant CRC}), and $\id_{@CRC}$
is constructed as follows.
Given $(λ,C)∈@CRC$, we send it to the zigzag
$$(λ,\rK \rN \frc C) → (λ, \frc C) ← (λ,C),$$
where the first map is the counit of $\frc C$ and the second map is the fibrant replacement map.

\proclaim Proposition.
^^={relative realization relative nerve equivalence}
^!{relative realization relative nerve} is correct and the zigzag of natural transformations~$ε$ is a zigzag of natural weak equivalences.

\proof Proof.
The naturality of the first transformation
follows from the naturality of counit maps
and the naturality of the second transformation
follows from the naturality of the fibrant replacement map $\id → \frc$.
The counit map $\rK \rN \frc C → \frc C$ is the derived counit map of a Quillen equivalence,
hence is a weak equivalence.
The fibrant replacement map is a weak equivalence by definition.

\proclaim Theorem.
^^={derived relative nerve equivalence}
The functor $\RN:@CRC→@PrL$ (^!{derived relative nerve}) is a ^{Dwyer–Kan equivalence}.

\proof Proof.
Combine ^!{relative nerve relative realization equivalence} and ^!{relative realization relative nerve equivalence}.

The following proposition is not used anywhere else in the article.
It shows that the more straightforward way to define a ^{Dwyer–Kan equivalence} $\RN_U:@CRC'_U→@PrL'_U$
is weakly equivalent to the functor $\RN$
under the ^{Dwyer–Kan equivalences} $"RInd_U:@CRC_U→@CRC'_U$ (^!{smaller universe CRCU}) and $"QInd_U:@PrL_U→@PrL'_U$ (^!{smaller universe PrLU}).

\proclaim Proposition.
^^={smaller CRCU universe PrLU}
Suppose $U$ is a strongly inaccessible cardinal.
Consider the functor $$\RN_U=\rN\frc:@CRC'_U→@PrL'_U,$$
where $\rN$ and $\frc$ are as in ^!{relative nerve notation}.
There is a zigzag of natural weak equivalences connecting
the functors $$\RN_U ∘ "RInd_U: @CRC_U→@PrL'_U$$
(^!{smaller RInd})
and $$"QInd_U ∘ \RN: @CRC_U → @PrL'_U,$$
with $"QInd_U$ as in ^!{smaller QInd} and $\RN$ as in ^!{derived relative nerve}, restricted to $@CRC_U$.

\proof Proof.
The natural weak equivalence that we need has the form
$$\RN_U("RInd_U(λ,C)) → "QInd_U(\RN(λ,C)).$$
Unfolding the definitions, we need a natural weak equivalence
$$\rN\frc("Reedy("MInd_U(λ,C))) → \rN\frc("Reedy("MInd_U(λ,\rK\rN\frc C))).$$
Such a natural weak equivalence is induced by the zigzag $\rK\rN\frc C→\frc C←C$ of ^!{relative realization relative nerve equivalence}.

\unsection References


\yearkeytrue

\refs

\bib\CHi
J.~H.~C.~Whitehead.
Combinatorial homotopy.  I.
Bulletin of the American Mathematical Society 55:3 (\y{1949}), 213–246.
\doi:10.1090/s0002-9904-1949-09175-9.

\bib\Kan
Daniel~M.~Kan.                                                                                                                                                                On~c.s.s.~categories.
Boletín de la Sociedad Matemática Mexicana 2 (\y{1957}), 82–94.
\https://dmitripavlov.org/scans/kan-on-css-categories.pdf.

\bib\HoAlg
Daniel~G.~Quillen.
Homotopical algebra.
Lecture Notes in Mathematics 43 (\y{1967}), Springer.
\doi:10.1007/bfb0097438.

\bib\LPK
Peter Gabriel, Friedrich Ulmer.
Lokal präsentierbare Kategorien.
Lecture Notes in Mathematics 221 (\y{1971}).
\doi:10.1007/BFb0059396.

\bib\SimpLoc
William~G.~Dwyer, Daniel~M.~Kan.
Simplicial localizations of categories.
Journal of Pure and Applied Algebra 17:3 (\y{1980}), 267--284.
\doi:10.1016/0022-4049(80)90049-3.

\bib\Calc
William~G.~Dwyer, Daniel~M.~Kan.
Calculating simplicial localizations.
Journal of Pure and Applied Algebra 18:1 (\y{1980}), 17--35.
\doi:10.1016/0022-4049(80)90113-9.

\bib\Func
William~G.~Dwyer, Daniel~M.~Kan.
Function complexes in homotopical algebra.
Topology 19:4 (\y{1980}), 427--440.
\doi:10.1016/0040-9383(80)90025-7.

\bib\AccCat
Michael Makkai, Robert Paré.
Accessible categories: the foundations of categorical model theory.
Contemporary Mathematics 104 (\y{1989}).
\doi:10.1090/conm/104.

\bib\LPAC
Jiří Adámek, Jiří Rosický.
Locally presentable and accessible categories.
London Mathematical Society Lecture Note Series 189 (\y{1994}), Cambridge University Press.
\doi:10.1017/cbo9780511600579.

\bib\Replacing[1998]
Daniel Dugger.
Replacing model categories with simplicial ones.
Transactions of the American Mathematical Society 353:12 (2001), 5003-5027.
\doi:10.1090/s0002-9947-01-02661-7.

\bib\SHT
Paul~G.~Goerss, John~F.~Jardine.
Simplicial homotopy theory.
Progress in Mathematics 174 (\y{1999}).
\doi:10.1007/978-3-0346-0189-4.

\bib\MC
Mark Hovey.
Model categories.
Mathematical Surveys and Monographs 63 (\y{1999}).
\doi:10.1090/surv/063.

\bib\Hovey
Mark Hovey.
Algebraic Topology Problem List.
Model categories.
March 6, \y{1999}.
{\font\small=cmr7 \small
\https://web.archive.org/web/20000830081819/http://claude.math.wesleyan.edu/~mhovey/problems/model.html.
}

\bib\SimStruc[2000]
Charles Rezk, Stefan Schwede, Brooke Shipley.
Simplicial structures on model categories and functors.
American Journal of Mathematics 123:3 (2001), 551--575.
\arXiv:math/0101162v1, \doi:10.1353/ajm.2001.0019.

\bib\Pres[2000]
Daniel Dugger.
Combinatorial model categories have presentations.
Advances in Mathematics 164:1 (2001), 177--201.
\arXiv:math/0007068v1, \doi:10.1006/aima.2001.2015.

\bib\SHM
Tibor Beke.
Sheafifiable homotopy model categories.
Mathematical Proceedings of the Cambridge Philosophical Society 129:3 (\y{2000}), 447--475.
\arXiv:math/0102087v1, \doi:10.1017/s0305004100004722.

\bib\QCKC
André Joyal.
Quasi-categories and Kan complexes.
Journal of Pure and Applied Algebra 175:1-3 (\y{2002}), 207–222.
\doi:10.1016/s0022-4049(02)00135-4.

\bib\MCL
Philip~S.~Hirschhorn.
Model categories and their localizations.
Mathematical Surveys and Monographs 99 (\y{2003}).
\doi:10.1090/surv/099.

\bib\DHKS
William~G.~Dwyer, Philip~S.~Hirschhorn, Daniel~M.~Kan, Jeffrey~H.~Smith.
Homotopy Limit Functors on Model Categories and Homotopical Categories.
Mathematical Surveys and Monographs 113 (\y{2004}).
\doi:10.1090/surv/113.

\bib\WESP[2002]
Daniel Dugger, Daniel~C.~Isaksen.
Weak equivalences of simplicial presheaves.
Contemporary Mathematics 346 (2004), 97--113.
\arXiv:math/0205025v1, \doi:10.1090/conm/346/06292.

\bib\Derivators[2006]
Olivier Renaudin.
Plongement de certaines théories homotopiques de Quillen dans les dérivateurs.
Journal of Pure and Applied Algebra 213:10 (2009), 1916--1935.
\arXiv:math/0603339v2, \doi:10.1016/j.jpaa.2009.02.014.

\bib\ThQ
Georges Maltsiniotis.
Le théorème de Quillen, d'adjonction des foncteurs dérivés, revisité.
Comptes Rendus Mathematique 344:9 (\y{2007}), 549--552.
\doi:10.1016/j.crma.2007.03.011.

\bib\LR[2007]
Clark Barwick.
On left and right model categories and left and right Bousfield localizations.
Homology, Homotopy and Applications 12:2 (2010), 245--320.
\arXiv:0708.2067v2, \doi:10.4310/hha.2010.v12.n2.a9.

\bib\RelCat[2010]
Clark Barwick, Daniel~M.~Kan.
Relative categories: Another model for the homotopy theory of homotopy theories.
Indagationes Mathematicae 23:1--2 (2012), 42--68.
\arXiv:1011.1691v2, \doi:10.1016/j.indag.2011.10.002.

\bib\ClassComb[2011]
Boris Chorny, Jiří Rosický.
Class-combinatorial model categories.
Homology, Homotopy and Applications 14:1 (2012), 263--280.
\arXiv:1110.4252v1, \doi:10.4310/hha.2012.v14.n1.a13.

\bib\TwoCat[2013]
Emily Riehl, Dominic Verity.
The 2-category theory of quasi-categories.
\arXiv:1306.5144v4.

\bib\Hobi[2013]
Zhen Lin Low.
The homotopy bicategory of $(∞,1)$-categories.
\arXiv:1310.0381v2.

\bib\DKLR[2013]
Vladimir Hinich.
Dwyer-Kan localization revisited.
Homology, Homotopy and Applications 18:1 (2016), 27–48.
\arXiv:1311.4128v4, \doi:10.4310/hha.2016.v18.n1.a3.

\bib\RectEC[2013]
Rune Haugseng.
Rectification of enriched $∞$-categories.
Algebraic \& Geometric Topology 15:4 (2015), 1931--1982.
\arXiv:1312.3881v4, \doi:10.2140/agt.2015.15.1931.

\bib\Heart[2014]
Zhen Lin Low.
The heart of a combinatorial model category.
Theory and Applications of Categories 31:2 (2016), 31–62.
\arXiv:1402.6659v4.

\bib\HTCQ[2014]
Karol Szumiło.
Homotopy theory of cocomplete quasicategories.
Algebraic \& Geometric Topology 17:2 (2017), 765–791.
\arXiv:1411.0303v1, \doi:10.2140/agt.2017.17.765.

\bib\QAdj[2015]
Aaron Mazel-Gee.
Quillen adjunctions induce adjunctions of quasicategories.
New York Journal of Mathematics 22 (2016), 57–93.
\arXiv:1501.03146v1.

\bib\NS[2015]
Thomas Nikolaus, Steffen Sagave.
Presentably symmetric monoidal $∞$-categories are represented by symmetric monoidal model categories.
Algebraic \& Geometric Topology 17:5 (2017), 3189--3212.
\arXiv:1506.01475v3, \doi:10.2140/agt.2017.17.3189.

\bib\Preder[2016]
Kevin Arlin.
On the $∞$-categorical Whitehead theorem and the embedding of quasicategories in prederivators.
Homology, Homotopy and Applications 22:1 (2020), 117--139.
\arXiv:1612.06980v4, \doi:10.4310/hha.2020.v22.n1.a8.

\bib\HTT
Jacob Lurie.
Higher Topos Theory.
April 9, \y{2017}.
\https://www.math.ias.edu/~lurie/papers/HTT.pdf.

\bib\HCHA
Denis-Charles Cisinski.
Higher Categories and Homotopical Algebra.
Cambridge Studies in Advanced Mathematics 180 (\y{2019}).
\doi:10.1017/9781108588737.

\bib\Clusters[2021]
Erwan Beurier, Dominique Pastor, René Guitart.
Presentations of clusters and strict free-cocompletions.
Theory and Applications of Categories 36:17 (2021), 492–513.
\http://www.tac.mta.ca/tac/volumes/36/17/36-17.pdf.

\bib\Kerodon[2021]
Jacob Lurie.
Kerodon.
\https://kerodon.net/.
